\documentclass[11pt,leqno]{article}
\usepackage{amssymb}
\usepackage[mathscr]{eucal}
\usepackage{amsmath,amssymb,latexsym,theorem,bbm}
\usepackage{color, url}
\usepackage{appendix}

\newcommand{\CC}{\mathsf{C}}
\newcommand{\DD}{\mathsf{D}}

\newcommand{\NN}{\mathbb{N}}

\newcommand{\RR}{\mathbb{R}}

\newcommand{\ZZ}{\mathbb{Z}}

\newcommand{\cA}{{\mathcal A}}
\newcommand{\cB}{{\mathcal B}}

\newcommand{\cC}{{\mathcal C}}
\newcommand{\cD}{{\mathcal D}}

\newcommand{\cG}{{\mathcal G}}
\newcommand{\cF}{{\mathcal F}}

\newcommand{\cL}{{\mathcal L}}

\newcommand{\cN}{{\mathcal N}}

\newcommand{\cX}{{\mathcal X}}
\newcommand{\cY}{{\mathcal Y}}

\newcommand{\cW}{{\mathcal W}}

\newcommand{\dd}{\mathrm{d}}
\newcommand{\ee}{\mathrm{e}}

\newcommand{\INARtwo}{\textup{INAR(2)}}

\DeclareMathOperator*{\argmin}{arg\,min}

\newcommand{\EE}{\operatorname{\mathbb{E}}}
\newcommand{\PP}{\operatorname{\mathbb{P}}}

\newcommand{\tX}{\widetilde{X}}
\newcommand{\tY}{\widetilde{Y}}

\renewcommand{\mid}{\,|\,}

\renewcommand{\leq}{\leqslant}
\renewcommand{\geq}{\geqslant}

\newcommand{\stoch}{\stackrel{\PP}{\longrightarrow}}
\newcommand{\distr}{\stackrel{\cL}{\longrightarrow}}
\newcommand{\distre}{\stackrel{\cL}{=}}

\newcommand{\as}{\stackrel{{\mathrm{a.s.}}}{\longrightarrow}}
\newcommand{\ase}{\stackrel{{\mathrm{a.s.}}}{=}}

\newcommand{\ns}{{\lfloor ns\rfloor}}

\newcommand{\proofend}{\hfill\mbox{$\Box$}}

\numberwithin{equation}{section}

\theoremstyle{change} \theorembodyfont{\em}
\newtheorem{Lem}{Lemma.}[section]
\newtheorem{Thm}[Lem]{Theorem.}
\newtheorem{Pro}[Lem]{Proposition.}
\newtheorem{Cor}[Lem]{Corollary.}
\newtheorem{Def}[Lem]{Definition.}

\theorembodyfont{\rm}
\newtheorem{Rem}[Lem]{Remark.}

\begin{document}

\begin{center}
 {\bfseries\Large
  On parameter estimation for critical affine processes} \\[5mm]
 {\sc\large
  M\'aty\'as $\text{Barczy}^{*}$,
  \ Leif $\text{D\"oring}$,
  \ Zenghu $\text{Li}$,
  \ Gyula $\text{Pap}$}
\end{center}


\renewcommand{\thefootnote}{}
\footnote{$*$ Corresponding author}
\vspace*{0.2cm}
\footnote{M. Barczy, Z. Li and G. Pap have been supported by the Hungarian
 Chinese Intergovernmental S \& T Cooperation Programme for 2011-2013
 under Grant No.\ 10-1-2011-0079.
M. Barczy and G. Pap have been partially supported by the Hungarian Scientific
 Research Fund under Grant No.\ OTKA T-079128.
L. D\"oring has been supported by the Foundation Science Mat\'ematiques de Paris. 
Z. Li has been partially supported by NSFC under Grant No.\ 11131003 and 973
 Program under Grant No.\ 2011CB808001.}

\vspace*{-10mm}

\begin{abstract}
First we provide a simple set of sufficient conditions for the weak convergence
 of scaled affine processes with state space $\mathbb{R}_+ \times \mathbb{R}^d$.
 We specialize our result to one-dimensional continuous state branching
 processes with immigration.
As an application, we study the asymptotic behavior of least squares estimators
 of some parameters of a two-dimensional critical affine diffusion process.
\end{abstract}

\section{Introduction}

In recent years quickly growing interest in pricing of credit-risky securities (e.g., defaultable bonds)
 has been seen in the mathematical finance literature.
One of the basic models (for applications see for instance Chen and Joslin \cite{CheJos}) is the following two-dimensional
 affine diffusion process:
 \begin{align}\label{DD}
  \begin{cases}
   \dd Y_t=(a-bY_t)\,\dd t+\sqrt{Y_t}\,\dd W_t,\\
   \dd X_t= (m-\theta X_t)\,\dd t+\sqrt{Y_t}\,\dd B_t,
  \end{cases} \quad t\geq 0,
 \end{align}
 where \ $a$, $b$, $\theta$ \ and \ $m$ \ are real parameters such that \ $a>0$ \ and \ $B$ \ and \ $W$ \ are
 independent standard Wiener processes.
Note that \ $Y$ \ is a Cox-Ingersol-Ross (CIR) process.
For practical use, it is important to estimate the appearing parameters from some discretely observed real data set.
In the case of the one-dimensional CIR process, the parameter estimation of \ $a$ \ and \ $b$ \ goes back to
 Overbeck and Ryd\'en \cite{OveRyd}, Overbeck \cite{Ove}, and see also the very recent papers of
 Ben Alaya and Kebaier \cite{BenKeb1,BenKeb2}.
For asymptotic results on discrete time critical branching processes with
 immigration, one may refer to Wei and Winnicki \cite{WeiWin1} and \cite{WeiWin2}.

The process \ $(Y,X)$ \ given by \eqref{DD} is a very special affine process.
The set of affine processes contains a large class of important Markov processes such as
 continuous state branching processes and Orstein-Uhlenbeck processes.
Further, a lot of models in financial mathematics are also special affine
 processes such as the Heston model \cite{Hes}, the model due to Barndorff-Nielsen
 and Shephard \cite{BarShe} or the model due to Carr and Wu \cite{CarWu}.
A precise mathematical formulation and complete characterization of regular
 affine processes are due to Duffie et al. \cite{DufFilSch}.
Later several authors have contributed to the study of properties of general affine processes:
 to name a few, Andersen and Piterbarg \cite{AndPit} (moment explosions in stochastic volatility models),
 Dawson and Li \cite{DawLi} (jump-type SDE representation for two-dimensional affine processes),
 Filipovi\'{c} and Mayerhofer \cite{FilMay} (applications to the pricing of bond and stock options),
 Glasserman and Kim \cite{GlaKim} (the range of finite exponential moments and the convergence to
 stationarity in affine diffusion models), Jena et al. \cite{JenKimXin}
 (long-term and blow-up behaviors of exponential moments in multi-dimensional
 affine diffusions), Keller-Ressel et al. \cite{KelSchTei1,KelSchTei2} (stochastically continuous, time-homogeneous affine
 processes with state space \ $\RR_+^n\times\RR^d$ \ or more general ones are regular).
We also refer to the overview articles Cuchiero et al. \cite{CucFilTei} and
 Friz and Keller-Ressel \cite{FriKel}.

To the best knowledge of the authors the parameter estimation problem for multi-dimensional affine processes has not been tackled so far.
Since affine processes are being used in financial mathematics very frequently, the question of parameter estimation for them
 is of high importance.
Our aim is to start the discussion with a simple non-trivial example: the two-dimensional affine diffusion process
 given by \eqref{DD}.

The article is divided into two parts and there are two appendices.
In Section~\ref{Section_scaling} we recall some notations,
 the definition of affine processes and some of their basic properties, and then a simple set of sufficient conditions
 for the weak convergence of scaled affine processes is presented.
Roughly speaking, given a family of affine processes
 \ $(Y^{(\theta)}(t), X^{(\theta)}(t))_{t\geq 0}$, \ $\theta>0$, \ such that
 the corresponding admissible parameters converge in an appropriate way
 (see Theorem \ref{Thm1}), the scaled process
 \ $\left( \theta^{-1} Y^{(\theta)}(\theta t), \theta^{-1}  X^{(\theta)}(\theta t) \right)_{t\geq 0}$ \
 converge weakly towards an affine diffusion process as \ $\theta \to \infty$.
\ We specialize our result for one-dimensional continuous state branching
 processes with immigration which generalizes Theorem 2.3 in Huang et al.\ \cite{HuaMaZhu}.
The scaling Theorem \ref{Thm1} is proved for quite general affine processes
 since it might have applications elsewhere later on.
In Section~\ref{Section_statistics} the scaling Theorem \ref{Thm1} is applied to study
 the asymptotic behavior of least squares and conditional least squares estimators of some parameters of a
 critical two-dimensional affine diffusion process  given by \eqref{DD}, see Theorems \ref{Thm2},
 \ref{Thm3} and \ref{Thm4}.
In Appendix \ref{AppendixA} we check that some integrals in the form of the infinitesimal generator of an affine process
 that we use are well-defined.
Appendix \ref{AppendixB} is devoted to show that the least squares estimator of \ $m$ \ cannot be asymptotically weakly consistent.

\section{A scaling theorem for affine processes}\label{Section_scaling}

Let \ $\NN$, \ $\ZZ_+$, \ $\RR$, \ $\RR_+$, \ $\RR_-$, \ $\RR_{++}$, \ and \ $\CC$
 \ denote the sets of positive integers, non-negative integers, real numbers,
 non-negative real numbers, non-positive real numbers, positive real numbers
 and complex numbers, respectively.
For \ $x , y \in \RR$, \ we will use the notations \ $x \land y := \min(x, y)$ \ and \ $x \vee y := \max(x, y)$.
\ For \ $x, y\in \CC^k$, $k\in\NN$, \ we write \ $\langle x,y\rangle := \sum_{i=1}^kx_iy_i$
 \ (notice that this is not the scalar product on \ $\CC^k$, \ however
 for \ $x\in\CC^k$ \ and \ $y\in\RR^k$, \ $\langle x,y\rangle$ \ coincides with the usual scalar product of
 \ $x$ \ and \ $y$).
\ By \ $\|x\|$ \ and \ $\|A\|$ \ we denote the Euclidean norm of a vector \ $x\in\RR^p$ \ and the induced matrix norm
 of a matrix \ $A\in\RR^{p\times p}$, \ respectively.
Further, let \ $U:=\{z_1+iz_2 : z_1\in\RR_-, \, z_2\in\RR\}\times(i\RR^d)$.
\ By \ $C^2_c(\RR_+\times\RR^d)$ \ ($C^\infty_c(\RR_+\times\RR^d)$) \ we denote
 the set of twice (infinitely) continuously differentiable complex-valued
 functions on \ $\RR_+\times\RR^d$ \ with compact support, where \ $d\in\NN$.
\ The set of c\`{a}dl\`{a}g functions from \ $\RR_+$ \ to \ $\RR_+\times\RR^d$
 \ will be denoted by \ $\DD(\RR_+,\RR_+\times\RR^d)$.
\ For a bounded function \ $g:\RR_+\times\RR^d\to\RR^p$, \ let
 \ $\|g\|_\infty := \sup_{x\in\RR_+\times\RR^d} \|g(x)\|$.
\ Convergence in distribution, in probability and almost sure convergence will be denoted by \ $\distr$, \ $\stoch$ \
 and \ $\as$, \ respectively.

Next we briefly recall the definition of affine processes with state space
 \ $\RR_+\times\RR^d$ \ based on Duffie et al.\ \cite{DufFilSch}.

\begin{Def}
A transition semigroup \ $(P_t)_{t\in\RR_+}$ \ with state space
 \ $\RR_+\times\RR^d$ \ is called a (general) affine semigroup if its
 characteristic function has the representation
 \begin{align}\label{affine_SG}
   \int_{\RR_+\times\RR^d} \ee^{\langle u,\xi\rangle} P_t(x,\dd\xi)
     = \ee^{\langle x,\psi(t,u)\rangle  +  \phi(t,u)}
 \end{align}
 for \ $x\in \RR_+\times\RR^d$, \ $u\in U$ \ and \ $t\in\RR_+$, \ where
 \ $\psi(t,\cdot)=(\psi_1(t,\cdot),\psi_2(t,\cdot))\in\CC\times\CC^d$ \ is a continuous
 \ $\CC^{1+d}$-valued function on \ $U$ \ and
 \ $\phi(t,\cdot)$ \ is a continuous \ $\CC$-valued function on \ $U$ \ satisfying \ $\phi(t,0)=0$.
\ The affine semigroup \ $(P_t)_{t\in\RR_+}$ \ defined by \eqref{affine_SG} is
 called regular if it is stochastically continuous (equivalently, for all \ $u\in U$,
 \ the functions \ $\RR_+\ni t\mapsto \Psi(t,u)$ \ and \ $\RR_+\ni t\mapsto \phi(t,u)$ \ are continuous)
 and \ $\partial_1\psi(0,u)$ \ and \ $\partial_1\phi(0,u)$ \ exist for all
 \ $u\in U$ \ and are continuous at \ $u=0$ \ (where \ $\partial_1\psi$ \ and
 \ $\partial_1\phi$ \ denote the partial derivatives of \ $\psi$ \ and
 \ $\phi$, \ respectively, with respect to the first variable).
\end{Def}

\begin{Rem}
We call the attention that Duffie et al. \cite{DufFilSch} in their Definition 2.1
 assume only that Equation \eqref{affine_SG} hold for \ $x\in \RR_+\times\RR^d$,
 $u\in \partial U=i\RR^{1+d}$, $t\in\RR_+$, \ i.e., instead of \ $u\in U$ \ they only require that
 \ $u$ \ should be an element of the boundary \ $\partial U$ \ of \ $U$.
\ However, by Proposition 6.4 in Duffie et al. \cite{DufFilSch}, one can formulate the definition of
 a regular affine process as we did.
Note also that this kind of definition was already given by Dawson and Li \cite[Definitions 2.1 and 3.3]{DawLi}.
Finally, we remark that every stochastically continuous affine semigroup is regular
 due to Keller-Ressel et al. \cite[Theorem 5.1]{KelSchTei1}.\proofend
\end{Rem}

\begin{Def}\label{Def_admissible}
A set of parameters \ $( a, \alpha, b, \beta, m, \mu)$ \ is called
 admissible if
 \renewcommand{\labelenumi}{{\rm(\roman{enumi})}}
 \begin{enumerate}
  \item
   $a = (a_{i,j})_{i,j=1}^{1+d} \in \RR^{(1+d)\times(1+d)}$ \ is a symmetric positive
    semidefinite matrix with \ $a_{1,1}=0$
    \ (hence \ $a_{1,k} = a_{k,1} = 0$ \ for all \ $k \in \{2,\ldots,1+d\}$),
  \item
   $\alpha = (\alpha_{i,j})_{i,j=1}^{1+d} \in \RR^{(1+d)\times(1+d)}$ \ is a symmetric
   positive semidefinite matrix,
  \item
   $b = (b_i)_{i=1}^{1+d} \in \RR_+\times\RR^d$,
  \item
   $\beta = (\beta_{i,j})_{i,j=1}^{1+d} \in \RR^{(1+d)\times(1+d)}$ \ with
    \ $\beta_{1,j}=0$ \ for all \ $j \in \{2,\ldots,1+d\}$,
  \item
   $m(\dd\xi) = m(\dd\xi_1, \dd\xi_2)$ \ is a \ $\sigma$-finite measure on
    \ $\RR_+\times\RR^d$ \ supported by
    \ $(\RR_+\times\RR^d)\setminus\{(0,0)\}$ \ such that
    \[
      \int_{\RR_+\times\RR^d}
       \left[ \xi_1 + (\|\xi_2\| \land \|\xi_2\|^2) \right] m(\dd\xi)
      <\infty,
    \]
  \item
   $\mu(\dd\xi) = \mu(\dd\xi_1, \dd\xi_2)$ \ is a $\sigma$-finite measure on
    \ $\RR_+\times\RR^d$ \ supported by
    \ $(\RR_+\times\RR^d)\setminus\{(0,0)\}$ \ such that
    \[
      \int_{\RR_+\times\RR^d} \|\xi\| \land \|\xi\|^2 \mu(\dd\xi)
      <\infty.
    \]
  \end{enumerate}
\end{Def}

\begin{Rem}
Note that our Definition \ref{Def_admissible} of the set of admissible parameters is not
 so general as Definition 2.6 in Duffie et al. \cite{DufFilSch}.
Firstly, the set of admissible parameters is defined only for affine process with state space \ $\RR_+\times\RR^d$,
 \ while Duffie et al. \cite{DufFilSch} consider affine processes with state space \ $\RR_+^n\times\RR^d$.
We restrict ourselves to this special case, since our scaling Theorem \ref{Thm1} is valid only in this case.
Secondly, our conditions (v) and (vi) of Definition \ref{Def_admissible} are stronger than that of
 (2.10) and (2.11) of Definition 2.6 in Duffie et al. \cite{DufFilSch}.
Thirdly, according to our definition, a set of admissible parameters does not contain parameters corresponding to killing,
 while in Definition 2.6 in Duffie et al. \cite{DufFilSch} such parameters are included.
Our definition of admissible parameters can be considered as a \ $(1+d)$-dimensional version of Definition 6.1
 in Dawson and Li \cite{DawLi}.
The reason for this definition is to have a more pleasant form of the infinitesimal generator of an affine process
 compared to that of Duffie et al. \cite[formula (2.12)]{DufFilSch}.
For more details, see Remark \ref{Rem4}.\proofend
\end{Rem}

\begin{Thm}{(Duffie et al. \cite[Theorem 2.7]{DufFilSch})}\label{Thm_Duffie_et_al}
Let \ $(a, \alpha, b, \beta, m, \mu)$ \ be a set of admissible parameters.
Then there exists a unique regular affine semigroup \ $(P_t)_{t\in\RR_+}$ \ with infinitesimal generator
  \begin{align}\label{affine_inf_gen}
  \begin{split}
   (\cA f)(x)
   &= \sum_{i,j=1}^{1+d} (a_{i,j} + \alpha_{i,j} x_1) f_{i,j}''(x)
      + \langle f'(x), b + \beta x \rangle \\
   &\phantom{\quad}
      + \int_{\RR_+\times \RR^d}
         (f(x+\xi) - f(x) - \langle f'_{(2)}(x), \xi_2 \rangle)
         \, m(\dd\xi)\\
   &\phantom{\quad}
      + \int_{\RR_+\times \RR^d}
         (f(x+\xi) - f(x) - \langle f'(x), \xi \rangle) x_1
         \, \mu(\dd\xi)
 \end{split}
 \end{align}
 for \ $x=(x_1,x_2)\in\RR_+\times\RR^d$ \ and \ $f\in C^2_c(\RR_+\times\RR^d)$, \ where \ $f_i'$,
 \ $i\in\{1,\ldots,1+d\}$, \ and \ $f_{i,j}''$, \ $i,j\in\{1,\ldots,1+d\}$,
 \ denote the first and second order partial derivatives of $f$ with respect to
 its \ $i$-th and \ $i$-th and \ $j$-th variables, and
 \ $f'(x) := (f_1'(x),\ldots,f_{1+d}'(x))^\top$,
 \ $f'_{(2)}(x) := (f_2'(x),\ldots,f_{1+d}'(x))^\top$.
Further, \ $\cC_c^\infty(\RR_+\times\RR^d)$ \ is a core of \ $\cA$.
\end{Thm}

\begin{Rem}\label{Rem4}
Note that the form of the infinitesimal generator \ $\cA$ \ in Theorem \ref{Thm_Duffie_et_al} is slightly different
 from the one given in (2.12) in Duffie et al. \cite{DufFilSch}.
Our formula \eqref{affine_inf_gen} is in the spirit of Dawson and Li \cite[formula (6.5)]{DawLi}.
On the one hand, the point is that under the conditions (v) and (vi) of Definition \ref{Def_admissible},
 one can rewrite (2.12) in Duffie et al. \cite{DufFilSch} into the form \eqref{affine_inf_gen},
 by changing the \ $2$-nd, $\ldots$, $(1+d)$-th coordinates of \ $b\in\RR_+\times\RR^d$ \ and the first column of
 \ $\beta\in\RR^{(1+d)\times(1+d)}$, \ respectively, in appropriate ways (see Appendix \ref{AppendixA}).
To see this, it is enough to check that the integrals in \eqref{affine_inf_gen} are well-defined
 (i.e., elements of \ $\CC$) \ under the conditions (v) and (vi) of Definition \ref{Def_admissible}.
For further details, see also Appendix~\ref{AppendixA}.
On the other hand, the killing rate (see page 995 in Duffie et al. \cite{DufFilSch}) of the affine semigroup
 \ $(P_t)_{t\in\RR_+}$ \ in Theorem \ref{Thm_Duffie_et_al} is identically zero.
This also implies that the affine processes that we will consider
later on will have lifetime infinity.\proofend\looseness=-1
\end{Rem}

\begin{Rem}
In dimension 2 (i.e., if \ $d=1$), by Theorem 6.2 in Dawson and Li \cite{DawLi} and Theorem 2.7 in Duffie et al. \cite{DufFilSch}
 (see also Theorem \ref{Thm_Duffie_et_al}), for an infinitesimal generator \ $\cA$ \ given by \eqref{affine_inf_gen} with \ $d=1$ \
 one can construct a two-dimensional system of jump type SDEs of which there exists a pathwise unique strong solution
 \ $(Y(t),X(t))_{t\in\RR_+}$ \ which is a regular affine Markov process with the given infinitesimal generator \ $\cA$.\proofend
\end{Rem}

The next lemma is simple but very useful.

\begin{Lem}\label{Lem_inf_gen}
Let \ $(Z(t))_{t\in\RR_+}$ \ be a time-homogeneous Markov process with state
 space \ $\RR_+\times\RR^d$ \ and let us denote its infinitesimal generator by
 \ $\cA_Z$.
\ Suppose that \ $\cC_c^2(\RR_+\times\RR^d)$ \ is a subset of the domain of
 \ $\cA_Z$.
\ Then for all \ $\theta\in\RR_{++}$, \ the time-homogeneous Markov process
 \ $(Z_\theta(t))_{t\in\RR_+} := (\theta^{-1}Z(\theta t))_{t\in\RR_+}$ \ has
 infinitesimal generator
 \begin{align*}
  (\cA_{Z_\theta}f)(x) = \theta (\cA_Z f_\theta)(\theta x), \qquad
  x\in \RR_+\times\RR^d, \quad f\in C_c^2(\RR_+\times\RR^d),
 \end{align*}
 where \ $f_\theta(x) := f(\theta^{-1}x)$, \ $x\in \RR_+\times\RR^d$.
\end{Lem}

\noindent{\bf Proof.}
By definition, the infinitesimal generator of \ $(Z_\theta(t))_{t\in\RR_+}$
 \ takes the form
 \begin{align*}
  (\cA_{Z_\theta}f)(x)
  &= \lim_{t\downarrow 0}
      \frac{\EE( f(\theta^{-1} Z(\theta t)) \mid \theta^{-1}Z(0)=x) - f(x)}{t}\\
  &= \lim_{t\downarrow 0}
      \frac{\theta \EE( f_\theta(Z(\theta t))\mid Z(0)=\theta x)
            - \theta f_\theta(\theta x)}
           {\theta t}\\
  &= \theta
     \lim_{t'\downarrow 0}
      \frac{\EE( f_\theta(Z(t')) \mid Z(0)=\theta x) - f_\theta(\theta x)}{t'}
    = \theta (\cA_Z f_\theta)(\theta x)
 \end{align*}
 for all \ $x\in \RR_+\times\RR^d$ \ and \ $f\in C_c^2(\RR_+\times\RR^d)$.
\proofend

\begin{Thm}\label{Thm1}
For all \ $\theta\in\RR_{++}$, \ let
 \ $(Y^{(\theta)}(t), X^{(\theta)}(t))_{t\in\RR_+}$ \ be a \ $(1+d)$-dimensional
 affine process with state space \ $\RR_+\times\RR^d$ \ and with admissible
 parameters
 \ $(a^{(\theta)}, \alpha^{(\theta)}, b^{(\theta)}, \beta^{(\theta)}, m, \mu)$
 \ such that additionally
 \begin{align}\label{cond1}
  \int_{\RR_+\times\RR^d} \|\xi\| \, m(\dd\xi)<\infty
    \qquad\text{and}\qquad
  \int_{\RR_+\times\RR^d} \|\xi\|^2 \, \mu(\dd\xi)<\infty.
 \end{align}
Let \ $a, \alpha, \beta \in \RR^{(1+d)\times(1+d)}$,
 \ $b\in\RR_+\times\RR^d$, \ and let
 \ $(Y(t), X(t))_{t\in\RR_+}$ \ be a
 \ $(1+d)$-dimensional affine process with state space \ $\RR_+\times\RR^d$
 \ and with the set of admissible parameters \ $(a, \widetilde{\alpha}, \widetilde{b}, \beta, 0, 0)$, \ where
 \[
   \widetilde{\alpha}
   := \alpha
      + \frac{1}{2} \int_{\RR_+\times\RR^d} \xi \xi^\top \, \mu(\dd\xi) ,
 \]
 and \ $\widetilde{b} = (\widetilde{b}_i)_{i=1}^{1+d}$ \ with
 \ $\widetilde{b}_i := b_i$ \ for \ $i \in \{2,\ldots,1+d\}$ \ and
 \[
   \widetilde{b}_1 := b_1 + \int_{\RR_+\times\RR^d} \xi_1 \, m(\dd\xi) .
 \]
If
 \begin{gather*}
  \theta^{-1} a^{(\theta)} \to a, \qquad
   \alpha^{(\theta)} \to \alpha, \qquad
   b^{(\theta)} \to b, \qquad
   \theta \beta^{(\theta)} \to \beta, \\
  \theta^{-1}(Y^{(\theta)}(0), X^{(\theta)}(0)) \distr (Y(0), X(0))
 \end{gather*}
 as \ $\theta\to\infty$, \ then
 \begin{align*}
  \left( Y^{(\theta)}_\theta(t), X^{(\theta)}_\theta(t)\right)_{t\in\RR_+}
  = \left(\theta^{-1}Y^{(\theta)}(\theta t),
                          \theta^{-1}X^{(\theta)}(\theta t)\right)_{t\in\RR_+}
  \distr (Y(t),X(t))_{t\in\RR_+}
 \end{align*}
 in \ $\DD(\RR_+,\RR_+\times\RR^d)$ \ as \ $\theta\to\infty$.
\end{Thm}

\begin{Rem}
(i) Note that the limit process \ $(Y(t), X(t))_{t\in\RR_+}$ \ in Theorem \ref{Thm1}
 has continuous sample paths almost surely.
However, this is not a big surprise, since in condition \eqref{cond1} of Theorem \ref{Thm1}
 we require finite second moment for the measure \ $\mu$.

\noindent (ii)
Note also that the matrix \ $\widetilde{\alpha}\in\RR^{(1+d)\times(1+d)}$ \ given in Theorem \ref{Thm1}
 is symmetric and positive semidefinite, since \ $\alpha$ \ is symmetric and positive semidefinite, and for all \ $z\in\RR^{1+d}$,
 \begin{align*}
 \left\langle \int_{\RR_+\times\RR^d} \xi \xi^\top \, \mu(\dd\xi) z,z \right\rangle
        = \int_{\RR_+\times\RR^d} (z^\top\xi)^2 \, \mu(\dd\xi)
       \geq 0.\\[-26pt]
 \end{align*}
\proofend
\end{Rem}

\noindent{\bf Proof of Theorem \ref{Thm1}.}
By Duffie et al. \cite[Theorem 2.7]{DufFilSch},
 \ $\cC_c^\infty(\RR_+\times\RR^d)$ \ is a core of the infinitesimal generator
 \ $\cA_{(Y,X)}$ \ of the process \ $(Y(t), X(t))_{t\in\RR_+}$, \ and hence
 \ $\{(f,\cA_{(Y,X)}f) : f\in\cD(\cA_{(Y,X)})\}$ \ coincides with the closure of
 \ $\{(f,\cA_{(Y,X)}f) : f\in\cC_c^\infty(\RR_+\times\RR^d)\}$, \ where
 \ $\cD(\cA_{(Y,X)})$ \ denotes the domain of \ $\cA_{(Y,X)}$, \ see, e.g.,
 Ethier and Kurtz \cite[page 17]{EthKur}.
In other words, the closure of
 \ $\{(f,\cA_{(Y,X)}f) : f\in\cC_c^\infty(\RR_+\times\RR^d)\}$ \ generates the
 affine semigroup corresponding to \ $\cA_{(Y,X)}$.

Next we show that for all \ $f\in \cC_c^\infty(\RR_+\times\RR^d)$, \ we have
 \begin{align}\label{help2}
  \lim_{\theta\to\infty} \sup_{x\in\RR_+\times\RR^d}
   \Big| (\cA_{(Y_\theta^{(\theta)},X_\theta^{(\theta)})}f)(x) - (\cA_{(Y,X)} f)(x) \Big|
   =0.
 \end{align}
First note that it is enough to prove \eqref{help2} for real-valued functions \ $f\in \cC_c^\infty(\RR_+\times\RR^d)$,
 \ since if \eqref{help2} holds for for real-valued functions \ $f\in \cC_c^\infty(\RR_+\times\RR^d)$,
 \ then, by decomposing \ $f$ \ into real and imaginary parts, the linearity of
 the infinitesimal generators in question and  triangular inequality yield
 \eqref{help2} for complex-valued functions \ $f\in \cC_c^\infty(\RR_+\times\RR^d)$.
Hence in what follows without loss of generality we can assume that
 \ $f\in \cC_c^\infty(\RR_+\times\RR^d)$ \ is real-valued.

\allowdisplaybreaks
For all \ $f\in\cC_c^\infty(\RR_+\times\RR^d)$, \ $\theta\in\RR_{++}$, \ and
 \ $x\in \RR_+\times\RR^d$, \ we have
 \begin{align}\label{help1}
  \begin{split}
   & f_\theta(x) = f(\theta^{-1}x),\\
   & (f_\theta)_i'(x) = \theta^{-1} f_i'(\theta^{-1}x),\qquad
     i\in\{1,\ldots,1+d\},\\
   & (f_\theta)_{i,j}'' (x) = \theta^{-2} f_{i,j}'' (\theta^{-1}x),\qquad
     i,j\in\{1,\ldots,1+d\}.
  \end{split}
 \end{align}
Then, by Lemma \ref{Lem_inf_gen}, \eqref{affine_inf_gen} and \eqref{help1},
 \begin{align*}
  (\cA_{(Y_\theta^{(\theta)}, X_\theta^{(\theta)})}f)(x)
  &= \theta (\cA_{(Y^{(\theta)}, X^{(\theta)})}f_\theta)(\theta x)\\
  &=\theta
    \Biggl[ \sum_{i,j=1}^{1+d}
             (a_{i,j}^{(\theta)} + \alpha_{i,j}^{(\theta)} \theta x_1)
             \theta^{-2} f_{i,j}''(\theta^{-1} \theta x)
            + \langle \theta^{-1} f'(\theta^{-1} \theta x), b^{(\theta)} + \beta^{(\theta)} \theta x \rangle \\
  &\quad 
            + \int_{\RR_+\times \RR^d}
               (f(\theta^{-1}(\theta x+\xi)) - f(\theta^{-1} \theta x)
                - \langle \theta^{-1} f'_{(2)}(\theta^{-1} \theta x),
                          \xi_2 \rangle)
               \, m(\dd\xi) \\
  &\quad 
            + \int_{\RR_+\times \RR^d}
               (f(\theta^{-1}(\theta x+\xi)) - f(\theta^{-1} \theta x)
                - \langle \theta^{-1} f'(\theta^{-1} \theta x), \xi \rangle)
               \theta x_1 \, \mu(\dd\xi) \Biggr]\\
  &=\sum_{i,j=1}^{1+d}
     (\theta^{-1} a_{i,j}^{(\theta)} + \alpha_{i,j}^{(\theta)} x_1) f_{i,j}''(x)
     + \langle f'(x), b^{(\theta)} + \theta \beta^{(\theta)} x \rangle \\
  &\phantom{\quad}
      + \int_{\RR_+\times \RR^d}
         (f(x+\theta^{-1}\xi) - f(x)
          - \langle f'_{(2)}(x), \theta^{-1} \xi_2 \rangle) \theta
         \, m(\dd\xi) \\
  &\phantom{\quad}
      + x_1 \int_{\RR_+\times \RR^d}
             (f(x+\theta^{-1}\xi) - f(x)
              - \langle f'(x), \theta^{-1} \xi \rangle) \theta^2
             \, \mu(\dd\xi)
 \end{align*}
for \ $f\in\cC_c^\infty(\RR_+\times\RR^d)$ \ and \ $x\in\RR_+\times\RR^d$.
\ Hence, for all \ $x=(x_1,x_2)\in\RR_+\times\RR^d$, \ using the triangular inequality and that
 \ $\vert \langle u,v\rangle\vert \leq \Vert u\Vert \Vert v\Vert$, $u,v\in\RR^p$, \ we have
 \begin{align*}
  &\Bigl| (\cA_{(Y_\theta^{(\theta)},X_\theta^{(\theta)})}f)(x)
          - (\cA_{(Y,X)} f)(x) \Bigr| \\
  &\leq \sum_{i,j=1}^{1+d}
         ( |\theta^{-1} a_{i,j}^{(\theta)} - a_{i,j}|
           + |\alpha_{i,j}^{(\theta)} - \alpha_{i,j}| x_1 ) |f_{i,j}''(x)|\\
  &\phantom{\quad}
     + ( \|b^{(\theta)} - b\| + \|\theta \beta^{(\theta)} -  \beta\| \|x\| )
          \|f'(x)\| \\
  &\phantom{\quad}
        + \left| \int_{\RR_+\times \RR^d}
                  \Big( f(x+\theta^{-1}\xi) - f(x)
                        - \theta^{-1} \langle f'(x), \xi \rangle \Big)
                  \theta \, m(\dd\xi) \right| \\
  &\phantom{\quad}
        + x_1 \Bigg| \int_{\RR_+\times \RR^d}
                      \Big( f(x+\theta^{-1}\xi) - f(x)
                            - \theta^{-1} \langle f'(x), \xi \rangle
                            - \frac{1}{2} \theta^{-2}
                              \langle f''(x) \xi, \xi \rangle \Big)
                      \theta^2 \, \mu(\dd\xi) \Bigg| ,
 \end{align*}
 where
 \[
   f''(x) := \begin{bmatrix}
              f_{1,1}''(x)  & \cdots & f_{1,1+d}''(x) \\
              \vdots       & \ddots & \vdots \\
              f_{1+d,1}''(x) & \cdots & f_{1+d,1+d}''(x)
             \end{bmatrix} .
 \]
Since \ $f\in\cC_c^\infty(\RR_+\times\RR^d)$, \ we have
 \begin{gather*}
  \sup_{x\in\RR_+\times\RR^d} x_1 |f_{i,j}''(x)| < \infty, \qquad
  \sup_{x\in\RR_+\times\RR^d} |f_{i,j}''(x)| < \infty, \qquad
    \forall\; i,j\in\{1,\ldots,1+d\}, \\
  \sup_{x\in\RR_+\times\RR^d} \|x\| \| f'(x) \| < \infty, \qquad
  \sup_{x\in\RR_+\times\RR^d} \| f'(x) \| < \infty,
 \end{gather*}
 and hence, by our assumptions, in order to prove \eqref{help2} it is enough to
 check that
 \begin{align}\label{help3}
  \lim_{\theta\to\infty} \sup_{x\in\RR_+\times\RR^d}
   \left| \int_{\RR_+\times \RR^d}
           \Big( f(x+\theta^{-1}\xi) - f(x)
                 - \theta^{-1} \langle f'(x), \xi \rangle \Big)
           \theta \, m(\dd\xi) \right|
  =0,
 \end{align}
 and
 \begin{align}\label{help4}
 \begin{split}
  \lim_{\theta\to\infty} \sup_{x\in\RR_+\times\RR^d}
   x_1 & \left| \int_{\RR_+\times \RR^d}
               \Big( f(x+\theta^{-1}\xi) - f(x)
                     - \theta^{-1} \langle f'(x), \xi \rangle\right.\\
      &\;\left. - \frac{1}{2} \theta^{-2}
                       \langle f''(x) \xi, \xi \rangle \Big)
                \theta^2 \, \mu(\dd\xi) \right|
  =0.
 \end{split}
 \end{align}

First we consider \eqref{help3}.
Let \ $\varepsilon\in\RR_{++}$ \ be fixed.
Let us choose an \ $M\in\RR_{++}$ \ such that
 \begin{align}\label{help6}
  2 \|f'\|_\infty
  \int_{(\RR_+\times\RR^d)\setminus ([0,M]\times [-M,M]^d)} \|\xi\| \, m(\dd\xi)
  < \frac{\varepsilon}{2}.
 \end{align}
In what follows, for abbreviation, \ $[0,M] \times [-M,M]^d$ \ will be denoted
 by \ $D_M$.
\ Such an \ $M$ \ can be chosen, since \ $f\in\cC_c^\infty(\RR_+\times\RR^d)$
 \ yields \ $\|f'\|_\infty<\infty$ \ and, by assumption \eqref{cond1},
 \ $\int_{\RR_+\times\RR^d} \|\xi\| \, m(\dd\xi)<\infty$.
\ By Taylor's theorem, for all \ $\theta\in\RR_{++}$, \ $x\in\RR_+\times\RR^d$
 \ and \ $\xi\in D_M$ \ there exists some
 \ $\tau = \tau(\theta, x, \xi) \in[0,1]$ \ such that
 \begin{align}\label{help15}
  f(x+\theta^{-1}\xi) - f(x)
  = \langle f'(x+\theta^{-1}\tau\xi), \theta^{-1}\xi \rangle.
 \end{align}
 Then
 \begin{multline*}
  \left| \int_{\RR_+\times \RR^d}
          \Bigl(f(x+\theta^{-1}\xi) - f(x)
                - \theta^{-1} \langle f'(x), \xi \rangle \Bigr)
          \theta \, m(\dd\xi) \right| \\
  \leq \int_{\RR_+\times \RR^d}
        \bigl| \langle f'(x+\theta^{-1}\tau\xi), \xi \rangle
               - \langle f'(x), \xi \rangle \bigr| \, m(\dd\xi)
  \leq A^{(1)}_{\theta,M}(x) + A^{(2)}_{\theta,M}(x)
 \end{multline*}
 for all \ $x\in\RR_+\times\RR^d$, \ where
 \begin{align*}
  &A^{(1)}_{\theta,M}(x)
   := \int_{D_M}
       | \langle f'(x+\theta^{-1}\tau\xi) - f'(x), \xi \rangle |
       \, m(\dd\xi) , \\
  &A^{(2)}_{\theta,M}(x)
   := \int_{(\RR_+\times\RR^d)\setminus D_M}
       \left( |\langle f'(x+\theta^{-1}\tau\xi), \xi \rangle|
              + |\langle f'(x), \xi \rangle| \right)
       m(\dd\xi) .
 \end{align*}
Here
 \begin{align*}
  A^{(1)}_{\theta,M}(x)
  &\leq \int_{D_M}
         \|f'(x+\theta^{-1}\tau\xi) - f'(x)\| \|\xi\| \, m(\dd\xi) \\
  &\leq \sup_{\xi\in D_M} \|f'(x+\theta^{-1}\tau\xi) - f'(x)\|
        \int_{\RR_+\times\RR^d} \|\xi\| \, m(\dd\xi) .
 \end{align*}
The convexity of \ $D_M$ \ implies
 \ $\tau \xi = \tau(\theta, x, \xi) \xi \in D_M$ \ for all
 \ $\theta \in \RR_{++}$, \ $x \in \RR_+\times\RR^d$ \ and \ $\xi \in D_M$,
 \ and, hence,
 \[
   \sup_{\xi\in D_M} \|f'(x+\theta^{-1}\tau\xi) - f'(x)\|
   \leq \sup_{\widetilde\xi\in D_M} \|f'(x+\theta^{-1}\widetilde\xi) - f'(x)\| .
 \]
Since \ $f'$ \ is uniformly continuous on \ $\RR_+\times\RR^d$ \ (which follows by mean value theorem using also that
 \ $\Vert f''\Vert_\infty<\infty$), \ there exists a \ $\theta_0\in\RR_{++}$ \ (depending on \ $\varepsilon$ \ and \ $M$)
 \ such that
 \[
   \sup_{\widetilde\xi\in D_M}
    \|f'(x+\theta^{-1}\widetilde\xi) - f'(x)\|
   \int_{\RR_+\times\RR^d} \|\xi\| \, m(\dd\xi)
   < \frac{\varepsilon}{2}
 \]
 for all \ $x\in\RR_+\times\RR^d$ \ and \ $\theta \in [\theta_0,\infty)$.
\ Further, by \eqref{help6}, we have
 \begin{align*}
  A^{(2)}_{\theta,M}(x)
  \leq 2 \|f'\|_\infty \int_{(\RR_+\times\RR^d)\setminus D_M} \|\xi\| \, m(\dd\xi)
  < \frac{\varepsilon}{2}
 \end{align*}
 for all \ $x\in\RR_+\times\RR^d$ \ and \ $\theta\in\RR_{++}$.
\ Putting the pieces together we have \eqref{help3}.

Now we turn to prove \eqref{help4} in a similar way.
Let \ $\varepsilon\in\RR_{++}$ \ be fixed again.
Let us now choose an \ $M\in\RR_{++}$ \ such that
 \begin{align}\label{2help6}
  2 \sup_{x \in \RR_+\times\RR^d} x_1 \|f''(x)\|
  \int_{(\RR_+\times\RR^d)\setminus D_M} \|\xi\|^2 \, \mu(\dd\xi)
  < \frac{\varepsilon}{2}.
 \end{align}
Such an $M$ can be chosen, since
 $\sup_{x \in \RR_+\times\RR^d} x_1 \|f''(x)\|<\infty$ for all $f\in\cC_c^\infty(\RR_+\times\RR^d)$ and,
 by assumption \eqref{cond1},
 \ $\int_{\RR_+\times\RR^d} \|\xi\|^2 \, \mu(\dd\xi)<\infty$.
\ By Taylor's theorem, for all \ $\theta\in\RR_{++}$, \ $x\in\RR_+\times\RR^d$
 \ and \ $\xi\in D_M$ \ there exists some
 \ $\tau = \tau(\theta, x, \xi) \in[0,1]$ \ such that
 \begin{align}\label{2help15}
  f(x+\theta^{-1}\xi) - f(x) - \langle f'(x), \theta^{-1}\xi \rangle
  = \frac{1}{2}
    \langle f''(x+\theta^{-1}\tau\xi) \theta^{-1}\xi, \theta^{-1}\xi \rangle.
 \end{align}
Then
 \begin{align*}
  &x_1 \left| \int_{\RR_+\times \RR^d}
              \Big(f(x+\theta^{-1}\xi) - f(x)
                   - \theta^{-1} \langle f'(x), \xi \rangle
                   - \frac{1}{2}
                     \theta^{-2} \langle f''(x) \xi , \xi \rangle \Big)
              \theta^2 \, \mu(\dd\xi) \right| \\
  &\quad\leq \frac{1}{2} x_1
       \int_{\RR_+\times \RR^d}
        \bigl| \langle f''(x+\theta^{-1}\tau\xi) \xi, \xi \rangle
               - \langle f''(x) \xi, \xi \rangle \bigr| \mu(\dd\xi)\\
  &\quad\leq \frac{1}{2} ( B^{(1)}_{\theta,M}(x) + B^{(2)}_{\theta,M}(x) )
 \end{align*}
 for all \ $x\in\RR_+\times\RR^d$, \ where
 \begin{align*}
  &B^{(1)}_{\theta,M}(x)
   := x_1 \int_{D_M}
           | \langle (f''(x+\theta^{-1}\tau\xi) - f''(x)) \xi, \xi \rangle |
           \, \mu(\dd\xi) , \\
  &B^{(2)}_{\theta,M}(x)
   := x_1 \int_{(\RR_+\times\RR^d)\setminus D_M}
           \left( |\langle f''(x+\theta^{-1}\tau\xi) \xi, \xi \rangle|
                  + |\langle f''(x) \xi, \xi \rangle| \right)
           \mu(\dd\xi) .
 \end{align*}
Here
 \begin{align*}
  B^{(1)}_{\theta,M}(x)
  &\leq x_1 \int_{D_M}
             \|f''(x+\theta^{-1}\tau\xi) - f''(x)\| \|\xi\|^2 \, \mu(\dd\xi) \\
  &\leq \sup_{\xi\in D_M} x_1 \|f''(x+\theta^{-1}\tau\xi) - f''(x)\|
        \int_{\RR_+\times\RR^d} \|\xi\|^2 \, \mu(\dd\xi) ,
 \end{align*}
 and note that \ $\Vert A-B\Vert\geq \Vert A\Vert - \Vert B\Vert$, $A,B\in\RR^{p\times p}$, \ yields that
 \begin{align*}
   x_1 \|f''(x+\theta^{-1}\tau\xi) - f''(x)\|
    \leq \|(x_1+\theta^{-1}\tau\xi_1) f''(x+\theta^{-1}\tau\xi) - x_1 f''(x)\|
         + \theta^{-1}\tau\xi_1 \|f''(x+\theta^{-1}\tau\xi)\| .
 \end{align*}
Further, we have again
 \[
   \{ \tau \xi = \tau(\theta, x, \xi) \xi
       : \text{$\theta \in \RR_{++}$, \ $x \in \RR_+\times\RR^d$
         \ and \ $\xi \in D_M$} \}
    \subset D_M,
  \]
  and hence
 \begin{align*}
   \sup_{\xi\in D_M}
    &\|(x_1+\theta^{-1}\tau\xi_1)f''(x+\theta^{-1}\tau\xi) - x_1f''(x)\|\\
    &\leq \sup_{\widetilde\xi\in D_M}
         \|(x_1+\theta^{-1}\widetilde\xi_1)f''(x+\theta^{-1}\widetilde\xi)
            - x_1f''(x)\| .
 \end{align*}
Since the function \ $x \mapsto x_1 f''(x)$ \ is uniformly continuous on
 \ $\RR_+\times\RR^d$, \ there exists a \ $\theta_1\in\RR_{++}$ \ (depending on
 \ $\varepsilon$ \ and \ $M$) \ such that
 \[
   \sup_{\widetilde\xi\in D_M}
    \|(x_1+\theta^{-1}\widetilde\xi_1)f''(x+\theta^{-1}\widetilde\xi)
      - x_1f''(x)\|
    \int_{\RR_+\times\RR^d} \|\xi\|^2 \, \mu(\dd\xi)
   < \frac{\varepsilon}{4}
 \]
 for all \ $x\in\RR_+\times\RR^d$ \ and \ $\theta \in [\theta_1,\infty)$.
\ Moreover, there exists a \ $\theta_2\in\RR_{++}$
 \ (depending on \ $\varepsilon$ \ and \ $M$) \ such that
 \begin{align*}
   \theta^{-1}\sup_{\xi\in D_M} \tau\xi_1 \|f''(x+\theta^{-1}\tau\xi)\|
   &\int_{\RR_+\times\RR^d} \|\xi\|^2 \, \mu(\dd\xi)\\
   &\leq \theta^{-1} \|f''\|_\infty \sup_{\xi\in D_M} \xi_1
        \int_{\RR_+\times\RR^d} \|\xi\|^2 \, \mu(\dd\xi)
    < \frac{\varepsilon}{4}
 \end{align*}
 for all \ $x\in\RR_+\times\RR^d$ \ and \ $\theta \in [\theta_2,\infty)$.
\ Consequently, \ $B^{(1)}_{\theta,M}(x) < \frac{\varepsilon}{2}$ \ for all
 \ $x\in\RR_+\times\RR^d$ \ and \ $\theta \in [\theta_1 + \theta_2, \infty)$.
\ Further,
 \begin{align*}
  B^{(2)}_{\theta,M}(x)
  &\leq x_1 \int_{(\RR_+\times\RR^d)\setminus D_M}
             \left( \|f''(x+\theta^{-1}\tau\xi)\| + \|f''(x)\| \right)
             \|\xi\|^2 \, \mu(\dd\xi) \\
  &\leq \sup_{\xi\in \RR_+\times\RR^d}
         x_1 \left( \|f''(x+\theta^{-1}\tau\xi)\| + \|f''(x)\| \right)
        \int_{(\RR_+\times\RR^d)\setminus D_M} \|\xi\|^2 \, \mu(\dd\xi) .
 \end{align*}
Here
 \[
   x_1 \left( \|f''(x+\theta^{-1}\tau\xi)\| + \|f''(x)\| \right)
   \leq (x_1+\theta^{-1}\tau\xi_1)\|f''(x+\theta^{-1}\tau\xi)\| + x_1\|f''(x)\| ,
 \]
 hence
 \[
   \sup_{\xi\in \RR_+\times\RR^d}
    x_1 \left( \|f''(x+\theta^{-1}\tau\xi)\| + \|f''(x)\| \right)
   \leq 2 \sup_{x\in \RR_+\times\RR^d} x_1\|f''(x)\| .
 \]
Consequently, by \eqref{2help6}, we have
 \begin{align*}
  B^{(2)}_{\theta,M}(x)
  \leq 2 \sup_{x\in \RR_+\times\RR^d} x_1\|f''(x)\| \int_{(\RR_+\times\RR^d)\setminus D_M} \|\xi\|^2 \, \mu(\dd\xi)
  < \frac{\varepsilon}{2}
 \end{align*}
 for all \ $x\in\RR_+\times\RR^d$ \ and \ $\theta\in\RR_{++}$.
Putting the pieces together we have \eqref{help4}.

Finally, Ethier and Kurtz \cite[Corollary 8.7 on page 232]{EthKur} yields our
 assertion.
Namely, with the notations of part (f) of this corollary (but replacing \ $n$
 \ by \ $\theta$), \ let
 \begin{itemize}
  \item
   $G_\theta:= E_\theta:= E:=\RR_+\times\RR^d$ \ for all \ $\theta\in\RR_{++}$,
  \item
   $C_a := \cC_c^\infty(\RR_+ \times \RR^d)$ \ which strongly separates
    points in \ $\RR_+\times\RR^d$ \ (indeed, for every
    \ $(x_1,x_2)\in \RR_+ \times \RR^d$ \ and \ $\delta\in\RR_{++}$, \ the bump
    function defined by
    \ $h_1(u_1,u_2) := \exp\bigl\{- \frac{1}{1-(u_1-x_1)^2} - \frac{1}{1-\|u_2-x_2\|^2} \bigr\}$ \
    if \ $|u_1-x_1|<1$ \ and \ $\|u_2-x_2\|<1$ \ with \ $(u_1,u_2)\in \RR_+ \times \RR^d$, \
    and \ $h_1(u_1,u_2) := 0$ \ otherwise, satisfies (4.7) on page 113 in Ethier and Kurtz \cite{EthKur}),
  \item
   $\eta_\theta : E_\theta \to E$ \ with \ $\eta_\theta(x_1,x_2):=(x_1,x_2)$,
    \ $(x_1,x_2)\in E_\theta$ \ for all \ $\theta\in\RR_{++}$,
  \item
   $\pi_\theta : E \to E_\theta$ \ with \ $\pi_\theta(x_1,x_2):= (x_1,x_2)$,
    \ $(x_1,x_2)\in E$ \ for all \ $\theta\in\RR_{++}$,
  \item
   for each \ $f \in \cC_c^\infty(\RR_+ \times \RR^d)$ \ one can choose
    \ $f_\theta:=f$ \ and
    \ $g_\theta := \cA_{(Y^{(\theta)}_\theta, X^{(\theta)}_\theta)} f$ \ for all
    \ $\theta\in\RR_{++}$
    \ (and hence
    \ $(f_\theta, g_\theta) \in \widehat{\cA}_{(Y^{(\theta)}_\theta,X^{(\theta)}_\theta)}$
    \ defined on page 24 in Ethier and Kurtz \cite{EthKur} by part (c) of
    Proposition 1.5 on page 9 in Ethier and Kurtz \cite{EthKur}),
  \item
   $(\cG^\theta_t)_{t\in\RR_+} := (\cF_t^{(Y^{(\theta)}_\theta,X^{(\theta)}_\theta)})_{t\in\RR_+}$, \ where
    \ $\cF_t^{(Y^{(\theta)}_\theta,X^{(\theta)}_\theta)}$ \ denotes the $\sigma$-algebra generated by \ $\{(Y^{(\theta)}_\theta(s),X^{(\theta)}_\theta(s)), s\in[0,t]\}$.
 \end{itemize}
Then, by our assumptions, convergence of the initial distributions holds,
 condition (8.35) on page 232 in Ethier and Kurtz \cite{EthKur} is
 automatically satisfied, and \eqref{help2} shows the validity of condition
 (8.36) on page 232 in Ethier and Kurtz~\cite{EthKur}.
\proofend

\begin{Rem}
 By giving an example, we shed some light on why we consider only
 \ $(1+d)$-dimensional affine processes with state space \ $\RR_+\times\RR^d$ \ in Theorem \ref{Thm1}
 instead of \ $(n+d)$-dimensional ones with state space \ $\RR_+^n\times\RR^d$, $n\in\NN$.
\ Let \ $(Y_t)_{t\in\RR_+}$ \ be a two-dimensional continuous state branching process with infinitesimal generator
 \begin{align*}
  (\cA_Y f)(y)
  = \sum_{i=1}^2
     y_i \int_{\RR_+^2 \setminus\{0\} } \Bigl( f(y+u) - f(y) - f'_i(y) u_i \Bigr) p_i(\dd u),
 \end{align*}
 for \ $f\in\cC_c^2(\RR_+^2)$ \ and \ $y=(y_1,y_2)\in\RR_+^2$, \ where \ $p_i$, $i=1,2$, \  are
 \ $\sigma$-finite measures on \ $\RR_+^2\setminus\{0\}$ \ such that
 \begin{align}\label{help48}
   \int_{\RR_+^2\setminus\{0\}} (u_1 +  \Vert u\Vert^2) p_2(\dd u) < \infty
     \qquad \text{and}\qquad
   \int_{\RR_+^2\setminus\{0\}} (u_2 +  \Vert u\Vert^2) p_1(\dd u) < \infty,
 \end{align}
 see, e.g., Duffie et al. \cite[Theorem 2.7]{DufFilSch}.
Note that \ $Y$ \ can be considered as a two-dimensional affine process with state space \ $\RR_+^2$ \
 (formally with \ $d=0$).
\ Then, by a simple modification of Lemma \ref{Lem_inf_gen}, for all \ $\theta>0$, \ $f\in\cC_c^2(\RR_+^2)$ \ and \ $y=(y_1,y_2)\in\RR_+^2$,
 \begin{align*}
  (\cA_{Y_\theta}f)(y)
   &= \theta(\cA_Y f_\theta)(\theta y)\\
  &\quad= \theta
     \sum_{i=1}^2
      \theta y_i
       \int_{\RR_+^2\setminus\{0\}}
        \Bigl( f(\theta^{-1}(\theta y +u)) - f(\theta^{-1}\theta y)
               - \theta^{-1} f'_i(\theta^{-1} \theta y) u_i \Bigr)
        p_i(\dd u) \\
  &\quad= \theta^2
     \sum_{i=1}^2
      y_i \int_{\RR_+^2\setminus\{0\}}
           \Bigl( f(y + \theta^{-1}u) - f(y)
                  - \langle f'(y), \theta^{-1}u\rangle \Bigr)
           p_i(\dd u)\\
  &\phantom{\quad=\;}
      + \theta \sum_{i=1}^2
              y_i f'_{3-i}(y) \int_{\RR_+^2\setminus\{0\}} u_{3-i} \, p_i(\dd u),
 \end{align*}
 where the last equality follows by \eqref{help48}.
\ Supposing that \ $f$ \ is real-valued, by Taylor's theorem,
 \begin{align*}
   f(y + \theta^{-1}u) - f(y) - \langle f'(y), \theta^{-1}u\rangle
      = \frac{1}{2}\langle f''(y+\tau\theta^{-1}u) \theta^{-1} u, \theta^{-1}u\rangle
      = \frac{\theta^{-2}}{2}\langle f''(y+\tau\theta^{-1}u) u, u\rangle
 \end{align*}
 with some \ $\tau=\tau(u,y)\in[0,1]$.
\ Hence, similarly to the proof of \eqref{help4}, we get
 \begin{align*}
  \lim_{\theta\to\infty}
   \theta^2
   &\sum_{i=1}^2
    y_i \int_{\RR_+^2\setminus\{0\}}
         \Bigl( f(y + \theta^{-1}u) - f(y)
                - \langle f'(y), \theta^{-1}u\rangle \Bigr)
         p_i(\dd u)\\
  & = \frac{1}{2} \sum_{i=1}^2 y_i \int_{\RR_+^2\setminus\{0\}} \langle f''(y)u,u\rangle p_i(\dd u)
 \end{align*}
 for real-valued \ $f\in\cC_c^2(\RR_+^2)$ \ and \ $y=(y_1,y_2)\in\RR_+^2$.
\ However, \ $(\cA_{Y_\theta}f)(y)$ \ does not converge as \ $\theta \to \infty$
 \ provided that
 \[
   \sum_{i=1}^2
              y_i f'_{3-i}(y) \int_{\RR_+^2\setminus\{0\}} u_{3-i} \, p_i(\dd u) \ne 0.
 \]
We also note that this phenomena is somewhat similar to that of Remark 2.1 in Ma \cite{Ma}.\proofend
\end{Rem}

In the next remark we formulate some special cases of Theorem \ref{Thm1}.

\begin{Rem}\label{Rem1}
(i) If \ $(Y(t),X(t))_{t\in\RR_+}$ \ is a \ $(1+d)$-dimensional affine process on
 \ $\RR_+\times\RR^d$ \ with admissible parameters \ $(a,\alpha,b,0,m,\mu)$ \ such that condition
 \eqref{cond1} is satisfied, then the conditions of Theorem \ref{Thm1} are
 satisfied for
 \ $(Y^{(\theta)}(t),X^{(\theta)}(t))_{t\in\RR_+} := (Y(t),X(t))_{t\in\RR_+}$,
 \ $\theta\in\RR_{++}$, \ and hence
 \begin{align*}
  \left( \theta^{-1} Y(\theta t), \theta^{-1} X(\theta t)\right)_{t\in\RR_+}
  \distr (\cY(t),\cX(t))_{t\in\RR_+}
  \qquad \text{as \ $\theta\to\infty$}
 \end{align*}
 in \ $\DD(\RR_+,\RR_+\times\RR^d)$, \ where \ $(\cY(t),\cX(t))_{t\in\RR_+}$ \ is a
 \ $(1+d)$-dimensional affine process on \ $\RR_+\times\RR^d$ \ with admissible
 parameters \ $(0,\widetilde{\alpha},\widetilde{b},0,0,0)$, \ where
 \ $\widetilde{\alpha}$ \ and \ $\widetilde{b}$ \ are given in Theorem
 \ref{Thm1}.

(ii) If \ $(Y(t),X(t))_{t\in\RR_+}$ \ is a \ $(1+d)$-dimensional affine process on
 \ $\RR_+\times\RR^d$ \ with \ $(Y(0),X(0))=(0,0)$ \ and with admissible
 parameters \ $(0,\alpha,b,0,0,0)$, \ then
 \[
   \left(\theta^{-1} Y(\theta t), \theta^{-1} X(\theta t) \right)_{t\in\RR_+}
   \distre (Y(t),X(t))_{t\in\RR_+}
   \qquad \text{for all \ $\theta\in\RR_{++}$,}
 \]
 where \ $\distre$ \ denotes equality in distribution.
Indeed, by Proposition 1.6 on page 161 in Ethier and Kurtz \cite{EthKur}, it is
 enough to check that the semigroups (on the Banach space of bounded Borel measurable functions
 on \ $\RR_+\times\RR^d$) \ corresponding to the processes in question
 coincide.
By the definition of a core, this follows from the equality of the
 infinitesimal generators of the processes in question on the core
 \ $\cC_c^\infty(\RR_+\times\RR^d)$, \ which has been shown in the proof of
 Theorem \ref{Thm1}.\proofend
\end{Rem}

Next we present a corollary of Theorem \ref{Thm1} which states weak convergence
 of appropriately normalized one-dimensional continuous state branching processes with immigration.
Our corollary generalizes Theorem 2.3 in Huang et al.\ \cite{HuaMaZhu} in the
 sense that we do not have to suppose that
 \ $\int_1^\infty \xi^2\,m(\dd\xi)<\infty$, \ only that
 \ $\int_1^\infty \xi \,m(\dd\xi)<\infty$
 \ (with the notations of Huang et al.\ \cite{HuaMaZhu}), and our proof defers
 from that of Huang et al. \cite{HuaMaZhu}.

\begin{Cor}
For all \ $\theta\in\RR_{++}$, \ let \ $(Y^{(\theta)}(t))_{t\in\RR_+}$ \ be a
 one-dimensional continuous state branching process with immigration on \ $\RR_+$ \ with branching mechanism
 \[
   R^{(\theta)}(z)
   := \beta^{(\theta)}z + \alpha^{(\theta)}z^2
      + \int_{\RR_+} (\ee^{-zu} - 1 + zu)\,p(\dd u),
   \qquad z\in\RR_+,
 \]
 and with immigration mechanism
 \[
   F^{(\theta)}(z):= b^{(\theta)}z + \int_{\RR_+} (1- \ee^{-zu} )\,n(\dd u),
   \qquad z\in\RR_+,
 \]
 where \ $\alpha^{(\theta)}\geq 0$, \ $b^{(\theta)}\geq 0$,
 \ $\beta^{(\theta)}\in\RR$ \ and \ $n$ \ and \ $p$ \ are measures on
 \ $(0,\infty)$ \ such that
 \[
  \int_{\RR_+} u \, n(\dd u)<\infty
  \qquad \text{and}\qquad
  \int_{\RR_+} u^2 \, p(\dd u)<\infty.
 \]
Let \ $\alpha, b, \beta \in \RR$, \ and let
 \ $(Y(t))_{t\in\RR_+}$ \ be a one-dimensional
 continuous state branching process with immigration on \ $\RR_+$ \ with branching mechanism
 \[
   R(z):= -\beta z
          + \left(\alpha + \frac{1}{2}\int_{\RR_+} u^2\, p(\dd u)\right)z^2,
   \qquad z\in\RR_+,
 \]
 and with immigration mechanism
 \[
   F(z):= \left(b + \int_{\RR_+} u\, n(\dd u) \right)z,\qquad z\in\RR_+.
 \]
If
 \[
   \lim_{\theta\to\infty} \alpha^{(\theta)} = \alpha, \qquad
   \lim_{\theta\to\infty} b^{(\theta)} = b, \qquad
   \lim_{\theta\to\infty}\theta \beta^{(\theta)}=\beta, \qquad
   Y^{(\theta)}(0)\distr Y(0)
 \]
 as \ $\theta\to\infty$, \ then
 \begin{align*}
  \left(\theta^{-1} Y^{(\theta)}(\theta t)\right)_{t\in\RR_+}
  \distr (Y(t))_{t\in\RR_+}
  \qquad \text{as \ $\theta\to\infty$}
 \end{align*}
 in \ $\DD(\RR_+,\RR_+)$.
\end{Cor}

\noindent{\bf Proof.}
For each \ $\theta \in \RR_{++}$ \ and \ $t \in \RR_+$, \ let
 \ $X^{(\theta)}(t) := 0$.
\ Then for each \ $\theta \in \RR_{++}$, \ the process
 \ $(Y^{(\theta)}(t),X^{(\theta)}(t))_{t\in\RR_+}$ \ is a two-dimensional affine
 process with admissible parameters
 \ $(0,\overline{\alpha}^{(\theta)},\overline{b}^{(\theta)},
     \overline{\beta}^{(\theta)},m,\mu)$, \ where
 \begin{gather*}
  \overline{\alpha}^{(\theta)} := \begin{bmatrix}
                                 \alpha^{(\theta)} & 0 \\
                                 0 & 0 \\
                                \end{bmatrix} , \qquad
  \overline{b}^{(\theta)} := \begin{bmatrix}
                            b^{(\theta)} \\
                            0
                           \end{bmatrix} , \qquad
  \overline{\beta}^{(\theta)} := \begin{bmatrix}
                                \beta^{(\theta)} & 0 \\
                                0 & 0 \\
                               \end{bmatrix} , \\
  \mu(\dd\xi) = \mu(\dd\xi_1,\dd\xi_2)
  := p(\dd\xi_1)\times \delta_0(\dd\xi_2),\\
  m(\dd\xi) = m(\dd\xi_1,\dd\xi_2) := n(\dd\xi_1)\times \delta_0(\dd\xi_2),
 \end{gather*}
 where \ $\delta_0$ \ denotes the Dirac measure concentrated on \ $0\in\RR$.
\ Then, by Theorem \ref{Thm1}, for the two-dimensional affine processes
 \ $(Y^{(\theta)}(t),X^{(\theta)}(t))_{t\in\RR_+}$, \ $\theta \in \RR_{++}$, \ we have
 \begin{align*}
  \left( \theta^{-1} Y^{(\theta)}(\theta t), \theta^{-1} X^{(\theta)}(\theta t) \right)_{t\in\RR_+}
  \distr (Y(t),X(t))_{t\in\RR_+}
  \qquad \text{as \ $\theta\to\infty$}
 \end{align*}
 in \ $\DD(\RR_+,\RR_+\times\RR)$, \ where \ $(Y(t),X(t))_{t\in\RR_+}$ \ is a
 two-dimensional affine process on \ $\RR_+\times\RR$ \ with infinitesimal
 generator
 \begin{align*}
  (\cA_{(Y,X)} f)(x)
     = \left(\alpha + \frac{1}{2}\int_{\RR_+}u^2\,p(\dd u)\right)x_1f_{1,1}''(x)
       + \left(b + \beta x_1 + \int_{\RR_+}u\,n(\dd u)\right)f_1'(x),
 \end{align*}
 for \ $x=(x_1,x_2)\in\RR_+\times\RR$ \ and \ $f\in\cC^2_c(\RR_+\times\RR)$.
\ Note that in fact \ $X$ \ is the identically zero process.
Finally, Theorem 9.30 in Li \cite{Li} yields the assertion.
\proofend

\section{Least squares estimator for a critical two-dimensional
            affine diffusion process}\label{Section_statistics}

In this section continuous time stochastic processes will be written as \ $(\xi_t)_{t\in\RR_+}$ \ instead of \ $(\xi(t))_{t\in\RR_+}$.
\ Let \ $(\Omega,\cF,(\cF_t)_{t\in\RR_+},\PP)$ \ be a filtered probability space satisfying the usual conditions,
 i.e., \ $(\Omega,\cF,\PP)$ \ is complete, the filtration $(\cF_t)_{t\in\RR_+}$ \ is right-continuous and
 \ $\cF_0$ \ contains all the \ $\PP$-null sets in \ $\cF$.
\ Let \ $(W_t)_{t\in\RR_+}$ \ and \ $(B_t)_{t\in\RR_+}$ \ be independent standard \ $(\cF_t)_{t\in\RR_+}$-Wiener processes.
Let us consider the following two-dimensional diffusion process given by the SDE
 \begin{align}\label{2dim_affine}
  \begin{cases}
   \dd Y_t = (a-bY_t)\,\dd t + \sqrt{Y_t}\,\dd W_t,\\
   \dd X_t = (m-\theta X_t)\,\dd t + \sqrt{Y_t}\,\dd B_t,
  \end{cases} \qquad t\in\RR_+,
 \end{align}
 where \ $a\in\RR_{++}$ \ and \ $b,\theta,m\in\RR$.

\subsection{Preparations and (sub)(super)criticality}

The next proposition is about the existence and uniqueness of a strong solution
 of the SDE \eqref{2dim_affine}.

\begin{Pro}\label{Pro_affine}
Let \ $(\eta,\zeta)$ \ be a random vector independent of \ $(W_t,B_t)_{t\in\RR_+}$ \ satisfying \ $\PP(\eta\geq 0)=1$.
\ Then, for all \ $a\in\RR_{++}$ \ and \ $b, m, \theta\in\RR$, \ there is a (pathwise) unique strong
 solution \ $(Y_t,X_t)_{t\in\RR_+}$ \ of the SDE \eqref{2dim_affine} such that \ $\PP((Y_0,X_0) = (\eta,\zeta))=1$
 \ and \ $\PP(\text{$Y_t\geq 0$ \ for all \ $t\in\RR_+$})=1$.
\ Further, for all \ $0\leq s\leq t$,
 \begin{align}\label{help27}
  Y_t = \ee^{-b(t-s)}
        \left( Y_s + a\int_s^t \ee^{-b(s-u)}\,\dd u
               + \int_s^t \ee^{-b(s-u)}\sqrt{Y_u} \,\dd W_u\right),
 \end{align}
 and
 \begin{align}\label{help22}
  X_t = \ee^{-\theta(t-s)}
        \left( X_s + m\int_s^t \ee^{-\theta(s-u)}\,\dd u
               + \int_s^t \ee^{-\theta(s-u)}\sqrt{Y_u} \,\dd B_u\right).
 \end{align}
\end{Pro}

\noindent{\bf Proof.}
By Ikeda and Watanabe \cite[Example 8.2, page 221]{IkeWat}, there is a pathwise
 unique non-negative strong solution \ $(Y_t)_{t\in\RR_+}$ \ of the first
 equation in \eqref{2dim_affine} with any initial value \ $\eta$ \ satisfying
 \ $\PP(\eta\geq 0)=1$.
\ Next, by applications of the It\^{o}'s formula to the processes
 \ $(\ee^{b t} Y_t)_{t\in\RR_+}$ \ and \ $(\ee^{\theta t} X_t)_{t\in\RR_+}$,
 \ respectively, we have
 \begin{align*}
  \dd(\ee^{b t} Y_t)
  &= b\ee^{b t}Y_t\,\dd t + \ee^{b t}\dd Y_t
   = b \ee^{b t}Y_t\,\dd t
     + \ee^{b t}\big( (a-b Y_t)\,\dd t + \sqrt{Y_t}\,\dd W_t\big)\\
  &= a\ee^{b t}\,\dd t + \ee^{b t}\sqrt{Y_t}\,\dd W_t,
  \qquad t\in\RR_+,
 \end{align*}
 and
 \begin{align*}
  \dd(\ee^{\theta t} X_t)
  &= \theta\ee^{\theta t}X_t\,\dd t + \ee^{\theta t}\dd X_t
   = \theta\ee^{\theta t}X_t\,\dd t
     + \ee^{\theta t}\big( (m-\theta X_t)\,\dd t + \sqrt{Y_t}\,\dd B_t\big)\\
  &= m\ee^{\theta t}\,\dd t + \ee^{\theta t}\sqrt{Y_t}\,\dd B_t,
  \qquad t\in\RR_+,
 \end{align*}
 which imply \eqref{help27} and \eqref{help22} in case of \ $s=0$.
\ If \ $0\leq s\leq t$, \ then, by
 \[
   \ee^{b t} Y_t
   = \ee^{b s} Y_s + a\int_s^t \ee^{b u}\,\dd u
     + \int_s^t \ee^{b u}\sqrt{Y_u} \,\dd B_u,
 \]
 and
 \[
   \ee^{\theta t} X_t
   = \ee^{\theta s} X_s + m\int_s^t \ee^{\theta u}\,\dd u
     + \int_s^t \ee^{\theta u}\sqrt{Y_u} \,\dd B_u,
 \]
 we have \eqref{help27} and \eqref{help22}.
Finally, we note that the existence of a pathwise unique strong solution
 \ $(Y_t,X_t)_{t\in\RR_+}$ \ of the SDE \eqref{2dim_affine} with
 \ $\PP(\text{$Y_t\geq 0$ \ for all \ $t\in\RR_+$})=1$ \ follows also by a
 general result of Dawson and Li \cite[Theorem 6.2]{DawLi}.
\proofend

Note that it is the assumption \ $a\in\RR_{++}$ \ that ensures \ $\PP(\text{$Y_t\,{\geq}\, 0$,
\ $\forall\,t\,{\in}\,\RR_+$})\,{=}\,1$.

Next we present a result about the first moment of \ $(Y_t,X_t)_{t\in\RR_+}$.

\begin{Pro}\label{Pro_moments}
Let \ $(Y_t,X_t)_{t\in\RR_+}$ \ be a strong solution of the SDE \eqref{2dim_affine}
 satisfying \ $\PP(Y_0\geq 0)=1$, \ $\EE(Y_0)<\infty$,
 \ and \ $\EE(X_0)<\infty$.
\ Then
 \begin{align*}
  \begin{bmatrix}
   \EE(Y_t) \\
   \EE(X_t) \\
  \end{bmatrix}
  = \begin{bmatrix}
     \ee^{-bt} & 0 \\
     0 & \ee^{-\theta t} \\
    \end{bmatrix}
    \begin{bmatrix}
     \EE(Y_0) \\
     \EE(X_0) \\
    \end{bmatrix}
    + \begin{bmatrix}
       \int_0^t \ee^{-bs} \, \dd s & 0 \\
       0 & \int_0^t \ee^{-\theta s} \, \dd s \\
      \end{bmatrix}
      \begin{bmatrix}
       a \\
       m \\
      \end{bmatrix},\qquad t\in\RR_+,
 \end{align*}
\end{Pro}

\noindent{\bf Proof.}
By Proposition \ref{Pro_affine}, we have
 \begin{align*}
   &Y_t = \ee^{-bt}
         \left( Y_0 + a\int_0^t\ee^{bu}\,\dd u
                + \int_0^t \ee^{bu} \sqrt{Y_u}\,\dd W_u\right),
         \qquad t\in\RR_+,\\
   &X_t = \ee^{-\theta t}
         \left( X_0 + m\int_0^t\ee^{\theta u}\,\dd u
                + \int_0^t \ee^{\theta u} \sqrt{Y_u}\,\dd B_u\right),
         \qquad t\in\RR_+,
 \end{align*}
 and so, taking expectations of both sides,
 \begin{align*}
  & \EE(Y_t)
    = \ee^{-bt}\EE(Y_0) + a \ee^{-bt} \int_0^t\ee^{bu}\,\dd u
    = \ee^{-bt}\EE(Y_0) + a \int_0^t\ee^{-bu}\,\dd u,
  \quad t\in\RR_+,\\
  & \EE(X_t)
    = \ee^{-\theta t}\EE(X_0) + m \ee^{-\theta t} \int_0^t\ee^{\theta u}\,\dd u
    = \ee^{-\theta t}\EE(X_0) + m \int_0^t\ee^{-\theta u}\,\dd u,
   \;\; t\in\RR_+,
 \end{align*}
 where we used that the processes
 \[
    \left( \int_0^t \ee^{bu} \sqrt{Y_u}\,\dd W_u \right)_{t\in\RR_+}
     \quad \text{and}\qquad
     \left( \int_0^t \ee^{\theta u} \sqrt{Y_u}\,\dd B_u \right)_{t\in\RR_+}
 \]
 are martingales which can be checked as follows.
First we check that they are local martingales with respect to the filtration \ $(\cF_t)_{t\in\RR_+}$.
\ Let us define the increasing sequence of stopping times by \ $\delta_n:=\inf\{t\geq 0:Y_t\geq n\}$, $n\in\NN$.
\ Since \ $Y$ \ has continuous trajectories almost surely, we have \ $\PP(\lim_{n\to\infty}\delta_n=\infty)=1$.
\ Using \ $(\delta_n)_{n\in\NN}$ \ as a localizing sequence, we have
 \[
   \EE\left(\int_0^{t\wedge \delta_n} \ee^{2bu} Y_u\,\dd u \right)
      \leq nt\max(1,\ee^{2bt}), \quad t\in\RR_+,\; n\in\NN.
 \]
The local martingale property of
 \ $\bigl( \int_0^t \ee^{bu} \sqrt{Y_u}\,\dd W_u \bigr)_{t\in\RR_+}$ \
 follows by Ikeda and Watanabe \cite[page 57]{IkeWat}.
Hence, using \eqref{help27} and that \ $a\in\RR_{++}$, \ we find that
 \[
  \EE(\ee^{b(t\wedge \delta_n)}Y_{t\wedge \delta_n})
     = \EE(Y_0) + a \EE\left( \int_0^{t\wedge \delta_n} \ee^{bu}\,\dd u \right)
     \leq \EE(Y_0) + at\max(1,\ee^{bt})
 \]
 for all \ $t\in\RR_+$ \ and \ $n\in\NN$, \ and then, by Fatou's lemma,
 \begin{align}\label{help42}
    \EE(\ee^{bt}Y_t)
      \leq \liminf_{n\to\infty} \EE(\ee^{b(t\wedge \delta_n)}Y_{t\wedge \delta_n})
      \leq \EE(Y_0) + at\max(1,\ee^{bt}),
     \qquad t\in\RR_+.
 \end{align}
 Next, we can deduce that \ $\bigl( \int_0^t \ee^{bu} \sqrt{Y_u}\,\dd W_u
 \bigr)_{t\in\RR_+}$ \ is indeed a martingale.\break
First, we note that a local martingale \ $M$ \ is a square integrable martingale if \ $E([M,M]_t)<\infty$ \ for all
 \ $t\in\RR_+$, \ where \ $([M,M]_t)_{t\in\RR_+}$ \  denotes the quadratic variation process of \ $M$, \
 see, e.g., Corollary 3 on page 73 in Protter \cite{Pro}.
Here the quadratic variation process
of \ $\bigl( \int_0^t \ee^{bu} \sqrt{Y_u}\,\dd W_u \bigr)_{t\in\RR_+}$ \
 takes the form
 \[
    E\left(\int_0^t \ee^{2bu} Y_u\,\dd u \right)
     <\infty, \qquad  t\in\RR_{+},
 \]
 where, for the inequality, we used Fubini's theorem, \eqref{help42} and our assumption \ $\EE(Y_0)<\infty$.
Replacing \ $b$ \ by \ $\theta$, \ we have the desired martingale property of
 \ $\bigl( \int_0^t \ee^{\theta u} \sqrt{Y_u}\,\dd B_u \bigr)_{t\in\RR_+}$, \ too.
\proofend

Next we show that the process \ $(Y_t,X_t)_{t\in\RR_+}$ \ given by the SDE \eqref{2dim_affine}
 is an affine process.

\begin{Pro}\label{Pro_affine_2}
Let \ $(Y_t,X_t)_{t\in\RR_+}$ \ be a strong solution of the SDE \eqref{2dim_affine}
  satisfying \ $\PP(Y_0\geq 0)=1$.
\ Then \ $(Y_t,X_t)_{t\in\RR_+}$ \ is an affine process with infinitesimal generator
\begin{align}\label{help50_affine_generator}
 \begin{split}
  (\cA_{(Y,X)} f)(x)
    = (a-bx_1)f_1'(x) + (m-\theta x_2)f_2'(x)
       + \frac{1}{2}x_1(f_{1,1}''(x) + f_{2,2}''(x))
  \end{split}
 \end{align}
 for \ $x=(x_1,x_2)\in\RR_+\times\RR$ \ and \ $f\in\cC^2_c(\RR_+\times\RR)$.
\end{Pro}

\noindent{\bf Proof.}
For calculating the infinitesimal generator of \ $(Y_t,X_t)_{t\in\RR_+}$, \ without loss of generality,
 we may suppose that \ $\PP((Y_0,X_0) = (y_0,x_0))=1$, \ where \ $(y_0,x_0)\in\RR_+\times\RR$.
\ By It\^{o}'s formula, for all real-valued functions \ $f\in\cC^2_c(\RR_+\times\RR)$ \ we have
 \begin{align*}
  f(Y_t,X_t) & = f(y_0,x_0)
                + \int_0^t f_1'(Y_s,X_s)\sqrt{Y_s} \,\dd W_s
                + \int_0^t f_2'(Y_s,X_s)\sqrt{Y_s} \,\dd B_s \\
            &\phantom{=\;} + \int_0^t f_1'(Y_s,X_s)(a-bY_s)\,\dd s
                 + \int_0^t f_2'(Y_s,X_s)(m-\theta X_s)\,\dd s \\
            &\phantom{=\;} + \frac{1}{2}\left(\int_0^t f_{1,1}''(Y_s,X_s)Y_s\,\dd s
                                   + \int_0^t f_{2,2}''(Y_s,X_s)Y_s\,\dd s \right) \\
           &= f(y_0,x_0) + \int_0^t (\cA_{(Y,X)} f)(Y_s,X_s)\,\dd s
             + M_t(f),\qquad t\in\RR_+,
 \end{align*}
 where
 \[
    M_t(f):= \int_0^t f_1'(Y_s,X_s)\sqrt{Y_s} \,\dd W_s
               + \int_0^t f_2'(Y_s,X_s)\sqrt{Y_s} \,\dd B_s,\qquad t\in\RR_+,
 \]
 and \ $\cA_{(Y,X)} f$ \ is given by \eqref{help50_affine_generator}.
Here \ $(M_t(f))_{t\in\RR_+}$ \ is a square integrable martingale
 with respect to the filtration \ $(\cF_t)_{t\in\RR_+}$, \ since
 \begin{align*}
   &\int_0^t \EE((f_1'(Y_s,X_s))^2Y_s) \,\dd s
       \leq C_1    \int_0^t \EE(Y_s) \,\dd s
       <\infty,\qquad t\in\RR_+,\\
   &   \int_0^t \EE((f_2'(Y_s,X_s))^2Y_s) \,\dd s
       \leq C_2    \int_0^t \EE(Y_s) \,\dd s
       <\infty,\qquad t\in\RR_+,
 \end{align*}
 with some constants \ $C_1>0$ \ and \ $C_2>0$, \ where the finiteness of the integrals
 follow by Proposition \ref{Pro_moments}.
Finally, if \ $f\in\cC^2_c(\RR_+\times\RR)$ \ is complex valued, then, by decomposing \ $f$ \ into
 real and imaginary parts, one can argue in the same way as above.
\proofend

By Proposition \ref{Pro_affine_2}, the process \ $(Y_t,X_t)_{t\in\RR_+}$ \
 given by \eqref{2dim_affine} is a two-dimen\-sional affine process with admissible parameters
  \[
   \left( \begin{bmatrix}
           0 & 0 \\
           0 & 0
          \end{bmatrix} , \,
          \frac{1}{2} \begin{bmatrix}
                       1 & 0 \\
                       0 & 1
                      \end{bmatrix} , \,
          \begin{bmatrix}
           a \\
           m
          \end{bmatrix} , \,
          \begin{bmatrix}
             -b & 0 \\
             0 & -\theta \\
            \end{bmatrix} , \,
          0 ,
          0 \right) .
 \]

In what follows we define and study criticality of the affine process given by the SDE \eqref{2dim_affine}.

\begin{Def}\label{Def_criticality}
Let \ $(Y_t,X_t)_{t\in\RR_+}$ \ be an affine diffusion process given by the SDE \eqref{2dim_affine}
 satisfying \ $\PP(Y_0\geq 0)=1$.
\ We call \ $(Y_t,X_t)_{t\in\RR_+}$ \ subcritical, critical or supercritical if the spectral radius of the matrix
 \[
   \begin{pmatrix}
    \ee^{-bt} & 0 \\
    0 & \ee^{-\theta t} \\
   \end{pmatrix}
 \]
 is less than \ $1$, \ equal to \ $1$ \ or greater than \ $1$, \ respectively.
\end{Def}

Note that, since the spectral radius of the matrix given in Definition \ref{Def_criticality}
 is \ $\max(\ee^{-bt},\ee^{-\theta t})$, \ the affine process given in Definition \ref{Def_criticality}
 is
 \begin{align*}
  \text{subcritical} \qquad & \text {if \ $b>0$ \ and \ $\theta>0$,}\\
  \text{critical} \qquad
  & \text{if \ $b=0$, \ $\theta\geq 0$ \ or \ $b\geq 0$, \ $\theta=0$,}\\
  \text{supercritical} \qquad & \text{if \ $b<0$ \ or \ $\theta<0$.}
 \end{align*}
Definition \ref{Def_criticality} of criticality is in accordance with the
 corresponding definition for one-dimensional continuous state branching processes, see, e.g.,
 Li \cite[page~58]{Li}.

In this section we will always suppose that
 \begin{align*}
 \textbf{Condition (C): } &\quad  (b,\theta)=(0,0),\; \PP(Y_0 \geq 0)=1,\\
                          &\quad  \text{$\EE(Y_0) < \infty$, \ and \ $\EE(X_0^2) < \infty$.}
 \end{align*}
For some explanations why we study only this special case, see Remarks \ref{Rem2}, \ref{Rem3} and \ref{Rem6}.
In the next sections under Condition (C) we will study asymptotic behaviour of least squares
 estimator of \ $\theta$ \ and \ $(\theta,m)$, \ respectively.
Before doing so we recall some critical models both in discrete and continuous time.

In general, parameter estimation for critical models has a long history.
A common feature of the estimators for parameters of critical models is that one may prove weak limit
 theorems for them by using norming factors that are usually different from the norming factors
 for the subcritical and supercritical models.
Further, it may happen that one has to use different norming factors for two different critical cases.

We recall some discrete time critical models.
If \ $(\xi_k)_{k\in\ZZ_+}$ \ is an AR(1) process, i.e.,
 \ $\xi_k = \varrho \xi_{k-1} +\zeta_k$, \ $k\in\NN$, \ with \ $\xi_0 = 0$ \ and
 an i.i.d.\ sequence \ $(\zeta_k)_{k\in\NN}$ \ having mean \ $0$ \ and positive
 variance, then the (ordinary) least squares estimator of the so-called
 stability  parameter \ $\varrho$ \ based on the sample \ $\xi_1,\ldots,\xi_n$
 \ takes the form
 \[
   \widetilde\varrho_n
   = \frac{\sum_{k=1}^n \xi_{k-1}\xi_k}{\sum_{k=1}^n \xi_k^2}, \qquad n\in\NN,
 \]
 see, e.g., Hamilton \cite[17.4.2]{Ham}.
In the critical case, i.e., when \ $\varrho=1$, \ by Hamilton \cite[17.4.7]{Ham},
 \[
   n(\widetilde\varrho_n-1)
   \distr \frac{\int_0^1 \cW_t \, \dd \cW_t}{\int_0^1 \cW_t^2 \, \dd t}
   \qquad \text{as \ $n\to\infty$,}
 \]
 where \ $(\cW_t)_{t\in\RR_+}$ \ is a standard Wiener process and \ $\distr$
 \ denotes convergence in distribution.
Here \ $n(\widetilde\varrho_n-1)$ \ is known as the Dickey-Fuller statistics.
We emphasize that the asymptotic behaviour of \ $\widetilde\varrho_n$ \ is completely different in the subcritical
 \ $(\vert\rho\vert<1)$ \ and supercritical \ $(\vert\rho\vert>1)$ \ cases, where it is asymptotically normal and asymptotically
 Cauchy, respectively, see, e.g., Mann and Wald \cite{ManWal}, Anderson \cite{And} and White \cite{Whi}.

For continuous time critical models, we recall that Huang et al.\ \cite[Theorem 2.4]{HuaMaZhu} studied
 asymptotic behaviour of weighted conditional least squares estimator of the
 drift parameters for discretely observed continuous time critical branching processes
 with immigration given by
 \begin{align*}
  \widetilde Y_t & = \widetilde Y_0 + \int_0^t(a+b \widetilde Y_s)\,\dd s
            + \sigma \int_0^t \sqrt{\widetilde Y_s}\,\dd \cW_s
            + \int_0^t \int_{[0,\infty)}\xi\,\cN_0(\dd s,\dd \xi)\\
      &\phantom{=\;}
            + \int_0^t \int_0^{\widetilde Y_{s-}}\int_{[0,\infty)}\xi\,(\cN_1(\dd s,\dd u,\dd \xi) - \dd s\,\dd u\,p(\dd\xi)),
            \qquad t\in\RR_+,
 \end{align*}
 where \ $\widetilde Y_0\geq 0$, \ $a\geq 0$, \ $b\in\RR$, $\sigma\geq 0$, \ $\cW$ \ is a standard Wiener process,
 \ $\cN_0(\dd s,\dd\xi)$ \ is a Poisson random measure on \ $(0,\infty)\times[0,\infty)$ \ with intensity
 \ $\dd s\,n(\dd\xi)$, \ $\cN_1(\dd s,\dd u,\dd\xi)$ \ is a Poisson random measure on
 \ $(0,\infty)\times (0,\infty) \times [0,\infty)$ \ with intensity \ $\dd s\,\dd u\,p(\dd\xi)$ \ such that
 the \ $\sigma$-finite measures \ $n$ \ and \ $p$ \ are supported by \ $(0,\infty)$ \ and
 \[
    \int_0^\infty \xi\, n(\dd\xi) + \int_0^\infty \xi\wedge\xi^2\,p(\dd\xi) <\infty.
 \]
Our technique differs from that of Huang et al.\ \cite{HuaMaZhu} and for
 completeness we note that the limit distribution and some parts of the proof
 of their Theorem 2.4 suffer from some misprints.
Furthermore, Hu and Long \cite{HuLon3} studied the problem of parameter estimation for
 critical mean-reverting $\alpha$-stable motions
 \[
   \dd \widetilde X_t = (m-\theta \widetilde X_t)\,\dd t + \dd Z_t,\qquad t\in\RR_+,
 \]
 where \ $Z$
 \ is an \ $\alpha$-stable L\'evy motion with \ $\alpha\in(0,2)$)
 \ observed at discrete instants.
A least squares estimator is obtained and its asymptotics is discussed in the
 singular case \ $(m,\theta)=(0,0)$.
\ We note that the forms of the limit distributions of least squares estimators for critical
 two-dimensional affine diffusion processes in our Theorems \ref{Thm2}
 and \ref{Thm3} are the same as that of the limit distributions in Theorems 3.2
 and 4.1 in Hu and Long \cite{HuLon3}, respectively.
We also recall that Hu and Long \cite{HuLon1} considered the problem of parameter estimation not only for
 critical mean-reverting $\alpha$-stable motions, but also for some subcritical ones \ ($m=0$ \ and \ $\theta>0$)
 \ by proving limit theorems for the least squares estimators that are completely different from the ones in the critical case.
\ Huang et al.\ \cite{HuaMaZhu} investigated the asymptotic behaviour of weighted conditional least squares estimator of the
 drift parameters not only for critical continuous time branching processes with immigration, but also for subcritical and supercritical ones.

Using our scaling Theorem \ref{Thm1} we can only handle a special critical affine diffusion model given by \eqref{DD}
 (for a more detailed discussion, see Remark \ref{Rem3}).
The other critical and non-critical cases are under investigation but different techniques are needed.

\subsection{Least squares estimator of \ $\theta$ \ when \ $m$ \ is known}
\label{section_LSE_theta}

The least squares estimator (LSE) of \ $\theta$ \ based on the observations
 \ $X_i$, $i=0,1,\ldots,n$, \ can be obtained by solving the extremum problem
 \[
   \widetilde\theta_n^{\mathrm{LSE}}
   := \argmin_{\theta\in\RR} \sum_{i=1}^n (X_i - X_{i-1} - (m-\theta X_{i-1}))^2.
 \]
This definition of LSE of \ $\theta$ \ can be considered as the counterpart of the one
 given in Hu and Long \cite[formula (1.2)]{HuLon1} for generalized Ornstein-Uhlenbeck processes
 driven by \ $\alpha$-stable motions, see also Hu and Long \cite[formulas (3.1) and (4.1)]{HuLon3}.
For a mathematical motivation of the definition of the LSE of \ $\theta$, \ see later on Remark \ref{Rem5}.
With the notation
 \ $f(\theta):=\sum_{i=1}^n (X_i - X_{i-1} - (m-\theta X_{i-1}))^2$,
 \ $\theta\in\RR$, \ the equation \ $f'(\theta)=0$ \ takes the form:
 \[
   2 \sum_{i=1}^n (X_i - X_{i-1} - (m-\theta X_{i-1})) X_{i-1} = 0.
 \]
Hence
 \[
   \left(\sum_{i=1}^n X_{i-1}^2\right)\theta
   = - \sum_{i=1}^n (X_i - X_{i-1} - m) X_{i-1},
 \]
 i.e.,
 \begin{align}\label{help11}
 \begin{split}
  \widetilde\theta_n^{\mathrm{LSE}}
  = - \frac{\sum_{i=1}^n (X_i - X_{i-1} - m) X_{i-1}}{\sum_{i=1}^n X_{i-1}^2}
  = - \frac{\sum_{i=1}^n (X_i - X_{i-1})X_{i-1} - \left(\sum_{i=1}^n X_{i-1}\right)m}
           {\sum_{i=1}^n X_{i-1}^2}
 \end{split}
 \end{align}
 provided that \ $\sum_{i=1}^n X_{i-1}^2>0$.
\ Since \ $f''(\theta)=2\sum_{i=1}^n X_{i-1}^2$, \ $\theta\in\RR$, \ we have
 \ $\widetilde\theta_n^{\mathrm{LSE}}$ \ is indeed the solution of the extremum
 problem provided that \ $\sum_{i=1}^n X_{i-1}^2>0$.

\begin{Thm}\label{Thm2}
Let us assume that Condition (C) holds.
Then \ $\PP(\sum_{i=1}^n X_{i-1}^2>0)=1$ \ for all \ $n\geq2$, \ and there exists
 a unique LSE \ $\widetilde\theta_n^{\mathrm{LSE}}$ \ which has the form given in
 \eqref{help11}.
Further,
 \begin{align}\label{help23}
  n\widetilde\theta_n^{\mathrm{LSE}}
  \distr
  - \frac{\int_0^1 \cX_t\,\dd \cX_t - m\int_0^1 \cX_t\,\dd t}
         {\int_0^1 \cX_t^2\,\dd t}
  \qquad \text{as \ $n\to\infty$,}
 \end{align}
 where \ $(\cX_t)_{t\in\RR_+}$ \ is the second coordinate of a two-dimensional
 affine process \ $(\cY_t,\cX_t)_{t\in\RR_+}$ \ given by the unique strong solution of the SDE
 \begin{align}\label{help16}
 \begin{split}
 \begin{cases}
   \dd\cY_t = a\,\dd t + \sqrt{\cY_t}\,\dd \cW_t,\\[2mm]
   \dd\cX_t = m\,\dd t + \sqrt{\cY_t}\,\dd \cB_t,
 \end{cases}\qquad t\in\RR_+,
 \end{split}
 \end{align}
 with initial value \ $(\cY_0,\cX_0)=(0,0)$, \
 where \ $(\cW_t)_{t\in\RR_+}$ \ and \ $(\cB_t)_{t\in\RR_+}$ \ are independent
 standard Wiener processes.
\end{Thm}

\begin{Rem}\label{Rem2}
(i)  The limit distributions in Theorem \ref{Thm2} have the same forms as those
 of the limit distributions in Theorem 3.2 in Hu and Long \cite{HuLon3}.

\noindent
(ii) The limit distribution of \ $n\widetilde\theta_n^{\mathrm{LSE}}$ \ as
 \ $n\to\infty$ \ in Theorem \ref{Thm2} can be written also in the form
 \[
   - \frac{\int_0^1 \cX_t\,\dd (\cX_t - mt)}{\int_0^1 \cX_t^2\,\dd t}
   = - \frac{\int_0^1 \cX_t\sqrt{\cY_t}\,\dd \cB_t}
            {\int_0^1 \cX_t^2\,\dd t}.
 \]
(iii) By Proposition \ref{Pro_affine_2}, the affine process \ $(\cY_t,\cX_t)_{t\in\RR_+}$ \ given in
 Theorem \ref{Thm2} has infinitesimal
 generator
 \begin{align*}
  (\cA_{(\cY,\cX)} f)(x)
   = \frac{1}{2}x_1f_{1,1}''(x) + \frac{1}{2}x_1f_{2,2}''(x) + af_1'(x) + mf_2'(x)
 \end{align*}
 where \ $x=(x_1,x_2)\in\RR_+\times\RR$ \ and \ $f\in\cC^2_c(\RR_+\times\RR)$.

\noindent(iv)
Under the Condition (C), by Theorem \ref{Thm2} and Slutsky's lemma, we get
 \ $\widetilde\theta_n^{\mathrm{LSE}}$ \ converges stochastically to the parameter \ $\theta=0$ \ as \ $n\to\infty$.\proofend
\end{Rem}

\noindent{\bf Proof of Theorem \ref{Thm2}.}
By \eqref{help22}, we have
 \[
   X_t = X_0 + m t + \int_0^t \sqrt{Y_s} \, \dd B_s , \qquad t \in \RR_+ .
 \]
Hence for all \ $t \in \RR_{++}$, \ the conditional distribution of
 \ $X_t$ \ given \ $X_0$ \ and \ $(Y_s)_{s\in[0,t]}$ \ is a normal distribution
 with mean \ $X_0 + mt$ \ and with variance \ $\int_0^t Y_s \, \dd s$.
\ Here the variance \ $\int_0^t Y_s \, \dd s$ \ is positive almost surely for
 all \ $t \in \RR_{++}$.
\ Indeed, let
 \ $A_t := \{ \omega \in \Omega
              : \text{$s \mapsto Y_s(\omega)$ \ is continuous on \ $[0, t]$} \}$.
\ Then \ $\PP(A_t) = 1$, \ and, since \ $\PP(Y_0\geq 0)=1$, \ for all \ $\omega \in A_t$,
 \ $\int_0^t Y_s(\omega) \, \dd s = 0$ \ if and only if
 \ $Y_s(\omega) = 0$ \ for all \ $s \in [0, t]$.
\ By \eqref{2dim_affine}, we have
 \[
   Y_s = Y_0 + a s + \int_0^s \sqrt{Y_u} \, \dd W_u , \qquad s \in \RR_+ .
 \]
The stochastic integral on the right hand side can be approximated as
 \[
   \sup_{s\in[0,t]}
    \left| \sum_{i=1}^\ns \sqrt{Y_{(i-1)/n}} (W_{i/n} - W_{(i-1)/n})
           - \int_0^s Y_u \, \dd W_u \right|
   \stoch 0 \qquad \text{as \ $n \to \infty$}
 \]
 for all \ $t\in\RR_+$, \ by Jacod and Shiryaev \cite[Theorem I.4.44]{JSh}.
Hence there exists a sequence \ $(n_k)_{k\in\NN}$ \ of positive integers such
 that
 \[
   \sup_{s\in[0,t]}
    \left| \sum_{i=1}^{\lfloor n_k t \rfloor}
            \sqrt{Y_{(i-1)/n_k}} (W_{i/n_k} - W_{(i-1)/n_k})
           - \int_0^s Y_u \, \dd W_u \right|
   \as 0 \qquad \text{as \ $k \to \infty$}
 \]
 for all \ $t\in\RR_+$.
\ Consequently, with the notation
 \[
  \widetilde A_t := \left\{ \omega \in \Omega : \int_0^t Y_s(\omega) \, \dd s = 0 \right\},
 \]
 we have
 \begin{align*}
  \widetilde A_t\cap A_t
    &\subset \left\{\widetilde A_t\bigcap\left\{\text{$\int_0^s Y_u \, \dd W_u = 0$ \ for all
                       \ $s \in [0, t]$}\right\}\right\} \\
  &\subset \left\{\widetilde A_t\cap \left\{\text{$Y_s = Y_0 + a s$ \ for all \ $s \in [0, t]$}\right\}\right\}\\
  &=\left\{\widetilde A_t\cap \left\{\text{$Y_0s + a s^2/2=0$ \ for all \ $s \in [0, t]$}\right\}\right\}\\
  &= \left\{\widetilde A_t\cap \left\{\text{$Y_0 = -a s/2$ \ for all \ $s \in [0, t]$}\right\}\right\}
   = \emptyset,
 \end{align*}
 implying \ $\PP\bigl(\int_0^t Y_s \, \dd s = 0\bigr) = 0$, \ and hence
 \ $\PP\bigl(\int_0^t Y_s \, \dd s > 0\bigr) = 1$.
\ It yields that
 \begin{align}\label{help25}
  \PP(X_t = 0) = \EE\big(\PP(X_t = 0 \mid X_0, (Y_s)_{s\in[0,t]})\big) = 0, \qquad
  t \in \RR_{++} ,
 \end{align}
 and hence \ $\PP(\sum_{i=1}^n X_{i-1}^2 > 0) = 1$ \ for all \ $n \geq 2$.

Now we turn to prove \eqref{help23}.
By It\^{o}'s formula, we have
 \ $\dd(\cX_t^2) = 2 \cX_t \dd \cX_t + \cY_t \, \dd t$, \ $t \in \RR_+$,
 \ and hence, using also \ $\cX_0 = 0$, \ we have
 \begin{align}\label{help49}
  \int_0^1 \cX_s \, \dd \cX_s
  = \frac{1}{2} \left( \cX_1^2 - \int_0^1 \cY_s \, \dd s \right) .
 \end{align}
For the process \ $(X_t)_{t\in\RR_+}$, \ a discrete version
 \begin{align}\label{help53}
   \sum_{i=1}^n (X_i - X_{i-1}) X_{i-1}
   = \frac{1}{2} \left( X_n^2 - X_0^2 - \sum_{i=1}^n (X_i - X_{i-1})^2 \right)
 \end{align}
 of the identity \eqref{help49} can be easily checked.
The aim of the following discussion is to prove
 \begin{align}\label{help10}
  \begin{split}
   &\biggl(\frac{1}{n^2} \sum_{i=1}^n X_{i-1} ,
          \frac{1}{n^3} \sum_{i=1}^n X_{i-1}^2 ,
           \frac{1}{n} X_n ,
           \frac{1}{n} X_0 ,
           \frac{1}{n^2} \sum_{i=1}^n (X_i - X_{i-1})^2 \biggr) \\
   &\qquad\distr
    \left( \int_0^1 \cX_t \, \dd t , \int_0^1 \cX_t^2 \, \dd t ,
           \cX_1 , 0, \int_0^1 \cY_t \, \dd t \right)
   \qquad \text{as \ $n\to\infty$.}
  \end{split}
 \end{align}
Let us consider the unique strong solution of the SDE
 \begin{align}\label{help16_b}
  \begin{split}
   \begin{cases}
    \dd \tY_t = a \, \dd t + \sqrt{\tY_t} \, \dd W_t , \\[2mm]
    \dd \tX_t = m \, \dd t + \sqrt{\tY_t} \, \dd B_t ,
   \end{cases}\qquad t \in \RR_+ ,
  \end{split}
 \end{align}
 with initial value \ $(\tY_0, \tX_0)=(0, 0)$, \ where
 \ $(W_t)_{t\in\RR_+}$ \ and \ $(B_t)_{t\in\RR_+}$ \ are the independent
 standard \ $(\cF_t)_{t\in\RR_+}$-Wiener processes appearing in the SDE \eqref{2dim_affine}.
First note that \ $(\tY_t, \tX_t)_{t\in\RR_+} \distre (\cY_t, \cX_t)_{t\in\RR_+}$,
 \ and, by Proposition \ref{Pro_affine_2}, it is an affine process having admissible parameters
 \[
   \left( \begin{bmatrix}
           0 & 0 \\
           0 & 0
          \end{bmatrix} , \,
          \frac{1}{2} \begin{bmatrix}
                       1 & 0 \\
                       0 & 1
                      \end{bmatrix} , \,
          \begin{bmatrix}
           a \\
           m
          \end{bmatrix} , \,
          \begin{bmatrix}
           0 & 0 \\
           0 & 0
          \end{bmatrix} , \,
          0 ,
          0 \right) ,
 \]
 and condition \eqref{cond1} is trivially fulfilled.
Hence, by part (ii) of Remark \ref{Rem1}, we have
 \begin{align}\label{help14}
  \left( n^{-1} \tY_{nt}, n^{-1} \tX_{nt} \right)_{t\in\RR_+}
  \distre (\cY_t, \cX_t)_{t\in\RR_+}
  \qquad \text{for all \ $n \in \NN$.}
 \end{align}
Consequently, for all \ $n \in \NN$, \ we have
 \begin{align*}
   \biggl(\frac{1}{n^2} \sum_{i=1}^n \tX_{i-1} ,
          \frac{1}{n^3} \sum_{i=1}^n \tX_{i-1}^2 ,
          \frac{1}{n} \tX_n ,
          \frac{1}{n} \tX_0 ,
          \frac{1}{n^2} \sum_{i=1}^n (\tX_i - \tX_{i-1})^2 \biggr) 
     \distre (A_n, B_n, C_n, D_n, E_n) ,
 \end{align*}
 where
 \begin{align}\nonumber
  A_n &:= \frac{1}{n} \sum_{i=1}^n \cX_{(i-1)/n}
       \as \int_0^1 \cX_t \, \dd t , \qquad \text{as \ $n \to \infty$,} \\ \nonumber
  B_n &:= \frac{1}{n} \sum_{i=1}^n \cX_{(i-1)/n}^2
       \as \int_0^1 \cX_t^2 \, \dd t , \qquad \text{as \ $n \to \infty$,} \\ \nonumber
  C_n &:= \cX_1 , \\ \nonumber
  D_n &:= \cX_0 , \\\label{help20}
  E_n &:= \sum_{i=1}^n (\cX_{i/n} - \cX_{(i-1)/n})^2
       \stoch \int_0^1 \cY_t \, \dd t , \qquad \text{as \ $n \to \infty$.}
 \end{align}
Here the first two convergences are consequences of the definition of the Riemann integral
 using also that \ $(\cX_t)_{t\in\RR_+}$ \ has continuous sample paths almost surely.
The convergence \eqref{help20} can be checked as follows.
With the notations of Jacod and Shiryaev \cite{JSh},
 \ $\bigl(\tau_n:=\bigl(\frac{i}{n}\land 1\bigr)_{i\in\NN}\bigr)_{n\in\NN}$
 \ is a Riemann sequence of (adapted) subdivisions and hence, by Jacod and
 Shiryaev \cite[Theorem I.4.47]{JSh}, the sequence of processes
 \[
   \left(\sum_{i=1}^n
          \left(\cX\left(\frac{i}{n}\land 1\land t\right)
                -\cX\left(\frac{i-1}{n}\land 1\land t\right)\right)^2
   \right)_{t\in\RR_+},
   \qquad n\in\NN,
 \]
 converges to the quadratic variation process of \ $\cX$ \ in probability,
 uniformly on every compact interval.
Especially, with \ $t=1$, \ using also the SDE \eqref{help16}, we have
 \eqref{help20}.
Hence, in order to prove \eqref{help10}, it suffices to show convergences
 \begin{align}
  &\frac{1}{n^2} \sum_{i=0}^{n-1} X_i - \frac{1}{n^2} \sum_{i=0}^{n-1} \tX_i
   \stoch 0 , \label{conv1} \\
  &\frac{1}{n^3} \sum_{i=0}^{n-1} X_i^2 - \frac{1}{n^3} \sum_{i=0}^{n-1} \tX_i^2
   \stoch 0 , \label{conv2} \\
  &\frac{1}{n} X_n - \frac{1}{n} \tX_n \stoch 0 , \label{conv3}\\
  &\frac{1}{n} X_0 - \frac{1}{n} \tX_0 \stoch 0 , \label{conv4} \\
  &\frac{1}{n^2} \sum_{i=1}^n (X_i - X_{i-1})^2
   - \frac{1}{n^2} \sum_{i=1}^n (\tX_i - \tX_{i-1})^2
   \stoch 0 , \label{conv5}
 \end{align}
 as \ $n \to \infty$.
\ Indeed, one can refer to Slutsky's lemma using also that for any random vectors
 \ $U_n$, \ $V_n$, \ $n\in\NN$, \ $U$, \ $V$ \ such that \ $U_n \stoch U$ \ and \ $V_n \stoch V$ \
 as \ $n \to \infty$, we have \ $(U_n,V_n) \stoch (U,V)$ \ as \ $n \to \infty$,
 \ see, e.g., van der Vaart \cite[Theorem 2.7, part (vi)]{Vaa}.

The convergence \eqref{conv4} is trivial.
Next we show
 \begin{equation}\label{EY-tY}
  \EE(|Y_t - \tY_t|) \leq \EE(Y_0) , \qquad t \in \RR_+ .
 \end{equation}
By \eqref{2dim_affine} and \eqref{help16_b}, we have
 \[
   Y_t - \tY_t = Y_0 + \int_0^t (\sqrt{Y_s} - \sqrt{\tY_s}) \, \dd W_s , \qquad
   t \in \RR_+ .
 \]
For each \ $n \in \NN$, \ there exists an even and twice continuously
 differentiable function \ $\psi_n : \RR \to \RR_+$ \ with \ $\vert\psi_n(x)\vert\leq \vert x\vert$,
 \ $|\psi_n'(x)| \leq 1$, \ $\psi_n(x) \uparrow |x|$ \ as \ $n \to \infty$ \
 for all \ $x \in \RR$, \ and
 \begin{align*}
   \psi_n''(x-y)(\sqrt{x} - \sqrt{y})^2
      \leq \frac{2(\sqrt{x} - \sqrt{y})^2}{n\vert x-y\vert}
      \leq \frac{2}{n}
 \end{align*}
 for all \ $n\in\NN$ \ and \ $x,y\in\RR_+$, \ see, e.g., in Karatzas and Shreve \cite[Proof of Proposition 5.2.13]{KarShr}.
By It\^o's formula,
 \begin{align}\label{help51}
  \begin{split}
  \psi_n(Y_t - \tY_t)
   &= \psi_n(Y_0)
     + \frac{1}{2}
       \int_0^t \psi_n''(Y_s - \tY_s) \left(\sqrt{Y_s} - \sqrt{\tY_s}\right)^2 \, \dd s\\
   &\phantom{=\;}  + \int_0^t \psi_n'(Y_s - \tY_s) \left(\sqrt{Y_s} - \sqrt{\tY_s}\right) \, \dd W_s
  \end{split}
 \end{align}
 for all \ $t \in \RR_+$ \ and \ $n \in \NN$.
\ The last term is an \ $(\cF_t)_{t\in\RR_+}$-martingale, since
 \begin{align*}
   \EE\left( \int_0^t
              |\psi_n'(Y_s - \tY_s)| \left(\sqrt{Y_s} - \sqrt{\tY_s}\right)^2
              \, \dd s \right)
   &\leq \EE\left( \int_0^t |Y_s - \tY_s| \, \dd s \right)\\
   &\leq \int_0^t (\EE(Y_s) + \EE(\tY_s)) \, \dd s
    < \infty,
 \end{align*}
 where the last inequality follows by Lemma \ref{Pro_moments}.
Thus the expectation of the last term on the right hand side of \eqref{help51} is zero,
 whereas the expectation of the second integral is bounded by \ $2t/n$.
\ We conclude
 \[
   \EE(\psi_n(Y_t - \tY_t)) \leq \EE(\psi_n(Y_0)) + \frac{t}{n} ,
   \qquad t \in \RR_+ , \quad n \in \NN .
 \]
By monotone convergence theorem, we have
 \begin{align*}
  &\EE(\vert Y_t - \tY_t\vert)
      = \EE(\lim_{n\to\infty} \psi_n(Y_t - \tY_t) )
       = \lim_{n\to\infty} \EE(\psi_n(Y_t - \tY_t))\\
  &\quad \leq \liminf_{n\to\infty} \left(\EE(\psi_n(Y_0)) + \frac{t}{n}\right)
       = \lim_{n\to\infty} \EE(\psi_n(Y_0))
       = \EE(\lim_{n\to\infty} \psi_n(Y_0))
       =\EE(\vert Y_0\vert),
 \end{align*}
 which yields \eqref{EY-tY}.

Next, we derive
 \begin{equation}\label{EX-tX}
  \EE(|X_t - \tX_t|) \leq \EE(|X_0|) + \sqrt{t \EE(Y_0)} , \qquad t \in \RR_+ .
 \end{equation}
Again by \eqref{2dim_affine} and \eqref{help16_b}, we have
 \begin{align}\label{help52}
   X_t - \tX_t = X_0 + \int_0^t \left(\sqrt{Y_s} - \sqrt{\tY_s}\right) \, \dd B_s , \qquad
   t \in \RR_+ ,
 \end{align}
 hence
 \begin{align*}
  \EE(|X_t - \tX_t|)
  &\leq \EE(|X_0|)
        + \sqrt{\EE\left(\left(\int_0^t
                                (\sqrt{Y_s} - \sqrt{\tY_s})
                                \, \dd B_s\right)^2\right)} \\
  &= \EE(|X_0|)
     + \sqrt{\EE\left(\int_0^t \left(\sqrt{Y_s} - \sqrt{\tY_s}\right)^2 \, \dd s\right)}\\
  & \leq \EE(|X_0|) + \sqrt{\EE\left(\int_0^t |Y_s - \tY_s| \, \dd s\right)} \\
  &= \EE(|X_0|) + \sqrt{\int_0^t \EE(|Y_s - \tY_s|) \, \dd s} ,
 \end{align*}
 which yields \eqref{EX-tX} by \eqref{EY-tY}.

By \eqref{EX-tX}, we have
 \[
   \EE\left(\left|\frac{1}{n} X_n - \frac{1}{n} \tX_n\right|\right)
   \leq \frac{1}{n} \left( \EE(|X_0|) + \sqrt{n \EE(Y_0)} \right)
   \to 0 \qquad \text{as \ $n \to \infty$,}
 \]
 hence we obtain \eqref{conv3}.
In a similar way,
 \[
   \EE\left(\left| \frac{1}{n^2} \sum_{i=0}^{n-1} X_i
                   - \frac{1}{n^2} \sum_{i=0}^{n-1} \tX_i \right|\right)
   \leq \frac{1}{n^2}
        \sum_{i=0}^{n-1} \left( \EE(|X_0|) + \sqrt{i \EE(Y_0)} \right)
   \to 0
 \]
 as \ $n \to \infty$, \ hence we obtain \eqref{conv1}.

We also have
 \begin{equation}\label{EX-tX2}
  \EE\left((X_t - \tX_t)^2\right) \leq 2 \EE(X_0^2) + 2 t \EE(Y_0) , \qquad t \in \RR_+ .
 \end{equation}
Indeed, by \eqref{help52}, using Minkowski inequality, we have
 \begin{align*}
  \sqrt{\EE\left((X_t - \tX_t)^2\right)}
  &\leq \sqrt{\EE(X_0^2)}
        + \sqrt{\EE\left(\left(\int_0^t
                                (\sqrt{Y_s} - \sqrt{\tY_s})
                                \, \dd B_s\right)^2\right)} \\
  &= \sqrt{\EE(X_0^2)}
     + \sqrt{\EE\left(\int_0^t \left(\sqrt{Y_s} - \sqrt{\tY_s}\right)^2 \, \dd s\right)}\\
  & \leq \sqrt{\EE(X_0^2)}
        + \sqrt{\EE\left(\int_0^t |Y_s - \tY_s| \, \dd s\right)} \\
  &= \sqrt{\EE(X_0^2)} + \sqrt{\int_0^t \EE(|Y_s - \tY_s|) \, \dd s}
   \leq \sqrt{\EE(X_0^2)} + \sqrt{t \EE(Y_0)}
 \end{align*}
 by \eqref{EY-tY}, which yields \eqref{EX-tX2}.
In a similar way,
 \begin{equation}\label{EX2}
  \EE(X_t^2) \leq 3 \EE(X_0^2) + 3 m^2t^2 + 3 t \EE(Y_0) + 3 a t^2 /2 , \qquad
  t \in \RR_+ .
 \end{equation}
 since, by \eqref{2dim_affine} and \eqref{Pro_moments},
 \begin{align*}
  \sqrt{\EE(X_t^2)}
  &\leq \sqrt{\EE(X_0^2)}
        + |m| t
        + \sqrt{\EE\left(\left(\int_0^t \sqrt{Y_s} \, \dd B_s\right)^2\right)}\\
  & = \sqrt{\EE(X_0^2)}
     + |m| t
     + \sqrt{\EE\left(\int_0^t Y_s \, \dd s\right)} \\
  &= \sqrt{\EE(X_0^2)}
     + |m| t
     + \sqrt{\int_0^t (\EE(Y_0) + a s) \, \dd s}\\
  &= \sqrt{\EE(X_0^2)} + |m| t + \sqrt{t \EE(Y_0) + a t^2 /2}
 \end{align*}
 for \ $t\in\RR_+$, \ which yields \eqref{EX2}.
Clearly, \eqref{EX2} implies also \ $\EE(\tX_t^2) \leq 3m^2t^2 + 3 a t^2 /2$
 \ for all \ $t \in \RR_+$, \ and hence, together with \eqref{EX-tX2} and
 \eqref{EX2}, we conclude
 \begin{align*}
   &\EE\left(\left| \frac{1}{n^3} \sum_{i=0}^{n-1} X_i^2
                   - \frac{1}{n^3} \sum_{i=0}^{n-1} \tX_i^2 \right|\right)
   \leq \frac{1}{n^3}
         \sum_{i=0}^{n-1} \EE\left(|(X_i - \tX_i) (X_i + \tX_i)|\right) \\
   &\quad\leq \frac{1}{n^3}
         \sum_{i=0}^{n-1} \sqrt{ \EE((X_i - \tX_i)^2) \EE((X_i + \tX_i)^2)} \\
   &\quad\leq \frac{1}{n^3}
         \sum_{i=0}^{n-1}
          \sqrt{ 2 \EE((X_i - \tX_i)^2) (\EE(X_i^2) + \EE(\tX_i^2))} \\
   &\quad\leq \frac{1}{n^3}
         \sum_{i=0}^{n-1}
          \sqrt{  12 (\EE(X_0^2) + i \EE(Y_0))
                 \left(\EE(X_0^2) + (\EE(Y_0)) i + (2m^2 + a)i^2 \right)}
    \to 0
 \end{align*}
 as \ $n \to \infty$, \ hence we obtain \eqref{conv2}.

Next, we show \eqref{conv5}.
Again by \eqref{2dim_affine} and \eqref{help16_b}, we have
 \[
   X_i - X_{i-1} = m + \int_{i-1}^i \sqrt{Y_s} \, \dd B_s , \qquad
   \tX_i - \tX_{i-1} = m + \int_{i-1}^i \sqrt{\tY_s} \, \dd B_s , \qquad
   i \in \NN ,
 \]
 and hence
 \begin{align*}
  &\frac{1}{n^2} \sum_{i=1}^n (X_i - X_{i-1})^2
   - \frac{1}{n^2} \sum_{i=1}^n (\tX_i - \tX_{i-1})^2 \\
  &\qquad= \frac{2m}{n^2} \int_0^n \left(\sqrt{Y_s} - \sqrt{\tY_s}\right) \, \dd B_s\\
  &\qquad\phantom{=\;}
     + \frac{1}{n^2}
       \sum_{i=1}^n
        \left[\left(\int_{i-1}^i \sqrt{Y_s} \, \dd B_s\right)^2
              - \left(\int_{i-1}^i \sqrt{\tY_s} \, \dd B_s\right)^2\right]\\
  &\qquad=: 2m R_n + S_n .
 \end{align*}
Here, by \eqref{EY-tY},
 \begin{align*}
  \EE(R_n^2)
  &= \frac{1}{n^4}
     \EE\left(\left(\int_0^n
                     \left(\sqrt{Y_s} - \sqrt{\tY_s}\right) \, \dd B_s\right)^2\right)\\
  & = \frac{1}{n^4}
     \EE\left(\int_0^n \left(\sqrt{Y_s} - \sqrt{\tY_s}\right)^2 \, \dd s\right)
   \leq \frac{1}{n^4} \EE\left(\int_0^n |Y_s - \tY_s| \, \dd s\right)\\
  & = \frac{1}{n^4} \int_0^n \EE(|Y_s - \tY_s|) \, \dd s
   \leq \frac{\EE(Y_0)}{n^3}
   \to 0 \qquad \text{as \ $n \to \infty$,}
 \end{align*}
 hence \ $R_n \stoch 0$ \ as \ $n \to \infty$.
\ Further, by \eqref{EY-tY},
 \begin{align*}
  &\EE(|S_n|)
   = \EE\left(\left| \frac{1}{n^2}
                     \sum_{i=1}^n
                      \int_{i-1}^i
                       \left(\sqrt{Y_s} - \sqrt{\tY_s}\right)
                       \, \dd B_s
                      \int_{i-1}^i
                       \left(\sqrt{Y_s} + \sqrt{\tY_s}\right)
                       \, \dd B_s \right|\right) \\
  &\leq \frac{1}{n^2}
        \sum_{i=1}^n
         \EE\left(\left| \int_{i-1}^i
                          \left(\sqrt{Y_s} - \sqrt{\tY_s}\right)
                          \, \dd B_s
                         \int_{i-1}^i
                          \left(\sqrt{Y_s} + \sqrt{\tY_s}\right)
                          \, \dd B_s \right|\right) \\
  &\leq \frac{1}{n^2}
        \sum_{i=1}^n
         \sqrt{\EE\left(\left( \int_{i-1}^i
                                (\sqrt{Y_s} - \sqrt{\tY_s})
                                \, \dd B_s\right)^2\right)
               \EE\left(\left( \int_{i-1}^i
                                (\sqrt{Y_s} + \sqrt{\tY_s})
                                \, \dd B_s\right)^2\right)} \\
  &= \frac{1}{n^2}
     \sum_{i=1}^n
      \sqrt{\int_{i-1}^i \EE\left((\sqrt{Y_s} - \sqrt{\tY_s})^2\right) \, \dd s
            \int_{i-1}^i \EE\left((\sqrt{Y_s} + \sqrt{\tY_s})^2\right) \, \dd s} \\
  &\leq \frac{1}{n^2}
        \sum_{i=1}^n
         \sqrt{\int_{i-1}^i \EE(|Y_s - \tY_s|) \, \dd s
               \int_{i-1}^i 2 (\EE(Y_s) + \EE(\tY_s)) \, \dd s}\\
  &\leq \frac{1}{n^2}
        \sum_{i=1}^n
         \sqrt{\EE(Y_0) \int_{i-1}^i 2(\EE(Y_0) + 2 a s) \, \dd s} \\
  &= \frac{1}{n^2}
     \sum_{i=1}^n
      \sqrt{2 \EE(Y_0) (\EE(Y_0) + (2 i - 1) a)}
   \to 0 \qquad \text{as \ $n \to \infty$,}
 \end{align*}
 thus \ $S_n \stoch 0$ \ as \ $n \to \infty$, \ and we obtain \eqref{conv5},
 and hence \eqref{help10}.

Finally, by \eqref{help10} and the continuous mapping theorem, and using that
 \begin{align*}
  n\widetilde\theta_n^{\mathrm{LSE}}
  = \frac{\frac{m}{n^2} \sum_{i=1}^n X_{i-1}
          - \frac{1}{2n^2} X_n^2
          + \frac{1}{2n^2} X_0^2
          + \frac{1}{2n^2}
            \sum_{i=1}^n (X_i - X_{i-1})^2}
         {\frac{1}{n^3} \sum_{i=1}^n X_{i-1}^2} ,
 \end{align*}
 we have the assertion.
Indeed, the function \ $g : \RR^5 \to \RR$, \ defined by
 \begin{align*}
  g(x, y, z, u, v)
  :=\begin{cases}
     \frac{m x - (z^2 - u^2 - v)/2}{y} & \text{if \ $y \ne 0$,} \\
     0 & \text{if \ $y = 0$,}
    \end{cases}
 \end{align*}
 is continuous on the set \ $\{(x, y, z, u, v) \in \RR^5 : y \ne 0\}$, \ and
 the limit distribution in \eqref{help10} is concentrated on this set since
 \ $\PP\bigl(\int_0^1 \cX_t^2 \, \dd t > 0\bigr) = 1$.
\ Indeed, if \ $\PP(\int_0^1\cX_t^2\,\dd t =0 )>0$ \ held, then, by the almost sure continuity of the
 sample paths of \ $\cX$, \ we would have \ $\PP(\cX_t = 0,\; \forall\, t\in[0,1])>0$.
\ Hence on the event
 \ $\{\omega\in\Omega : \cX_t(\omega) = 0,  \; \forall\, t\in[0,1]\}$, \ the quadratic variation of
 \ $\cX$ \ would be identically zero.
Since \ $\dd \cX_t = m\,\dd t + \sqrt{\cY_t}\,\dd \cB_t$, $t\in\RR_+$, \ the quadratic variation of
 \ $\cX$ \ is the process \ $\bigl(\int_0^t\cY_u\,\dd u\bigr)_{t\in\RR_+}$, \
 and then we would have
 \[
   \int_0^t \cY_u\,\dd u = 0 \quad \text{for all \ $t\in[0,1]$}
 \]
 on the event \ $\{\omega\in\Omega : \cX_t(\omega) = 0,  \; \forall\, t\in[0,1]\}$.
\ This yields us to a contradiction similarly as at the beginning of the proof
 due to that \ $a\in\RR_{++}$ \ and \ $\dd\cY_t = a\,\dd t + \sqrt{\cY_t}\,\dd\cW_t$, $t\in\RR_+$.
\ Hence the continuous mapping theorem
 (see, e.g., Theorem 2.3 in van der Vaart \cite{Vaa}) yields
 \begin{multline*}
  g\left( \frac{1}{n^2} \sum_{i=1}^n X_{i-1} ,
          \frac{1}{n^3} \sum_{i=1}^n X_{i-1}^2 ,
          \frac{1}{n} X_n ,
          \frac{1}{n} X_0 ,
          \frac{1}{n^2} \sum_{i=1}^n (X_i-X_{i-1})^2 \right) \\
  \distr
  g\left( \int_0^1 \cX_t \, \dd t ,
          \int_0^1 \cX_t^2 \, \dd t ,
          \cX_1 ,
            0 ,
          \int_0^1 \cY_t \, \dd t \right)
  \qquad \text{as \ $n \to \infty$.}
 \end{multline*}
\ Since
 \begin{align*}
  &\PP\left( n\widetilde\theta_n^{\mathrm{LSE}}
     = g\left( \frac{1}{n^2} \sum_{i=1}^n X_{i-1} ,
            \frac{1}{n^3} \sum_{i=1}^n X_{i-1}^2 ,
            \frac{1}{n} X_n ,
            \frac{1}{n} X_0 ,
            \frac{1}{n^2} \sum_{i=1}^n (X_i - X_{i-1})^2 \right) \right)\\
  &\qquad\geq \PP\left(\sum_{i=1}^n X_{i-1}^2 > 0 \right)
    = 1
 \end{align*}
 for all \ $n \geq 2$, \ the assertion follows using \eqref{help49} and that if \ $\xi_n,$
 \ $\eta_n$, \ $n \in \NN$, \ and \ $\xi$ \ are random variables such that
 \ $\xi_n \distr \xi$ \ as \ $n \to \infty$ \ and
 \ $\lim_{n\to\infty} \PP(\xi_n = \eta_n) = 1$, \ then \ $\eta_n \distr \xi$ \ as
 \ $n \to \infty$, \ see, e.g., Barczy et al.\ \cite[Lemma 3.1]{BarIspPap1}.
\proofend

\begin{Rem}\label{Rem3}
If the affine diffusion process given by the SDE \eqref{2dim_affine} is
 critical but \ $(b,\theta)\ne(0,0)$ \ (i.e., \ $b=0$, $\theta>0$ \ or \ $b>0$,
 $\theta=0$), \ then the asymptotic behaviour of the LSE
 \ $\widetilde\theta_n^{\mathrm{LSE}}$ \ cannot be studied using Theorem \ref{Thm1}
 since its condition \ $\lim_{\theta\to\infty} \theta \beta^{(\theta)} = \beta$ \ is
 not satisfied.
\proofend
\end{Rem}

\subsection{Least squares estimator of \ $(\theta,m)$}
\label{section_LSE_M_theta}

The LSE of \ $(\theta,m)$ \ based on the observations \ $X_i$,
 \ $i=0,1,\ldots,n$, \ can be obtained by solving the extremum problem
 \[
   (\widehat\theta_n^{\mathrm{LSE}}, \widehat m_n^{\mathrm{LSE}})
   := \argmin_{(\theta,m)\in\RR^2}
       \sum_{i=1}^n (X_i - X_{i-1} - (m-\theta X_{i-1}))^2.
 \]
We need to solve the following system of equations with respect to
 \ $(\theta,m)$:
 \begin{align*}
  &2\sum_{i=1}^n (X_i - X_{i-1} - (m-\theta X_{i-1})) X_{i-1} = 0,\\
  &2\sum_{i=1}^n (X_i - X_{i-1} - (m-\theta X_{i-1})) = 0,
 \end{align*}
 which can be written also in the form
 \[
   \begin{bmatrix}
    \sum_{i=1}^n X_{i-1}^2 & - \sum_{i=1}^n X_{i-1} \\
    - \sum_{i=1}^n X_{i-1} & n
   \end{bmatrix}
   \begin{bmatrix}
    \theta \\
    m
   \end{bmatrix}
   =
   \begin{bmatrix}
    - \sum_{i=1}^n (X_i - X_{i-1}) X_{i-1} \\
    \sum_{i=1}^n (X_i - X_{i-1})
   \end{bmatrix}.
 \]
Then one can check that
 \begin{align}\label{help12}
  \widehat\theta_n^{\mathrm{LSE}}
  = - \frac{n\sum_{i=1}^n (X_i - X_{i-1})X_{i-1}
            -  \sum_{i=1}^n X_{i-1} \sum_{i=1}^n (X_i-X_{i-1})}
           {n\sum_{i=1}^n X_{i-1}^2 - \left(\sum_{i=1}^n X_{i-1}\right)^2},
 \end{align}
 and
 \begin{align}\label{help13}
  \widehat m_n^{\mathrm{LSE}}
  = \frac{\sum_{i=1}^n X_{i-1}^2 \sum_{i=1}^n (X_i-X_{i-1})
          -  \sum_{i=1}^n X_{i-1} \sum_{i=1}^n (X_i-X_{i-1})X_{i-1}}
         {n\sum_{i=1}^n X_{i-1}^2 - \left(\sum_{i=1}^n X_{i-1}\right)^2}
 \end{align}
 provided that \ $n\sum_{i=1}^n X_{i-1}^2 - \left(\sum_{i=1}^n X_{i-1}\right)^2>0$.
\ Since the matrix
 \[
   \begin{bmatrix}
    2\sum_{i=1}^n X_{i-1}^2 & -2\sum_{i=1}^n X_{i-1} \\
    -2\sum_{i=1}^n X_{i-1} & 2n
   \end{bmatrix}
 \]
 which consists of the second order partial derivatives of the function
 $\RR^2\ni(\theta,m)\mapsto\sum_{i=1}^n (X_i - X_{i-1} - (m-\theta X_{i-1}))^2$
 \ is positive definite provided that
 \ $n\sum_{i=1}^n X_{i-1}^2 - \left(\sum_{i=1}^n X_{i-1}\right)^2>0$, \ we have
 \ $(\widehat\theta_n^{\mathrm{LSE}},\widehat m_n^{\mathrm{LSE}})$ \ is indeed the
 solution of the extremum problem provided that
 \ $n\sum_{i=1}^n X_{i-1}^2 - \left(\sum_{i=1}^n X_{i-1}\right)^2>0$.

\begin{Thm}\label{Thm3}
Let us assume that Condition (C) holds.
Then
 \begin{align}\label{help26}
  \PP\left(n\sum_{i=1}^n X_{i-1}^2 - \left(\sum_{i=1}^n X_{i-1}\right)^2>0\right)=1
  \qquad \text{for all \ $n\geq 2$,}
 \end{align}
 and there exists a unique LSE
 \ $(\widehat\theta_n^{\mathrm{LSE}},\widehat m_n^{\mathrm{LSE}})$ \ which has the form
 given in \eqref{help12} and \eqref{help13}.
Further,
 \begin{align*}
  n\widehat\theta_n^{\mathrm{LSE}}
  \distr
  - \frac{\int_0^1 \cX_t\,\dd \cX_t - \cX_1\int_0^1 \cX_t\,\dd t}
         {\int_0^1 \cX_t^2\,\dd t - \left(\int_0^1 \cX_t\,\dd t\right)^2}
  \qquad \text{as \ $n\to\infty$,}
 \end{align*}
 and
 \begin{align*}
  \widehat m_n^{\mathrm{LSE}}
  \distr
  \frac{\cX_1\int_0^1 \cX_t^2\,\dd t
        - \int_0^1 \cX_t\,\dd t \int_0^1 \cX_t\,\dd \cX_t}
       {\int_0^1 \cX_t^2\,\dd t - \left(\int_0^1 \cX_t\,\dd t\right)^2}
    \qquad \text{as \ $n\to\infty$,}
 \end{align*}
 where \ $(\cX_t)_{t\in\RR_+}$ \ is the second coordinate of a two-dimensional
 affine process\ $(\cY_t,\cX_t)_{t\in\RR_+}$ \ given by the unique strong solution of the SDE
 \begin{align*}
 \begin{cases}
  \dd\cY_t = a\,\dd t + \sqrt{\cY_t}\,\dd \cW_t,\\[2mm]
  \dd\cX_t = m\,\dd t + \sqrt{\cY_t}\,\dd \cB_t,
 \end{cases}
 \qquad t\in\RR_+,
 \end{align*}
 with initial value \ $(\cY_0,\cX_0)=(0,0)$, \ where \ $(\cW_t)_{t\in\RR_+}$ \ and
 \ $(\cB_t)_{t\in\RR_+}$ \ are independent standard Wiener processes.
\end{Thm}

\begin{Rem}\label{Rem6}
(i) The limit distributions in Theorem \ref{Thm3} have the same forms as those
 of the limit distributions in Theorem 4.1 in Hu and Long \cite{HuLon3}.

\noindent
(ii) By Proposition \ref{Pro_affine_2}, the affine process \ $(\cY_t,\cX_t)_{t\in\RR_+}$ \ given in
 Theorem \ref{Thm3} has infinitesimal generator
 \begin{align*}
  (\cA_{(\cY,\cX)} f)(x)
  = \frac{1}{2}x_1f_{1,1}''(x) + \frac{1}{2}x_1f_{2,2}''(x) + af_1'(x) + mf_2'(x),
 \end{align*}
 where \ $x=(x_1,x_2)\in\RR_+\times\RR$ \ and \ $f\in\cC^2_c(\RR_+\times\RR)$.

\noindent(iii)
Under the Condition (C), by Theorem \ref{Thm3} and Slutsky's lemma, we get
 \ $\widehat\theta_n^{\mathrm{LSE}}$ \ converges stochastically to the parameter \ $\theta=0$ \
 as \ $n\to\infty$, \ and one can show that \ $\widehat m_n^{\mathrm{LSE}}$ \ does not converge stochastically
 to the parameter \ $m$ \ as \ $n\to\infty$, \ see Appendix \ref{AppendixB}.
\proofend
\end{Rem}

\noindent{\bf Proof of Theorem \ref{Thm3}.}
By an easy calculation,
\begin{align*}
  n\sum_{i=1}^n X_{i-1}^2 - \left(\sum_{i=1}^n X_{i-1}\right)^2
    = n \sum_{i=1}^n \left( X_{i-1} - \frac{1}{n}\sum_{j=1}^n X_{j-1} \right)^2
    \geq 0,
 \end{align*}
 and equality holds if and only if
 \[
   X_{i-1} = \frac{1}{n} \sum_{j=1}^n X_{j-1},\quad i=1,\ldots,n
   \qquad \Longleftrightarrow \qquad
    X_0=X_1=\cdots=X_{n-1}.
 \]
By \eqref{help22}, for all \ $n\geq 2$,
 \begin{align*}
 \PP(X_0=X_1=\cdots=X_{n-1})
   &\leq \PP(X_0 = X_1) = \PP\left(\int_0^1\sqrt{Y_s}\,\dd B_s = m\right)\\
   &= \EE\left(\PP\left( \int_0^1\sqrt{Y_s}\,\dd B_s = m \;\Big\vert\; (Y_s)_{s\in[0,1]}\right)\right)=0,
 \end{align*}
 where we used that the conditional distribution of \ $\int_0^1\sqrt{Y_s}\,\dd B_s$ \ given
 \ $(Y_s)_{s\in[0,1]}$ \ is a normal distribution with mean \ $0$ \ and with variance \ $\int_0^1 Y_s\,\dd s$.
\ Here the variance \ $\int_0^1 Y_s\,\dd s$ \ is positive almost surely, see the proof of Theorem \ref{Thm2}.
This yields \eqref{help26}.

By \eqref{help12} and \eqref{help13}, we have
 \begin{align*}
  &n\widehat\theta_n^{\mathrm{LSE}}
   = - \frac{\frac{1}{n^2}\sum_{i=1}^n (X_i - X_{i-1})X_{i-1}
             - \frac{1}{n^2} \sum_{i=1}^n X_{i-1} \frac{1}{n}(X_n-X_0)}
            {\frac{1}{n^3}\sum_{i=1}^n X_{i-1}^2
             - \left(\frac{1}{n^2}\sum_{i=1}^n X_{i-1}\right)^2},\\
  &\widehat m_n^{\mathrm{LSE}}
   = \frac{\frac{1}{n^3}\sum_{i=1}^n X_{i-1}^2 \frac{1}{n}(X_n-X_0)
           - \frac{1}{n^2}\sum_{i=1}^n X_{i-1}
             \frac{1}{n^2}\sum_{i=1}^n (X_i - X_{i-1})}
          {\frac{1}{n^3}\sum_{i=1}^n X_{i-1}^2
           - \left(\frac{1}{n^2}\sum_{i=1}^n X_{i-1}\right)^2},
 \end{align*}
 and using \eqref{help53} and \eqref{help10}, as in the proof of Theorem \ref{Thm2}, the continuous
 mapping theorem yields the assertion.
We only remark that
 \begin{align}\label{help37}
  \PP\left(\int_0^1 \cX_t^2\,\dd t - \left(\int_0^1 \cX_t\,\dd t\right)^2
         > 0\right) = 1.
 \end{align}
Indeed,
 \begin{align*}
   \int_0^1 \cX_t^2\,\dd t - \left(\int_0^1 \cX_t\,\dd t\right)^2
     = \int_0^1 \left(\cX_t - \int_0^1 \cX_s\,\dd s\right)^2
     \geq 0,
 \end{align*}
 and equality holds if and only if
 \[
   \cX_t = \int_0^1 \cX_s\,\dd s
    \quad \text{a.e. \ $t\in[0,1]$.}
 \]
Since \ $\cX$ \ has continuous sample paths almost surely,
 \begin{align}\label{help41}
  \PP\left( \int_0^1 \cX_t^2\,\dd t - \left(\int_0^1 \cX_t\,\dd t\right)^2 = 0 \right)>0
 \end{align}
 holds if and only if
 \[
   \PP\left( \cX_t = \int_0^1 \cX_s\,\dd s,\; \forall\,t\in[0,1] \right)>0.
 \]
Hence, since \ $\cX_0=0$, \ we have \eqref{help41} holds if and only if
 \ $\PP(\cX_t = 0,\;\forall\,t\in[0,1] )>0$, \ which is a contradiction due to our assumption
 \ $a\in\RR_{++}$ \ (for more details, see the proof of Theorem \ref{Thm2}).
\proofend

\subsection{Conditional least squares estimator of \ $(\theta,m)$}

For all \ $t\in\RR_+$, \ let \ $\cF^{(Y,X)}_t$ \ be the \ $\sigma$-algebra
 generated by \ $(Y_s,X_s)_{s\in[0,t]}$.
\ The conditional least squares estimator (CLSE) of \ $(\theta,m)$ \ based on
 the observations \ $X_i$, $i=0,1,\ldots,n$, \ can be obtained by solving the
 extremum problem
 \[
   (\widehat\theta_n^{\mathrm {CLSE}},\widehat m_n^{\mathrm{CLSE}})
   := \argmin_{(\theta,m)\in\RR^2} \sum_{i=1}^n
       \Big(X_i - \EE(X_i\mid \cF_{i-1}^{(Y,X)})\Big)^2.
 \]
By \eqref{help22}, for all \ $(y_0,x_0)\in\RR_+\times\RR$, \ we have
  \begin{align*}
    \EE\big( X_t \mid (Y_0,X_0) = (y_0,x_0)\big)
        = \ee^{-\theta t}x_0 + m\int_0^t\ee^{-\theta(t-u)}\,\dd u,
    \qquad t\geq 0,
 \end{align*}
 where we used that \ $\bigl( \int_0^t\ee^{\theta u}\sqrt{Y_u}\,\dd B_u \bigr)_{t\in\RR_+}$ \
 is a martingale (which follows by the proof of Proposition \ref{Pro_moments}).
Using that \ $(Y_t,X_t)_{t\in\RR_+}$ \ is a time-homogeneous Markov process, we have
 \begin{align*}
    \EE( X_t \mid \cF^{(Y,X)}_s)
     = \EE( X_t \mid (Y_s,X_s))
    = \ee^{-\theta(t-s)}X_s + m\int_s^t\ee^{-\theta(t-u)}\,\dd u
 \end{align*}
 for \ $0\leq s\leq t$.
\ Then
 \begin{align*}
    X_i - \EE( X_i \mid \cF^{(Y,X)}_{i-1})
      & = X_i - \ee^{-\theta} X_{i-1} - m \int_{i-1}^i\ee^{-\theta(i-u)}\,\dd u\\
      &  = X_i - \ee^{-\theta} X_{i-1} - m \int_0^1\ee^{-\theta v}\,\dd v\\
      & = X_i - \gamma X_{i-1} - \delta,
       \qquad i=1,\ldots,n,
 \end{align*}
 where
 \[
   \gamma:= \ee^{-\theta} \qquad \text{and} \qquad
   \delta:=m\int_0^1\ee^{-\theta v}\,\dd v
          =\begin{cases}
             m\frac{1-\ee^{-\theta}}{\theta} & \text{if \ $\theta\ne0$,}\\
             m                               & \text{if \ $\theta=0$.}
           \end{cases}
 \]
Hence for all \ $n\in\NN$,
 \begin{align}\label{help29}
  \begin{split}
   & \widehat\gamma_n^{\mathrm{CLSE}} =  \ee^{-\widehat\theta_n^{\mathrm{CLSE}}},\\
   & \widehat\delta_n^{\mathrm{CLSE}}
     = \widehat m_n^{\mathrm{CLSE}}
       \int_0^1\ee^{-\widehat \theta_n^{\mathrm{CLSE}} v}\,\dd v,
  \end{split}
 \end{align}
 where \ $(\widehat\gamma_n^{\mathrm{CLSE}}, \widehat\delta_n^{\mathrm{CLSE}})$ \ is a
 CLSE of \ $(\gamma,\delta)$ \ based on the observations
 \ $X_i$, $i=0,1,\ldots,n$, \ which can be obtained by solving the extremum
 problem
 \begin{align}\label{help40}
   (\widehat\gamma_n^{\mathrm{CLSE}},\widehat\delta_n^{\mathrm{CLSE}})
   := \argmin_{(\gamma,\delta)\in\RR^2} \sum_{i=1}^n
       (X_i - \gamma X_{i-1} - \delta)^2.
 \end{align}
Indeed, the function \ $A:\RR^2\to\RR^2$,
 \[
   \RR^2\ni (\theta',m') \mapsto A(\theta',m')
            := \begin{bmatrix}
               \gamma'\\
               \delta' \\
              \end{bmatrix}
            =: \begin{bmatrix}
               \ee^{-\theta'} \\
               m'\int_0^1\ee^{-\theta'v}\,\dd v \\
              \end{bmatrix}
              \in \RR_+\times\RR
 \]
 is bijective and measurable, and then there is a bijection between the set of CLSEs
 of the parameters \ $(\theta,m)$ \ and the set\vadjust{\vfill\eject} of CLSEs of the parameters \ $A(\theta,m)$.
\ This follows easily, since for
all \ $n\in\NN$, \ $(x_0,x_1,\ldots,x_n)\in\RR^{n+1}$ \ and \ $(\gamma',
\delta')\in\RR_+\times\RR$,
 \begin{align*}
   \sum_{i=1}^n (x_i - \gamma' x_{i-1} - \delta')^2
       = \sum_{i=1}^n \left(x_i -
               \begin{bmatrix}
                 \gamma' \\
                 \delta' \\
                \end{bmatrix}^\top
                \begin{bmatrix}
                 x_{i-1} \\
                 1 \\
                \end{bmatrix}
                \right)^2 
      = \sum_{i=1}^n \left(x_i -
                \left(A(\theta',m')\right)^\top
                \begin{bmatrix}
                 x_{i-1} \\
                 1 \\
                \end{bmatrix}
                \right)^2,
 \end{align*}
 hence \ $(\widehat\theta_n^{\mathrm{CLSE}},\widehat m_n^{\mathrm{CLSE}})$ \ is a CLSE of \ $(\theta,m)$ \ if and only if
 \ $A(\widehat\theta_n^{\mathrm{CLSE}},\widehat m_n^{\mathrm{CLSE}})$ \ is a CLSE of \ $A(\theta,m)$.

For the extremum problem \eqref{help40}, we need to solve the following system of equations with respect to
 \ $(\gamma,\delta)$:
 \begin{align*}
  &2\sum_{i=1}^n ( X_i - \gamma X_{i-1} - \delta ) X_{i-1} = 0,\\
  &2\sum_{i=1}^n ( X_i - \gamma X_{i-1} - \delta ) = 0,
 \end{align*}
 which can be written also in the form
 \[
   \begin{bmatrix}
     \sum_{i=1}^n X_{i-1}^2 & \sum_{i=1}^n X_{i-1} \\
     \sum_{i=1}^n X_{i-1} & n \\
   \end{bmatrix}
   \begin{bmatrix}
     \gamma \\
     \delta \\
   \end{bmatrix}
   =
   \begin{bmatrix}
     \sum_{i=1}^n  X_{i-1}X_i \\
     \sum_{i=1}^n  X_i \\
   \end{bmatrix}.
 \]
Then
 \begin{align}\label{help43}
  \widehat\gamma_n^{\mathrm{CLSE}}
  = \frac{n\sum_{i=1}^n X_{i-1}X_i  -  \sum_{i=1}^n X_{i-1} \sum_{i=1}^n X_i}
         {n\sum_{i=1}^n X_{i-1}^2 - \left(\sum_{i=1}^n X_{i-1}\right)^2},
 \end{align}
 and
 \begin{align}\label{help44}
  \widehat \delta_n^{\mathrm{CLSE}}
  = \frac{\sum_{i=1}^n X_{i-1}^2 \sum_{i=1}^n X_i
          - \sum_{i=1}^n X_{i-1} \sum_{i=1}^n X_{i-1}X_i}
         {n\sum_{i=1}^n X_{i-1}^2 - \left(\sum_{i=1}^n X_{i-1}\right)^2},
 \end{align}
 provided that
 \ $n\sum_{i=1}^n X_{i-1}^2 - \left(\sum_{i=1}^n X_{i-1}\right)^2\ne 0$.
\ Since the matrix
 \[
   \begin{bmatrix}
    2\sum_{i=1}^n X_{i-1}^2 & 2\sum_{i=1}^n X_{i-1} \\
    2\sum_{i=1}^n X_{i-1} & 2n
   \end{bmatrix}
 \]
 consisting of the second order partial derivatives of the function
 \ $\RR^2\ni(\gamma,\delta)\mapsto \sum_{i=1}^n (X_i - \gamma X_{i-1} - \delta)^2$
 \ is positive definite provided that
 \ $n\sum_{i=1}^n X_{i-1}^2 - \left(\sum_{i=1}^n X_{i-1}\right)^2>0$, \ we have
 \ $(\widehat\gamma_n^{\mathrm{CLSE}},\widehat \delta_n^{\mathrm{CLSE}})$ \ is indeed
 the solution of the extremum problem provided that
 \ $n\sum_{i=1}^n X_{i-1}^2 - \left(\sum_{i=1}^n X_{i-1}\right)^2>0$.

\begin{Rem}\label{Rem5}
Using the definition of CLSE of \ $(\theta,m)$ \ we give a mathematical motivation of the definition of the LSE
 \ $\widetilde\theta_n$ \ of \ $\theta$ \ introduced in Section \ref{section_LSE_theta}.
Note that if \ $\theta=0$, \ then
 \[
  X_i - \EE( X_i \mid \cF^{(Y,X)}_{i-1}) = X_i - X_{i-1} - m, \quad i=1,\ldots,n.
\]
If \ $\theta\ne 0$, \ then, by Taylor's theorem, \ $1-\ee^{-\theta} = \ee^{-\tau\theta}\theta$ \ with some
 \ $\tau=\tau(\theta)\in[0,1]$, \ and hence
 \begin{align*}
  X_i - \EE( X_i \mid \cF^{(Y,X)}_{i-1})
   &= X_i - \ee^{-\theta} X_{i-1} - m \int_0^1\ee^{-\theta v}\,\dd v\\
   &= X_i - X_{i-1} + \ee^{-\tau\theta}\theta X_{i-1} - m\ee^{-\tau\theta}
 \end{align*}
 for \ $i=1,\ldots,n-1$.
\ Hence for small values of \ $\theta$ \ one can approximate \ $X_i - \EE( X_i \mid \cF^{(Y,X)}_{i-1})$ \ by
 \ $X_i - X_{i-1} + \theta X_{i-1} - m = X_i - X_{i-1} - (m-\theta X_{i-1})$, $i=1,\ldots,n$.
\ Based on this, for small values of \ $\theta$, \ in the definition of the LSE of \ $\theta$, \ the sum
 \ $\sum_{i=1}^n (X_i - X_{i-1} - (m-\theta X_{i-1}))^2$ \ can be considered as an approximation of the sum
 \ $\sum_{i=1}^n (  X_i - \EE( X_i \mid \cF^{(Y,X)}_{i-1}) )^2$ \ in the definition of the CLSE of
 \ $(\theta,m)$.
\proofend
\end{Rem}

\begin{Thm}\label{Thm4}
Let us assume that Condition (C) holds.
Then
 \begin{align}\label{help28}
  \PP\left(n\sum_{i=1}^n X_{i-1}^2 - \left(\sum_{i=1}^n X_{i-1}\right)^2>0\right)=1
  \qquad \text{for all \ $n\geq 2$,}
 \end{align}
 and there exists a unique CLSE
 \ $(\widehat\theta_n^{\mathrm{CLSE}},\widehat m_n^{\mathrm{CLSE}})$ \ which has the
 form given in \eqref{help29}.
Further,
 \begin{align}\label{help30}
  n\widehat\theta_n^{\mathrm{CLSE}}
  \distr
  - \frac{\int_0^1 \cX_t\,\dd \cX_t - \cX_1\int_0^1 \cX_t\,\dd t}
         {\int_0^1 \cX_t^2\,\dd t - \left(\int_0^1 \cX_t\,\dd t\right)^2}
  \qquad \text{as \ $n\to\infty$,}
 \end{align}
 and
 \begin{align}\label{help31}
  \widehat m_n^{\mathrm{CLSE}}
  \distr
  \frac{\cX_1\int_0^1 \cX_t^2\,\dd t
        - \int_0^1 \cX_t\,\dd t \int_0^1 \cX_t\,\dd \cX_t}
       {\int_0^1 \cX_t^2\,\dd t - \left(\int_0^1 \cX_t\,\dd t\right)^2}
  \qquad \text{as \ $n\to\infty$,}
 \end{align}
 where \ $(\cX_t)_{t\in\RR_+}$ \ is the second coordinate of a two-dimensional
 affine process\ $(\cY_t,\cX_t)_{t\in\RR_+}$ \ given by the unique strong solution of the SDE
 \begin{align*}
  \begin{cases}
   \dd\cY_t = a\,\dd t + \sqrt{\cY_t}\,\dd \cW_t,\\[2mm]
   \dd\cX_t = m\,\dd t + \sqrt{\cY_t}\,\dd \cB_t,
  \end{cases}
 \qquad t\in\RR_+,
 \end{align*}
 with initial value \ $(\cY_0,\cX_0)=(0,0)$, \ where \ $(\cW_t)_{t\in\RR_+}$ \
 and \ $(\cB_t)_{t\in\RR_+}$ \ are independent standard Wiener processes.
\end{Thm}

\noindent{\bf Proof.}
By the proof of Theorem \ref{Thm3}, we have \eqref{help28}.
By \eqref{help43} and \eqref{help44}, for all \ $n\geq 2$ \ we have
 \begin{align*}
  \widehat\gamma_n^{\mathrm{CLSE}} - 1
  = \frac{n\sum_{i=1}^n (X_i - X_{i-1})X_{i-1} - X_n\sum_{i=1}^n X_{i-1}}
         {n\sum_{i=1}^n X_{i-1}^2  - \left( \sum_{i=1}^n X_{i-1} \right)^2},
 \end{align*}
 and
 \begin{align*}
  \widehat\delta_n^{\mathrm{CLSE}}
  = \frac{X_n\sum_{i=1}^n X_{i-1}^2
          - \sum_{i=1}^n X_{i-1}\sum_{i=1}^n (X_i - X_{i-1})X_{i-1}}
         {n\sum_{i=1}^n X_{i-1}^2  - \left( \sum_{i=1}^n X_{i-1} \right)^2}.
 \end{align*}
Using \eqref{help53} and \eqref{help10}, the continuous mapping theorem, by the same technique as
 in the proof of Theorem \ref{Thm3}, we get
 \begin{align}\label{help32}
  n(\widehat\gamma_n^{\mathrm{CLSE}} - 1)
  \distr
  \frac{\int_0^1 \cX_t\,\dd \cX_t - \cX_1\int_0^1 \cX_t\,\dd t}
       {\int_0^1 \cX_t^2\,\dd t - \left(\int_0^1 \cX_t\,\dd t\right)^2}
  \qquad \text{as \ $n\to\infty$,}
 \end{align}
 and
 \begin{align}\label{help33}
  \widehat \delta_n^{\mathrm{CLSE}}
  \distr
  \frac{\cX_1\int_0^1 \cX_t^2\,\dd t
        - \int_0^1 \cX_t\,\dd t \int_0^1 \cX_t\,\dd \cX_t}
       {\int_0^1 \cX_t^2\,\dd t - \left(\int_0^1 \cX_t\,\dd t\right)^2}
  \qquad \text{as \ $n\to\infty$.}
 \end{align}
By Slutsky's lemma, we also have
 \ $\widehat\gamma_n^{\mathrm{CLSE}} \stoch 1$ \ as \ $n\to\infty$.
\ Hence, by Taylor's theorem using also that \ $\widehat\gamma_n^{\mathrm{CLSE}}>0$, $n\in\NN$ \
 (due to its definition given in \eqref{help29}), we have
 \begin{align}\label{help34}
   \widehat\theta_n^{\mathrm{CLSE}}
      = - \log(\widehat\gamma_n^{\mathrm{CLSE}})
      = - \log(\widehat\gamma_n^{\mathrm{CLSE}})  - \log(1)
      = -\frac{1}{\xi_n}(\widehat\gamma_n^{\mathrm{CLSE}} - 1),
 \end{align}
 where \ $\xi_n$ \ is in the interval with endpoints \ $1$ \ and \ $\widehat\gamma_n^{\mathrm{CLSE}}$.
\ Since \ $\widehat\gamma_n^{\mathrm{CLSE}}\stoch 1$ \ as \ $n\to\infty$, \ we have \ $\xi_n\stoch 1$ \ as \ $n\to\infty$,
 \ and hence using the decomposition
 \begin{align*}
   n \widehat\theta_n^{\mathrm{CLSE}}
     = - \frac{1}{\xi_n}n(\widehat\gamma_n^{\mathrm{CLSE}} - 1),
        \qquad n\in\NN,
 \end{align*}
 Slutsky's lemma and \eqref{help32}, we get \eqref{help30}.

Next we turn to prove \eqref{help31}.
For this, by \eqref{help29}, \eqref{help33} and by Slutsky's lemma, it is enough to check that
 \[
  \int_0^1 \ee^{-\widehat\theta_n^{\mathrm{CLSE}}v}\,\dd v \stoch 1
   \qquad \text{as \ $n\to\infty$.}
 \]
Since
 \begin{align*}
   \int_0^1 \ee^{-\widehat\theta_n^{\mathrm{CLSE}}v}\,\dd v
    = \begin{cases}
       \frac{1-\ee^{-\widehat\theta_n^{\mathrm{CLSE}}}}
            {\widehat\theta_n^{\mathrm{CLSE}}}
            & \text{if \ $\widehat \theta_n^{\mathrm{CLSE}}\ne 0 $,}\\
       1   & \text{if \ $\widehat \theta_n^{\mathrm {CLSE}}=0$,}
      \end{cases}
 \end{align*}
 by \eqref{help34}, for all \ $\varepsilon>0$ \ we have
 \begin{align*}
  \PP\left( \left\vert \int_0^1 \ee^{-\widehat\theta_n^{\mathrm {CLSE}}v}\,\dd v -1 \right\vert \geq \varepsilon\right)
  &= \PP\left( \left\vert  \frac{1-\widehat\gamma_n^{\mathrm {CLSE}}}{\widehat\theta_n^{\mathrm {CLSE}}} -1 \right\vert \geq \varepsilon
              \,\Big\vert\, \widehat\theta_n^{\mathrm {CLSE}}\ne 0 \right)
        \PP(\widehat\theta_n^{\mathrm {CLSE}}\ne 0)\\
  &\phantom{=\;} + \PP\left(  \vert  1-1 \vert \geq \varepsilon \,\big \vert\, \widehat\theta_n^{\mathrm {CLSE}}= 0 \right)
         \PP(\widehat\theta_n^{\mathrm {CLSE}}= 0)\\
  & = \PP\left(  \vert  \xi_n-1 \vert \geq \varepsilon \,\big \vert\, \widehat\theta_n^{\mathrm {CLSE}}\ne 0 \right)
         \PP(\widehat\theta_n^{\mathrm {CLSE}}\ne 0) \stoch 0,
 \end{align*}
 since \ $\xi_n\stoch 0$ \ as \ $n\to\infty$.
\proofend

\begin{Rem}
(i) We do not consider the CLSE of \ $\theta$ \ supposing that \ $m$
 \ is known since the corresponding extremum problem is rather complicated, and
 from statistical point of view it has less importance.

\noindent(ii)
Under the Condition (C), by Theorem \ref{Thm4} and Slutsky's lemma, we get
 \ $\widehat\theta_n^{\mathrm{CLSE}}$ \ converges stochastically to the parameter \ $\theta=0$ \
 as \ $n\to\infty$, \ and one can show that \ $\widehat m_n^{\mathrm{CLSE}}$ \ does not converge stochastically
 to the parameter \ $m$ \ as \ $n\to\infty$, \ see Appendix \ref{AppendixB}.
\proofend
\end{Rem}

\medskip

\appendix

\noindent{\bf\Large Appendix}

\section{The integrals in \eqref{affine_inf_gen}}\label{AppendixA}

We check that the integrals in \eqref{affine_inf_gen} are well-defined, i.e., elements of \ $\CC$,
 \ under the conditions (v) and (vi) of Definition \ref{Def_admissible}.
For this, by decomposing a complex-valued function to real and imaginary parts, it is enough to verify that
 \begin{align}\label{help46}
   & \int_{\RR_+\times \RR^d}
      \big(f(x+\xi) - f(x) - \langle f'(x), \xi \rangle\big) x_1 \, \mu(\dd\xi) < \infty,\\ \label{help47}
   & \int_{\RR_+\times \RR^d}
         (f(x+\xi) - f(x) - \langle f'_{(2)}(x), \xi_2 \rangle)
         \, m(\dd\xi) < \infty
 \end{align}
 for all \ $x=(x_1,x_2)\in\RR_+\times\RR^d$ \ and real-valued \ $f\in\cC_c^2(\RR_+\times\RR^d)$.

First we check \eqref{help46}.
If \ $x=(x_1,x_2)\in\RR_+\times\RR^d$ \ and \ $\xi\in\RR_+\times\RR^d$ \ with \ $\Vert\xi\Vert\leq 1$, \ then, by
 \eqref{2help15} with \ $\theta=1$, \
 \begin{align*}
  \vert f(x+\xi) - f(x) - \langle f'(x),\xi\rangle\vert
   = \frac{1}{2}\vert \langle f''(x+\tau\xi)\xi,\xi \rangle \vert
   \leq \frac{1}{2}\Vert f''\Vert_\infty \Vert\xi\Vert^2,
 \end{align*}
 where \ $\tau=\tau(x,\xi)\in[0,1]$.
\ If \ $x=(x_1,x_2)\in\RR_+\times\RR^d$ \ and \ $\xi\in\RR_+\times\RR^d$ \ with \ $\Vert\xi\Vert>1$, \ then, by
 \eqref{help15} with \ $\theta=1$, \
 \begin{align*}
  \vert f(x+\xi) - f(x) - \langle f'(x),\xi\rangle\vert
    & \leq \vert f(x+\xi) - f(x) \vert + \vert \langle f'(x),\xi\rangle \vert
      = \vert \langle f'(x+\tau\xi),\xi\rangle \vert + \vert \langle f'(x),\xi\rangle \vert\\
    & \leq 2 \Vert f'\Vert_\infty \Vert\xi\Vert,
 \end{align*}
 where \ $\tau=\tau(x,\xi)\in[0,1]$.
Hence
 \begin{align*}
  &\left\vert
    \int_{\RR_+\times \RR^d}
      \big(f(x+\xi) - f(x) - \langle f'(x), \xi \rangle\big) x_1 \, \mu(\dd\xi)
   \right\vert\\
  &\quad \leq
   \int_{\{ \xi\in\RR_+\times \RR^d \,:\, \Vert\xi\Vert\leq 1\} }
     \vert(f(x+\xi) - f(x) - \langle f'(x), \xi \rangle\vert x_1 \, \mu(\dd\xi) \\
  &\phantom{\quad\leq\;}
    + \int_{\{ \xi\in\RR_+\times \RR^d \,:\, \Vert\xi\Vert>1 \}}
      \vert(f(x+\xi) - f(x) - \langle f'(x), \xi \rangle\vert x_1 \, \mu(\dd\xi)\\
  &\quad\leq \frac{1}{2}\Vert f''\Vert_\infty x_1
        \int_{\{\xi\in\RR_+\times \RR^d \,:\, \Vert\xi\Vert\leq 1\}}
           \Vert\xi\Vert^2 \, \mu(\dd\xi)\\
  &\phantom{\quad\leq\;}
       + 2 \Vert f'\Vert_\infty x_1
         \int_{\{\xi\in\RR_+\times \RR^d \,:\, \Vert\xi\Vert>1\}}
           \Vert\xi\Vert \, \mu(\dd\xi)
    <\infty,
 \end{align*}
 where the last inequality follows by assumption (vi) of Definition \ref{Def_admissible}.

Next we check \eqref{help47}.
If \ $x=(x_1,x_2)\in\RR_+\times\RR^d$ \ and \ $\xi\in\RR_+\times\RR^d$ \ with \ $\Vert\xi_2\Vert\leq 1$, \ then, by
 \eqref{2help15} with \ $\theta=1$, \
 \begin{align*}
  \vert f(x+\xi) - f(x) - \langle f'_{(2)}(x),\xi_2\rangle\vert
    & = \vert f(x+\xi) - f(x) - \langle f'(x),\xi\rangle + \langle f'_1(x),\xi_1\rangle \vert \\
    & \leq \frac{1}{2}\vert\langle f''(x+\tau\xi)\xi,\xi\rangle\vert + \vert\langle f'_1(x),\xi_1\rangle\vert\\
    & \leq \frac{1}{2}\Vert f''\Vert_\infty \Vert\xi\Vert^2
           + \Vert f_1''\Vert_\infty \xi_1,
 \end{align*}
 where \ $\tau=\tau(x,\xi)\in[0,1]$.
\ If \ $x=(x_1,x_2)\in\RR_+\times\RR^d$ \ and \ $\xi\in\RR_+\times\RR^d$ \ with \ $\Vert\xi_2\Vert>1$, \ then, by
 \eqref{help15} with \ $\theta=1$, \
 \begin{align*}
  \vert f(x+\xi) - f(x) & - \langle f'_{(2)}(x),\xi_2\rangle\vert
    = \vert \langle f'(x+\tau\xi),\xi\rangle  - \langle f'_{(2)}(x),\xi_2\rangle \vert \\
   & = \vert \langle f'_1(x+\tau\xi),\xi_1\rangle + \langle f'_{(2)}(x+\tau\xi),\xi_2\rangle - \langle f'_{(2)}(x),\xi_2\rangle \vert\\
   & \leq \Vert f_1''\Vert_\infty \xi_1 + \vert  \langle f'_{(2)}(x+\tau\xi) - f'_{(2)}(x),\xi_2\rangle \vert\\
   & \leq \Vert f_1''\Vert_\infty \xi_1 + 2 \Vert f_{(2)}'\Vert \Vert\xi_2\Vert,
 \end{align*}
 where \ $\tau=\tau(x,\xi)\in[0,1]$.
Using assumption (v) of Definition \ref{Def_admissible}, the finiteness of the integral in \eqref{help47}
 follows as for the integral in \eqref{help46}.

Having proved that the integrals in \eqref{affine_inf_gen} are well-defined, we check that
 under the conditions (v) and (vi) of Definition \ref{Def_admissible}, one can rewrite (2.12)
 in Duffie et al. \cite{DufFilSch} into the form \eqref{affine_inf_gen},
 by changing the \ $2$-nd, $\ldots$, $(1+d)$-th coordinates of \ $b\in\RR_+\times\RR^d$ \ and the first column of
 \ $\beta\in\RR^{(1+d)\times(1+d)}$, \ respectively.
More precisely, with the notations
 \begin{align*}
  \chi(\xi):=(\chi_1(\xi),\ldots,\chi_{1+d}(\xi))
   \qquad \text{and} \qquad
  \chi_{(2)}(\xi):=(\chi_2(\xi),\ldots,\chi_{1+d}(\xi))
 \end{align*}
 for \ $\xi\in\RR^{1+d}$, \ where
 \[
   \chi_k(\xi)
     := \begin{cases}
             (1\wedge \vert \xi_k\vert)\frac{\xi_k}{\vert \xi_k\vert} & \text{if \ $\xi_k\ne0$,}\\
             0 & \text{if \ $\xi_k=0$,}
        \end{cases}
      \qquad k=1,\ldots,1+d,
 \]
 for all \ $x=(x_1,x_2)\in\RR_+\times\RR^d$ \ and \ $f\in C^2_c(\RR_+\times\RR^d)$, \ we have
 \begin{align*}
   (\cA f)(x)
    &= \sum_{i,j=1}^{1+d} (a_{i,j} + \alpha_{i,j} x_1) f_{i,j}''(x)
        + \langle f'(x), b + \beta x \rangle\\
    &\phantom{=\;}  + \int_{\RR_+\times \RR^d}
         (f(x+\xi) - f(x) - \langle f'_{(2)}(x), \xi_2 \rangle)
         \, m(\dd\xi)\\
   &\phantom{=\;}
      + \int_{\RR_+\times \RR^d}
         (f(x+\xi) - f(x) - \langle f'(x), \xi \rangle) x_1
          \, \mu(\dd\xi)\\
   &= \sum_{i,j=1}^{1+d} (a_{i,j} + \alpha_{i,j} x_1) f_{i,j}''(x)
      + \langle f'(x), \widetilde b + \widetilde \beta x \rangle\\
   &\phantom{=\;} + \int_{\RR_+\times \RR^d}
         (f(x+\xi) - f(x) - \langle f'_{(2)}(x), \chi_{(2)}(\xi) \rangle)
         \, m(\dd\xi)\\
   &\phantom{=\;}
      + \int_{\RR_+\times \RR^d}
         (f(x+\xi) - f(x) - \langle f'(x), \chi(\xi) \rangle) x_1
          \, \mu(\dd\xi) ,
 \end{align*}
 where \ $\widetilde b = (\widetilde b_1, \widetilde b_2)\in\RR_+\times\RR^d$ \ and
 \ $\widetilde\beta = (\widetilde\beta_{i,j})_{i,j=1}^{1+d}\in\RR^{(1+d)\times (1+d)}$ \ with
 \begin{align*}
  &\widetilde b_1 := b_1,\\
  &\widetilde b_2 := b_2 + \int_{\RR_+\times\RR^d} (\chi_{(2)}(\xi) - \xi_2)\,m(\dd\xi),\\
  &(\widetilde\beta_{i,1})_{i=1}^{1+d}
      := (\beta_{i,1})_{i=1}^{1+d}
         + \int_{\RR_+\times\RR^d} (\chi(\xi) - \xi)\,\mu(\dd\xi),\\
  &(\widetilde\beta_{i,j})_{i=1,\ldots,1+d}^{j=2,\ldots,1+d}
      := (\beta_{i,1})_{i=1,\ldots,1+d}^{j=2,\ldots,1+d}.
 \end{align*}

Note also that there is another way for checking that the integrals in \eqref{affine_inf_gen} are well-defined
 under the conditions (v) and (vi) of Definition \ref{Def_admissible}.
Namely, using that the integrals in (2.12) in Duffie et al. \cite{DufFilSch} are well-defined, the assertion follows since
 \begin{align*}
  \left\Vert \int_{\RR_+\times\RR^d} (\chi_{(2)}(\xi) - \xi_2) \,m(\dd\xi) \right\Vert
    &= \left( \sum_{i=2}^{1+d} \left( \int_{\RR_+\times\RR^d} (\chi_i(\xi) - \xi_i)\,m(\dd\xi) \right)^2 \right)^{1/2} \\
  &\phantom{=\;} = \left( \sum_{i=2}^{1+d} \left( \int_{\RR_+\times\RR^d} \vert \xi_i\vert
               \mathbf 1_{\{\vert \xi_i\vert >1 \}} \,m(\dd\xi) \right)^2 \right)^{1/2}\\
  &\phantom{=\;} \leq \sqrt{d}\int_{\RR_+\times\RR^d} \Vert \xi\Vert \mathbf 1_{\{\Vert \xi\Vert >1 \}} \,m(\dd\xi)
    <\infty,
 \end{align*}
 and similarly
 \[
    \left\Vert \int_{\RR_+\times\RR^d} (\chi(\xi) - \xi) \,\mu(\dd\xi) \right\Vert <\infty.
 \]

\section{On consistency properties of the LSE and CLSE of~$(\theta,m)$}
\label{AppendixB}


Let us suppose that Condition (C) holds.
Using Slutsky's lemma and Theorems \ref{Thm3} and \ref{Thm4} we get
 \ $\widehat\theta_n^{\mathrm{LSE}}$ \ and \ $\widehat\theta_n^{\mathrm{CLSE}}$ \ converge stochastically to the parameter \ $\theta=0$ \
 as \ $n\to\infty$, \ respectively, and in what follows we show that \ $\widehat m_n^{\mathrm{LSE}}$ \ and
 \ $\widehat m_n^{\mathrm{CLSE}}$ \ do not converge stochastically to the parameter \ $m$ \ as \ $n\to\infty$, \ respectively.
For this it is enough to check that the weak limits of \ $\widehat m_n^{\mathrm{LSE}}$ \ and \ $\widehat m_n^{\mathrm{CLSE}}$ \
 given in Theorems \ref{Thm3} and \ref{Thm4} do not equal to \ $m$ \ almost surely, respectively.
Since the weak limits in question are the same, we can give a common proof.
First note that
 \begin{align*}
  &\frac{\cX_1 \int_0^1 \cX_t^2 \, \dd t
        - \int_0^1 \cX_t \, \dd t \int_0^1 \cX_t \, \dd \cX_t}
       {\int_0^1 \cX_t^2 \, \dd t - \left( \int_0^1 \cX_t \, \dd t \right)^2}
    - m \\
  & \qquad = \frac{(\cX_1 - m) \int_0^1 \cX_t^2 \, \dd t
          - \int_0^1 \cX_t \, \dd t \left(\int_0^1 \cX_t \, \dd \cX_t - m \int_0^1 \cX_t \, \dd t\right)}
        {\int_0^1 \cX_t^2 \, \dd t - \left( \int_0^1 \cX_t \, \dd t \right)^2}\\
  & \qquad = \frac{\int_0^1 \cX_t^2 \, \dd t \int_0^1 \dd [\cX_t - mt]
          - \int_0^1 \cX_t \, \dd t \int_0^1 \cX_t \, \dd [\cX_t - mt]}
         {\int_0^1 \cX_t^2 \, \dd t
          - \left( \int_0^1 \cX_t \, \dd t \right)^2} ,
 \end{align*}
 and hence
 \[
   \frac{\cX_1 \int_0^1 \cX_t^2 \, \dd t
         - \int_0^1 \cX_t \, \dd t \int_0^1 \cX_t \, \dd \cX_t}
        {\int_0^1 \cX_t^2 \, \dd t - \left( \int_0^1 \cX_t \, \dd t \right)^2}
   \ase m
 \]
 if and only if
 \[
   J := \int_0^1 \cX_t^2 \,\dd t \int_0^1 \dd [\cX_t - mt]
        - \int_0^1 \cX_t \, \dd t \int_0^1 \cX_t \, \dd [\cX_t - mt]
   \ase 0 ,
 \]
 where \ $\ase$ \ denotes equality almost surely.
Here \ $J$ \ can be written in the form
 \[
   J = \int_0^1 \cX_s^2 \left( \int_0^1 \, \dd [\cX_t - mt] \right) \dd s
       - \int_0^1 \cX_s
          \left( \int_0^1 \cX_t \, \dd [\cX_t - mt] \right) \dd s ,
 \]
 and hence \ $\EE(J)$ \ takes the following form
 \begin{align*}
    \int_0^1 \EE\left( \cX_s^2 \int_0^1 \, \dd [\cX_t - mt] \right) \dd s
     - \int_0^1 \EE\left( \cX_s \int_0^1 \cX_t \, \dd [\cX_t - mt] \right) \dd s .
 \end{align*}
Here
 \ $(\cX_t - mt)_{t \in \RR_+}
    = \bigl(\int_0^t \sqrt{\cY_u} \, \dd \cB_u\bigr)_{t \in \RR_+}$
 \ is a square integrable martingale (see the proof of Proposition \ref{Pro_moments}) and
 \ $\cX_s = ms + \int_0^s \sqrt{\cY_u} \, \dd \cB_u$ \ for \ $s \in \RR_+$,
 \ thus, for all \ $s \in [0,1]$, \ we have
 \begin{align*}
  &\EE\left( \cX_s^2 \int_0^1 \, \dd [\cX_t - mt]
            \, \Bigg| \, \cF_1^\cY \right) \\
  &= m^2 s^2
      \EE\left( \int_0^1 \sqrt{\cY_t} \, \dd \cB_t
                \, \Bigg| \, \cF_1^\cY \right)
       + 2 ms \EE\left( \int_0^s \sqrt{\cY_u} \, \dd \cB_u
                  \int_0^1 \sqrt{\cY_t} \, \dd \cB_t
                  \, \Bigg| \, \cF_1^\cY \right)\\
  &\phantom{=\;} + \EE\left( \left( \int_0^s \sqrt{\cY_u} \, \dd \cB_u \right)^2
                  \int_0^1 \sqrt{\cY_t} \, \dd \cB_t
                  \, \Bigg| \, \cF_1^\cY \right)\\
  & = 2 ms \int_0^s \cY_u \, \dd u,
  \end{align*}
 where \ $\cF_1^\cY$ \ denotes the \ $\sigma$-algebra generated by \ $(\cY_u)_{u\in[0,1]}$.
\ For the last equality above, we used that conditionally on the \ $\sigma$-algebra \ $\cF_1^\cY$,
 \ the stochastic process  \ $\bigl(\int_0^t \sqrt{\cY_u} \, \dd \cB_u\bigr)_{t \in [0,1]}$ \ is a Gauss process
 with mean function identically \ $0$ \ and with covariance function
 \begin{align*}
   \EE\left(\int_0^s \sqrt{\cY_u} \, \dd \cB_u\int_0^t \sqrt{\cY_u} \, \dd \cB_u
             \, \Bigg| \, \cF_1^\cY \right)
        = \int_0^{s\wedge t} \cY_u \, \dd u,
        \qquad  s,t\in[0,1],
 \end{align*}
 (for the mean and covariation function, see Karatzas and Shreve \cite[formulas (3.2.21) and (3.2.23)]{KarShr}),
 and we also used that the third moment of a centered normally distributed random variable is \ $0$.
\ Similarly, for all \ $s \in [0,1]$, \ we have
 \begin{align*}
   &\EE\left( \cX_s \int_0^1 \cX_t \, \dd [\cX_t - mt]
            \, \Bigg| \, \cF_1^\cY \right) \\
    &= ms
      \EE\left( \int_0^1 \cX_t \sqrt{\cY_t} \, \dd \cB_t
                \, \Bigg| \, \cF_1^\cY \right)
      + \EE\left( \int_0^s \sqrt{\cY_u} \, \dd \cB_u
                  \int_0^1 \cX_t \sqrt{\cY_t} \, \dd \cB_t
                  \, \Bigg| \, \cF_1^\cY \right)\\
   &\; = \EE\left( \int_0^s \cX_u \cY_u \, \dd u \, \Bigg| \, \cF_1^\cY \right)
           = \int_0^s \EE\left( \cX_u \cY_u \mid \cF_1^\cY \right) \, \dd u
           = m \int_0^s u \cY_u \, \dd u .
 \end{align*}
Thus
 \begin{align*}
  &\EE\left( \cX_s^2 \int_0^1 \, \dd [\cX_t - mt] \right)
  = 2 ms \int_0^s \EE(\cY_u) \, \dd u , \\
 &\EE\left( \cX_s \int_0^1 \cX_t \, \dd [\cX_t - mt] \right)
  = m \int_0^s u \EE(\cY_u) \, \dd u .
 \end{align*}
\allowdisplaybreaks%
Consequently,
 \begin{align*}
  \EE(J) &= 2m \int_0^1 s \left( \int_0^s \EE(\cY_u) \, \dd u \right) \dd s
            - m \int_0^1
                 \left( \int_0^s u \EE(\cY_u) \, \dd u \right) \dd s \\
         &= 2m \int_0^1 \left( \int_u^1 s \, \dd s \right) \EE(\cY_u) \, \dd s
            - m \int_0^1
                 \left( \int_u^1 \, \dd s \right) u \EE(\cY_u) \, \dd u \\
         &= m \int_0^1 (1-u) \EE(\cY_u) \, \dd u
          = ma \int_0^1 (1-u) u \, \dd u.
 \end{align*}
Hence if \ $m\ne 0$, \ then \ $\EE(J) \ne 0$, \ which clearly yields that \ $J \ase 0$ \ is impossible.

If \ $m=0$, \ then \ $\EE(J) = 0$, \ and hence for proving \ $\PP(J=0)<1$, \ it is enough to show that \ $\EE(J^2) > 0$.
\ Now \ $J$ \ can be written in the form \ $J = J_1 - J_2$ \ with
 \[
   J_1 := \int_0^1 \cX_s^2 \, \dd s \int_0^1 \, \dd \cX_t ,
   \qquad
   J_2 := \int_0^1 \cX_s \, \dd s \int_0^1 \cX_t \, \dd \cX_t .
 \]
Clearly, \ $\EE(J^2) = \EE(J_1^2) - 2\EE(J_1J_2) + \EE(J_2^2)$.
\ Here \ $J_1 J_2$ \ can be written in the form
 \[
   J_1 J_2 = \int_0^1
             \int_0^1
              \cX_{s_1}^2 \cX_{s_2}
              \left( \int_0^1 \, \dd \cX_t \right)
              \left( \int_0^1 \cX_t \, \dd \cX_t \right) \dd s_1 \, \dd s_2 ,
 \]
 hence
 \[
   \EE(J_1J_2)
   = \int_0^1
      \int_0^1
       \EE\left( \cX_{s_1}^2 \cX_{s_2}
                 \left( \int_0^1 \, \dd \cX_t \right)
                 \left( \int_0^1 \cX_t \, \dd \cX_t \right)
          \right)
       \dd s_1 \, \dd s_2 .
 \]
If \ $s_1, s_2 \in [0,1]$, \ then, by \eqref{help49},
 \begin{align*}
  &\EE\left( \cX_{s_1}^2 \cX_{s_2}
             \left( \int_0^1 \, \dd \cX_t \right)
             \left( \int_0^1 \cX_t \, \dd \cX_t \right)
             \, \bigg| \, \cF_1^\cY \right)\\
  &\qquad = \frac{1}{2}
          \EE\left( \cX_{s_1}^2 \cX_{s_2} \cX_1
             \left( \cX_1^2  - \int_0^1 \cY_s \, \dd s \right)
             \, \bigg| \, \cF_1^\cY \right)\\
   &\qquad = \frac{1}{2}
          \EE\left( \cX_{s_1}^2 \cX_{s_2} \cX_1^3 \, \bigg| \, \cF_1^\cY \right)
    - \frac{1}{2}
          \EE\left( \cX_{s_1}^2 \cX_{s_2} \cX_1 \int_0^1\cY_s\,\dd s \, \bigg| \, \cF_1^\cY \right)\\
  &\qquad = \frac{1}{2}
          \EE\left( \cX_{s_1}^2 \cX_{s_2} \cX_1^3 \, \bigg| \, \cF_1^\cY \right)
            - \frac{1}{2} \left(\int_0^1\cY_s\,\dd s\right)
          \EE\left( \cX_{s_1}^2 \cX_{s_2} \cX_1 \, \bigg| \, \cF_1^\cY \right).
 \end{align*}
Similarly to the proof of Proposition \ref{Pro_affine} one can check that conditionally of \ $\cF_1^{\cY}$, \ the process
 \ $(\cX_t)_{t\in[0,1]}$ \ is a centered Gauss process (especially having independent increments).
Hence if \ $s_1, s_2 \in [0,1]$ \ with \ $s_1 \leq s_2$, \ then
 \begin{align*}
       &\EE\left( \cX_{s_1}^2 \cX_{s_2} \cX_1^3 \, \bigg| \, \cF_1^\cY \right)\\
       & = \EE\left( \cX_{s_1}^6 \, \bigg| \, \cF_1^\cY \right)
           + 4\EE\left( \cX_{s_1}^3 (\cX_{s_2} - \cX_{s_1})^3 \, \bigg| \, \cF_1^\cY \right)
           + \EE\left( \cX_{s_1}^3 (\cX_1 - \cX_{s_2})^3 \, \bigg| \, \cF_1^\cY \right)\\
       &\phantom{=\;}  + 4\EE\left( \cX_{s_1}^5 (\cX_{s_2} - \cX_{s_1}) \, \bigg| \, \cF_1^\cY \right)
                       + 3\EE\left( \cX_{s_1}^5 (\cX_1 - \cX_{s_2}) \, \bigg| \, \cF_1^\cY \right)\\
       &\phantom{=\;}  + 6\EE\left( \cX_{s_1}^4 (\cX_{s_2} - \cX_{s_1})^2 \, \bigg| \, \cF_1^\cY \right)
                       + 3\EE\left( \cX_{s_1}^4 (\cX_1 - \cX_{s_2})^2 \, \bigg| \, \cF_1^\cY \right)\\
       &\phantom{=\;}  + 6\EE\left( \cX_{s_1}^3 (\cX_{s_2} - \cX_{s_1})(\cX_1 - \cX_{s_2})^2 \, \bigg| \, \cF_1^\cY \right)\\
        &\phantom{=\;}+ 3\EE\left( \cX_{s_1}^3 (\cX_{s_2} - \cX_{s_1})^2(\cX_1 - \cX_{s_2}) \, \bigg| \, \cF_1^\cY \right)\\
        &\phantom{=\;} + 9\EE\left( \cX_{s_1}^4 (\cX_{s_2} - \cX_{s_1})(\cX_1 - \cX_{s_2}) \, \bigg| \, \cF_1^\cY \right)
                       + \EE\left( \cX_{s_1}^2 (\cX_{s_2} - \cX_{s_1})^4 \, \bigg| \, \cF_1^\cY \right)\\
       &\phantom{=\;}  + \EE\left( \cX_{s_1}^2 (\cX_{s_2} - \cX_{s_1})(\cX_1 - \cX_{s_2})^3 \, \bigg| \, \cF_1^\cY \right)\\
       &\phantom{=\;}  + 3\EE\left( \cX_{s_1}^2 (\cX_{s_2} - \cX_{s_1})^2(\cX_1 - \cX_{s_2})^2 \, \bigg| \, \cF_1^\cY \right)\\
       &\phantom{=\;}  + 3\EE\left( \cX_{s_1}^2 (\cX_{s_2} - \cX_{s_1})^3(\cX_1 - \cX_{s_2}) \, \bigg| \, \cF_1^\cY \right)\\
       &\phantom{=\;}  + 6\EE\left( \cX_{s_1}^3 (\cX_{s_2} - \cX_{s_1})^2(\cX_1 - \cX_{s_2}) \, \bigg| \, \cF_1^\cY \right).
 \end{align*}
Using that \ $\cX$ \ has independent increments and that the odd moments of a centered normally distributed random
 variable are \ $0$, \ if \ $s_1, s_2 \in [0,1]$ \ with \ $s_1 \leq s_2$, \ then
 \begin{align*}
    &\EE\left( \cX_{s_1}^2 \cX_{s_2} \cX_1^3 \, \bigg| \, \cF_1^\cY \right)\\
    &= \EE\left( \cX_{s_1}^6 \, \bigg| \, \cF_1^\cY \right)
        + 6\EE\left( \cX_{s_1}^4 \, \bigg| \, \cF_1^\cY \right) \EE\left((\cX_{s_2} - \cX_{s_1})^2 \, \bigg| \, \cF_1^\cY \right) \\
    &\phantom{=\;} + 3\EE\left( \cX_{s_1}^4 \, \bigg| \, \cF_1^\cY \right) \EE\left( (\cX_1 - \cX_{s_2})^2 \, \bigg| \, \cF_1^\cY \right)\\
    &\phantom{=\;} + \EE\left( \cX_{s_1}^2 \, \bigg| \, \cF_1^\cY \right) \EE\left((\cX_{s_2} - \cX_{s_1})^4 \, \bigg| \, \cF_1^\cY \right)\\
    &\phantom{=\;} + 3\EE\left( \cX_{s_1}^2 \, \bigg| \, \cF_1^\cY \right)
          \EE\left((\cX_{s_2} - \cX_{s_1})^2  \, \bigg| \, \cF_1^\cY \right)
          \EE\left((\cX_1 - \cX_{s_2})^2 \, \bigg| \, \cF_1^\cY \right)\\
    &= 15\left(\int_0^{s_1}\cY_u\,\dd u\right)^3
        + 18\left(\int_0^{s_1}\cY_u\,\dd u\right)^2 \left(\int_{s_1}^{s_2}\cY_u\,\dd u\right)\\
   &\phantom{=\;} + 9\left(\int_0^{s_1}\cY_u\,\dd u\right)^2 \left(\int_{s_2}^{1}\cY_u\,\dd u\right)
        + 3\left(\int_0^{s_1}\cY_u\,\dd u\right) \left(\int_{s_1}^{s_2}\cY_u\,\dd u\right)^2\\
     &\phantom{=\;}
        + 3\left(\int_0^{s_1}\cY_u\,\dd u\right) \left(\int_{s_1}^{s_2}\cY_u\,\dd u\right)
         \left(\int_{s_2}^1\cY_u\,\dd u\right).
 \end{align*}
One can also check that if \ $s_1, s_2 \in [0,1]$ \ with \ $s_1 \leq s_2$, \ then
 \begin{align*}
   \EE\left( \cX_{s_1}^2 \cX_{s_2} \cX_1 \, \bigg| \, \cF_1^\cY \right)
      & = \EE\left( \cX_{s_1}^4 \, \bigg| \, \cF_1^\cY \right)
        + \EE\left( \cX_{s_1}^2 (\cX_{s_2} - \cX_{s_1})^2 \, \bigg| \, \cF_1^\cY \right)\\
     & = 3\left(\int_0^{s_1}\cY_u\,\dd u\right)^2 + \left(\int_0^{s_1}\cY_u\,\dd u\right)\left(\int_{s_1}^{s_2}\cY_u\,\dd u\right) .
 \end{align*}
Hence if \ $s_1, s_2 \in [0,1]$ \ with \ $s_1 \leq s_2$, \ then
 \begin{align*}
  &\EE\left( \cX_{s_1}^2 \cX_{s_2}
             \left( \int_0^1 \, \dd \cX_t \right)
             \left( \int_0^1 \cX_t \, \dd \cX_t \right)
             \, \bigg| \, \cF_1^\cY \right) \\
   &\qquad = 6\left(\int_0^{s_1}\cY_u\,\dd u\right)^3
               + 7 \left(\int_0^{s_1}\cY_u\,\dd u\right)^2 \left(\int_{s_1}^{s_2}\cY_u\,\dd u\right)\\
   &\phantom{=\;}\qquad  + 3\left(\int_0^{s_1}\cY_u\,\dd u\right)^2 \left(\int_{s_2}^1\cY_u\,\dd u\right)
                         + \left(\int_0^{s_1}\cY_u\,\dd u\right) \left(\int_{s_1}^{s_2}\cY_u\,\dd u\right)^2\\
   &\phantom{=\;}\qquad  + \left(\int_0^{s_1}\cY_u\,\dd u\right)\left(\int_{s_1}^{s_2}\cY_u\,\dd u\right)
                           \left(\int_{s_2}^1\cY_u\,\dd u\right).
 \end{align*}
Similar expression hold in case of \ $s_1, s_2 \in [0,1]$ \ with \ $s_1 \geq s_2$, \ we have to change  \ $s_1$ \ by \ $s_2$.

Moreover, \ $J_1^2$ \ can be written in the form
 \[
   J_1^2 = \int_0^1
           \int_0^1
            \cX_{s_1}^2 \cX_{s_2}^2
            \left( \int_0^1 \, \dd \cX_t \right)^2 \dd s_1 \, \dd s_2 ,
 \]
 hence
 \[
   \EE(J_1^2)
   = \int_0^1
      \int_0^1
       \EE\left( \cX_{s_1}^2 \cX_{s_2}^2
                 \left( \int_0^1 \, \dd \cX_t \right)^2 \right)
       \dd s_1 \, \dd s_2 .
 \]
If \ $s_1, s_2 \in [0,1]$ \ with \ $s_1 \leq s_2$, \ then
 \begin{align*}
  &\EE\left( \cX_{s_1}^2 \cX_{s_2}^2
             \left( \int_0^1 \, \dd \cX_t \right)^2
             \, \bigg| \, \cF_1^\cY \right)\\
& = \EE\left( \cX_{s_1}^2
                \cX_{s_2}^2
                \cX_1^2
                \, \big| \, \cF_1^\cY \right)\\
   &= \EE\left( \cX_{s_1}^6
                \, \big| \, \cF_1^\cY \right)
      + \EE\left( \cX_{s_1}^4
                  \, \big| \, \cF_1^\cY \right)
        \EE\left( \left( \cX_1 - \cX_{s_2} \right)^2
                  \, \big| \, \cF_1^\cY \right) \\
   &\quad\,
      + 6 \EE\left( \cX_{s_1}^4
                    \, \big| \, \cF_1^\cY \right)
        \EE\left( \left( \cX_{s_2} - \cX_{s_1} \right)^2
                  \, \big| \, \cF_1^\cY \right) \\
   &\quad\,
      + \EE\left( \cX_{s_1}^2
                  \, \big| \, \cF_1^\cY \right)
        \EE\left( \left( \cX_{s_2} - \cX_{s_1} \right)^4
                  \, \big| \, \cF_1^\cY \right) \\
   &\quad\,
      + \EE\left( \cX_{s_1}^2
                  \, \big| \, \cF_1^\cY \right)
        \EE\left( \left( \cX_{s_2} - \cX_{s_1} \right)^2
                  \, \big| \, \cF_1^\cY \right)
        \EE\left( \left( \cX_1 - \cX_{s_2} \right)^2
                  \, \big| \, \cF_1^\cY \right)\\
  &= 15 \left(\int_0^{s_1}\cY_u\,\dd u\right)^3
     + 3 \left(\int_0^{s_1}\cY_u\,\dd u\right)^2 \left(\int_{s_2}^1\cY_u\,\dd u\right)\\
  &\phantom{=\;} + 18 \left(\int_0^{s_1}\cY_u\,\dd u\right)^2 \left(\int_{s_1}^{s_2}\cY_u\,\dd u\right)
                 + 3 \left(\int_0^{s_1}\cY_u\,\dd u\right) \left(\int_{s_1}^{s_2}\cY_u\,\dd u\right)^2\\
  &\phantom{=\;} + \left(\int_0^{s_1}\cY_u\,\dd u\right) \left(\int_{s_1}^{s_2}\cY_u\,\dd u\right)
       \left(\int_{s_2}^1\cY_u\,\dd u\right).
 \end{align*}
Similar expression hold in case of \ $s_1, s_2 \in [0,1]$ \ with \ $s_1 \geq s_2$, \ we have to change  \ $s_1$ \ by \ $s_2$.

Furthermore, using \eqref{help49},
 \ $J_2^2$ \ can be written in the form
 \[
   J_2^2 = \frac{1}{4}\int_0^1
           \int_0^1
            \cX_{s_1}\cX_{s_2}
            \left( \cX_1^2 - \int_0^1\cY_u\,\dd u \right)^2 \dd s_1 \, \dd s_2 ,
 \]
 hence
 \[
   \EE(J_2^2)=
     \frac{1}{4}\int_0^1 \int_0^1
            \EE\left(\cX_{s_1}\cX_{s_2}
                     \left( \cX_1^2 - \int_0^1\cY_u\,\dd u \right)^2\right) \dd s_1 \, \dd s_2.
 \]
Here if \ $s_1, s_2 \in [0,1]$ \ with \ $s_1 \leq s_2$, \ then we have
 \begin{align*}
  &\EE\left(\cX_{s_1}\cX_{s_2}
              \left( \cX_1^2 - \int_0^1\cY_u\,\dd u \right)^2
                 \, \bigg| \, \cF_1^\cY \right)  \\
  & = \EE\left(\cX_{s_1}\cX_{s_2}\cX_1^4 \, \bigg| \, \cF_1^\cY \right)
     - 2\left(\int_0^1\cY_u\,\dd u\right) \EE\left(\cX_{s_1}\cX_{s_2}\cX_1^2 \, \bigg| \, \cF_1^\cY \right)  \\
  &\phantom{=\;}
     + \left(\int_0^1\cY_u\,\dd u\right)^2 \EE\left(\cX_{s_1}\cX_{s_2} \, \bigg| \, \cF_1^\cY \right),
 \end{align*}
 where, using the arguments as above, one can check that
 \begin{align*}
  &\EE\left(\cX_{s_1}\cX_{s_2}\cX_1^4 \, \bigg| \, \cF_1^\cY \right)\\
  &= \EE\left( \cX_{s_1}^6 \, \bigg| \, \cF_1^\cY \right)
       + 5 \EE\left( \cX_{s_1}^2 (\cX_{s_2} - \cX_{s_1})^4 \, \bigg| \, \cF_1^\cY \right)\\
  &\phantom{=\;} + \EE\left( \cX_{s_1}^2 (\cX_1 - \cX_{s_2})^4 \, \bigg| \, \cF_1^\cY \right)
                  + 10 \EE\left( \cX_{s_1}^4 (\cX_{s_2} - \cX_{s_1})^2\, \bigg| \, \cF_1^\cY \right)\\
  &\phantom{=\;}  + 6\EE\left( \cX_{s_1}^4 (\cX_1 - \cX_{s_2})^2 \, \bigg| \, \cF_1^\cY \right)
                  + 18 \EE\left( \cX_{s_1}^2 (\cX_{s_2} -
                  \cX_{s_1})^2(\cX_1 - \cX_{s_2})^2 \, \bigg| \, \cF_1^\cY \right)\\
  &= 15\left(\int_0^{s_1}\cY_u\,\dd u\right)^3
       + 15 \left(\int_0^{s_1}\cY_u\,\dd u\right) \left(\int_{s_1}^{s_2}\cY_u\,\dd u\right)^2\\
  &\phantom{=\;} + 3\left(\int_0^{s_1}\cY_u\,\dd u\right) \left(\int_{s_2}^1\cY_u\,\dd u\right)^2
                 + 30 \left(\int_0^{s_1}\cY_u\,\dd u\right)^2 \left(\int_{s_1}^{s_2}\cY_u\,\dd u\right)\\
  &\phantom{=\;} + 18 \left(\int_0^{s_1}\cY_u\,\dd u\right)^2 \left(\int_{s_2}^1\cY_u\,\dd u\right)\\
  &\phantom{=\;}
    + 18 \left(\int_0^{s_1}\cY_u\,\dd u\right) \left(\int_{s_1}^{s_2}\cY_u\,\dd u\right)\left(\int_{s_2}^1\cY_u\,\dd u\right),
 \end{align*}
 and
 \begin{align*}
  \EE\left(\cX_{s_1}\cX_{s_2}\cX_1^2 \, \bigg| \, \cF_1^\cY \right)
  &= \EE\left( \cX_{s_1}^4 \, \bigg| \, \cF_1^\cY \right)
       + 3 \EE\left( \cX_{s_1}^2 (\cX_{s_2} - \cX_{s_1})^2 \, \bigg| \, \cF_1^\cY \right)\\
  &\phantom{=\;} + \EE\left( \cX_{s_1}^2 (\cX_1 - \cX_{s_2})^2 \, \bigg| \, \cF_1^\cY \right)\\
  &= 3\left(\int_0^{s_1}\cY_u\,\dd u\right)^2
     + 3 \left(\int_0^{s_1}\cY_u\,\dd u\right) \left(\int_{s_1}^{s_2}\cY_u\,\dd u\right)\\
  &\phantom{=\;} + \left(\int_0^{s_1}\cY_u\,\dd u\right) \left(\int_{s_2}^1\cY_u\,\dd u\right),
 \end{align*}
 and
 \begin{align*}
  \EE\left(\cX_{s_1}\cX_{s_2} \, \bigg| \, \cF_1^\cY \right)
  = \EE\left( \cX_{s_1}^2 \, \bigg| \, \cF_1^\cY \right)
  =  \int_0^{s_1}\cY_u\,\dd u.
 \end{align*}
Hence, by an easy calculation, if \ $s_1, s_2 \in [0,1]$ \ with \ $s_1 \leq s_2$, \ then
 \begin{align*}
   \EE&\left(\cX_{s_1}\cX_{s_2}
              \left( \cX_1^2 - \int_0^1\cY_u\,\dd u \right)^2
                 \, \bigg| \, \cF_1^\cY \right)  \\
      &  = 10 \left(\int_0^{s_1}\cY_u\,\dd u\right)^3
         + 10 \left(\int_0^{s_1}\cY_u\,\dd u\right) \left(\int_{s_1}^{s_2}\cY_u\,\dd u\right)^2\\
     &\phantom{=\;} + 20 \left(\int_0^{s_1}\cY_u\,\dd u\right)^2 \left(\int_{s_1}^{s_2}\cY_u\,\dd u\right)\\
     &\phantom{=\;} + 12 \left(\int_0^{s_1}\cY_u\,\dd u\right)^2 \left(\int_{s_2}^1\cY_u\,\dd u\right)\\
     &\phantom{=\;} + 2 \left(\int_0^{s_1}\cY_u\,\dd u\right)\left(\int_{s_2}^1\cY_u\,\dd u\right)^2\\
     &\phantom{=\;} + 12 \left(\int_0^{s_1}\cY_u\,\dd u\right)\left(\int_{s_1}^{s_2}\cY_u\,\dd u\right)
            \left(\int_{s_2}^1\cY_u\,\dd u\right).
 \end{align*}
Similar expression hold in case of \ $s_1, s_2 \in [0,1]$ \ with \ $s_1 \geq s_2$, \ we have to change  \ $s_1$ \ by \ $s_2$.

Hence if \ $s_1, s_2 \in [0,1]$ \ with \ $s_1 \leq s_2$, \ then we have
 \begin{align*}
   &\EE\Bigg( \cX_{s_1}^2 \cX_{s_2}^2
             \left( \int_0^1 \, \dd \cX_t \right)^2 - 2\cX_{s_1}^2 \cX_{s_2}
             \left( \int_0^1 \, \dd \cX_t \right)
             \left( \int_0^1 \cX_t \, \dd \cX_t\right)\\
   &\phantom{\EE\Bigg(\,}
             +\frac{1}{4}\cX_{s_1}\cX_{s_2}\left(\cX_1^2 - \int_0^1 \cY_u\,\dd u\right)^2
             \, \bigg| \, \cF_1^\cY
             \Bigg)\\
   &= \frac{11}{2}\left(\int_0^{s_1}\cY_u\,\dd u\right)^3
               + 9 \left(\int_0^{s_1}\cY_u\,\dd u\right)^2 \left(\int_{s_1}^{s_2}\cY_u\,\dd u\right)\\
   &\phantom{=\;} + \frac{7}{2} \left(\int_0^{s_1}\cY_u\,\dd u\right) \left(\int_{s_1}^{s_2}\cY_u\,\dd u\right)^2
                 + \frac{1}{2}\left(\int_0^{s_1}\cY_u\,\dd u\right) \left(\int_{s_2}^1\cY_u\,\dd u\right)^2\\
   &\phantom{=\;} + 2\left(\int_0^{s_1}\cY_u\,\dd u\right) \left(\int_{s_1}^{s_2}\cY_u\,\dd u\right)
                   \left(\int_{s_2}^1\cY_u\,\dd u\right),
 \end{align*}
 and a similar expression hold in case of \ $s_1, s_2 \in [0,1]$ \ with \ $s_1 \geq s_2$,
 \ we have to change  \ $s_1$ \ by \ $s_2$.
\ Then \ $\EE(J^2) = \EE(J_1^2) - 2\EE(J_1J_2) + \EE(J_2^2)$ \ takes the form
 \begin{align*}
  & \int_0^1\int_0^1
             \EE\Bigg(
                \frac{11}{2}\left(\int_0^{s_1\wedge s_2}\cY_u\,\dd u\right)^3
             + 9 \left(\int_0^{s_1\wedge s_2}\cY_u\,\dd u\right)^2 \left(\int_{s_1\wedge s_2}^{s_1\vee s_2}\cY_u\,\dd u\right)\\
  &\phantom{= }
      + \frac{7}{2} \left(\int_0^{s_1\wedge s_2}\cY_u\,\dd u\right)\!\! \left(\int_{s_1\wedge s_2}^{s_1\vee s_2}\cY_u\,\dd u\right)^2
      + \frac{1}{2}\left(\int_0^{s_1\wedge s_2}\cY_u\,\dd u\right) \!\!\left(\int_{s_1\vee s_2}^1\cY_u\,\dd u\right)^2\\
  &\phantom{= }
               + 2\left(\int_0^{s_1\wedge s_2}\cY_u\,\dd u\right) \left(\int_{s_1\wedge s_2}^{s_1\vee s_2}\cY_u\,\dd u\right)
                   \left(\int_{s_1\vee s_2}^1\cY_u\,\dd u\right)
               \Bigg)\,\dd s_1\dd s_2 >0,
 \end{align*}
 where for the last inequality we used that
 \[
   \EE\left( \left(\int_0^{s_1\wedge s_2}\cY_u\,\dd u\right)^i \left(\int_{s_1\wedge s_2}^{s_1\vee s_2}\cY_u\,\dd u\right)^j
              \left(\int_{s_1\vee s_2}^1\cY_u\,\dd u\right)^k \right)
             >0
 \]
 for \ $i,j,k\in\{0,1,2,3\}$, \ which follows by \ $a\in\RR_{++}$ \ and that \ $\PP(Y_t\geq 0 \;\;\text{for all \ $t\in\RR_+$})=1$.
\ Consequently, we conclude \ $\EE(J^2) > 0$, \ which clearly yields that \ $J \ase 0$ \ is impossible.

\indent
   {\sc M\'aty\'as Barczy},
   Faculty of Informatics, University of Debrecen,
   Pf.~12, H--4010 Debrecen, Hungary.
   Tel.: +36-52-512900, Fax: +36-52-512996,
   {\sl e--mail: barczy.matyas@inf.unideb.hu}

\indent
    {\sc Leif D\"oring},
    Laboratoire de Probabilit\'es et Mod\'eles Al\'eatoires,
    Universit\'e Paris 6, 4 place Jussieu, 75252 Paris, France.
    Tel.: +33-0-1-44275319, Fax: +33-0-1-44277223,\\
    {\sl e--mail: leif.doering@googlemail.com}

\indent
     {\sc Zenghu Li},
     School of Mathematical Sciences, Beijing Normal University,
     Beijing 100875,  People's Republic of China.
     Tel.: +86-10-58802900, Fax: +86-10-58808202,
     {\sl e--mail: lizh@bnu.edu.cn}

\indent
      {\sc Gyula Pap},
      Bolyai Institute, University of Szeged,
      Aradi v\'ertan\'uk tere 1, H--6720 Szeged, Hungary.
      Tel.: +36-62-544033, Fax: +36-62-544548,
      {\sl e--mail: papgy@math.u-szeged.hu}

\end{document}